\documentclass[11pt]{amsart} 
\title{\texttt{mcode.sty} Demo}
\usepackage[utf8]{inputenc}
\usepackage{amsmath}
\usepackage{epsfig}
\usepackage{siunitx} 
\usepackage{color,graphics}
\usepackage{booktabs}
\usepackage{multirow}
\usepackage{epstopdf}
\usepackage[all]{xy}
\usepackage{graphicx}
\usepackage{xcolor}
\usepackage{pagecolor}
\pagecolor{white}
\usepackage{amsmath, amssymb}
\numberwithin{equation}{section}
\usepackage[vlined, ruled]{algorithm2e} 
\usepackage{fixmath, bm, amssymb, amsfonts, latexsym}
\usepackage{xcolor}
\usepackage{caption}
\usepackage{subcaption}
\usepackage{ragged2e}

 \usepackage[a4paper,
            width=18.6cm,
            height=26.3cm,
            centering]{geometry}
\usepackage{hyperref}
\usepackage{xcolor}
\hypersetup{
  colorlinks   = true, 
  urlcolor     = red, 
  linkcolor    = blue, 
  citecolor   = blue 
}

\newtheorem{theorem}{Theorem}[section]

\newtheorem{definition}[theorem]{Definition}

\newtheorem{remark}[theorem]{Remark}

\newcommand{\baa}{\begin{eqnarray*}}
\newcommand{\eaa}{\end{eqnarray*}}
\newcommand{\ba}{\begin{equation}}
\newcommand{\ea}{\end{equation}}

\pagecolor{white}

\usepackage{lineno, blindtext}
\usepackage{float}

\usepackage{amsmath}

\author{Sudipta Sahu$^{\S}$, Rathan Samala$^{\S}$*}
\title{Integral Invariants of Vasudeva Murthy's Relaxation Systems: Analysis and Numerical Validation}
\thanks{ 
$^{\S}$Department of Humanities and Sciences, Indian Institute of Petroleum and Energy-Visakhapatnam,
India-531035({sudipta.sahu@iipe.ac.in, rathans.math@iipe.ac.in})
\newline
$^*$  Corresponding author: Rathan Samala}
\date{\today}
\begin{document}
\maketitle
\begin{abstract}
\vspace{-1.0cm}
Vasudeva Murthy’s relaxation approach [A.S. Vasudeva Murthy, J. Comput. Appl. Math., 203(2), pp. 437-443, 2007], originally proposed for the Jin-Xin relaxation model, provides an alternative formulation with invariant properties that is consistent and retains the semilinear structure incomparison to the standard one.  In this work,  Vasudeva Murthy’s relaxation approach for various relaxation systems are proposed such as the shallow-water equations, the Broadwell model, the Euler equations with heat transfer and two-dimensional Jin-Xin model. For  proposed relaxation models, the associated integral invariants are rigorously established at the theoretical level. The main advantage of the integral invariant is that it provides a conserved quantity for the relaxation system by incorporating the coupled contributions of the solution variables in vector form.   To validate the analytical results, numerical simulations are carried out for each model using three second-order numerical schemes: CS-EBT2, a semi-implicit second-order central finite-volume scheme for hyperbolic systems with relaxation source terms [S. Sahu, E. Macca, and R. Samala, J. Comput. Phys., 563, 115100, 2026]; UCS2, a finite-volume central relaxation-type scheme [S. F. Liotta, V. Romano, and G. Russo, SIAM J. Numer. Anal., 38(4), 1337-1356, 2000]; and IMEX-RK2, a second-order Implicit-Explicit Runge-Kutta scheme [Pareschi and Russo, J. Sci. Comput., 25, 129–155, 2005]. Numerical results, compared to exact or finely resolved reference solutions, confirm that the models preserve integral invariants, remain stable under CFL restrictions, and exhibit robust and accurate behavior across all benchmark systems tested.
\end{abstract}
\bigskip
\noindent 2000 AMS(MOS) Classification:
41A10, 65M06
\medskip
\noindent
Keywords: Hyperbolic systems of balance laws, Stiff source terms, Relaxation model, Implicit-Explicit scheme, Central scheme, Integral invariance.
\pagestyle{myheadings} \thispagestyle{plain} \markboth{Sudipta Sahu, Rathan Samala}{Numerical Scheme to Vasudeva Murthy's Alternative Relaxation Approach}
\section{Introduction}
Hyperbolic balance laws with stiff source terms provide a mathematical framework for describing physical phenomena in which convective transport and relaxation mechanisms coexist. Such systems arise in numerous applications, including compressible fluid dynamics, kinetic theory, and shallow-water flows \cite{LeVeque2002,Toro2009}. A general one-dimensional hyperbolic balance law can be written as
\begin{equation}
U_t + F(U)_x = S(U),
\label{eq:balance-law}
\end{equation}
where $U$ denotes the vector of conserved variables, $F(U)$ is the flux function, and $S(U)$ represents the source term. The presence of stiff source terms introduces multiple time scales into the governing equations, making both their mathematical analysis and numerical approximation considerably more challenging. Consequently, numerical methods must simultaneously maintain stability and accurately capture the interaction between convective transport and relaxation processes.

Although the present study primarily focuses on one-dimensional problems, hyperbolic balance laws naturally admit multidimensional extensions \cite{KatsoulakisTzavaras1997, SW_friction}. In two spatial dimensions, the governing equations take the form
\begin{equation}
U_t + F_1(U)_x + F_2(U)_y = S(U),
\label{eq:2D-balance-law}
\end{equation}
where $F_1(U)$ and $F_2(U)$ denote the fluxes in the $x$- and $y$-directions, respectively. Multidimensional balance laws arise in a variety of practical applications, including shallow-water flows and gas dynamics, and therefore constitute an important extension of the one-dimensional setting.

The numerical approximation of hyperbolic balance laws has received considerable attention over the past several decades, leading to the development of numerous finite-volume methods. Among these, the central schemes introduced by Nessyahu and Tadmor \cite{NT} have become a popular class of Riemann-solver-free methods owing to their simplicity, robustness, and ability to achieve high-resolution approximations without characteristic decomposition or exact and approximate Riemann solvers. Building upon this framework, several high-resolution central schemes have subsequently been proposed for non-homogeneous hyperbolic balance laws \cite{Liotta,russo,Pareschi,SB}. In particular, the second-order CS-EBT2 \cite{2026} and UCS2 \cite{Liotta} schemes retain the simplicity of the original central formulation while providing accurate and essentially non-oscillatory solutions. Their Riemann-solver-free formulation makes them particularly attractive for problems involving stiff relaxation source terms. An alternative and widely used strategy for treating stiffness is provided by implicit-explicit (IMEX) Runge-Kutta methods \cite{ImexBosca,Boscarino-Filbet,PareschiRusso}. These methods discretize the convective flux explicitly while integrating the stiff source term implicitly, thereby achieving stable and efficient time integration for relaxation systems without sacrificing accuracy.

In parallel with advances in numerical discretization, considerable effort has been devoted to the analytical development of relaxation approximations for hyperbolic balance laws. A successful relaxation model should not only provide an accurate approximation of the original balance law but also preserve its essential analytical properties. In this direction, Vasudeva Murthy  in 2007 \cite{murthy2007alternative} proposed an alternative relaxation formulation for the linear Jin-Xin model. Compared with the classical Jin-Xin relaxation system \cite{Sjin_1995}, the proposed formulation preserves the semi linear structure while satisfying an integral invariance property, the standard sub-characteristic condition, an associated $L^\infty$-estimate and the existence of a convex entropy \cite{Jin1996, natalini1996convergence}. Furthermore, its Chapman-Enskog expansion yields the corresponding parabolic limit and provides a rigorous justification of the sub-characteristic condition. The main advantage of the integral invariant is that it provides a conserved quantity for the relaxation system by incorporating the coupled contributions of the solution variables in vector form. It also provides a useful measure of the conservation and consistency of numerical schemes and naturally connects the relaxation system with its limiting model as the relaxation parameter tends to zero. Despite these attractive theoretical properties, the numerical performance and practical applicability of this alternative relaxation formulation have not yet been systematically investigated using IMEX and finite volume central Nessyahu-Tadmor type numerical schemes.

Motivated by this observation, the present work extends Vasudeva Murthy's alternative relaxation framework to several representative hyperbolic balance law systems, including the shallow-water equations, the Broadwell model, and the Euler equations with heat transfer. The proposed methodology is further generalized to the two-dimensional Jin-Xin relaxation model. While the alternative relaxation formulation for the one-dimensional linear Jin-Xin model was established in \cite{murthy2007alternative}, this work develops analogous alternative relaxation formulations for the remaining systems. Moreover, the associated integral invariance properties are rigorously established for each proposed relaxation model. To assess the accuracy and robustness of the proposed relaxation formulations, numerical experiments are conducted for first-time as per author's best knowledge using the second-order finite-volume central schemes CS-EBT2 \cite{2026} and UCS2 \cite{Liotta}, together with a second-order IMEX Runge-Kutta method \cite{PareschiRusso}. Whenever analytical solutions are available, quantitative comparisons are performed directly. Otherwise, highly resolved reference solutions computed with the IMEX Runge-Kutta scheme are employed. The numerical results demonstrate that the proposed relaxation formulations provide accurate and robust approximations for both one- and two-dimensional hyperbolic balance law systems.

The remainder of this paper is organized as follows. Section~\ref{sec:relaxation} introduces the relaxation framework, develops the alternative relaxation models, and establishes their integral invariance properties. Section~\ref{sec:numerical} describes the numerical methods and presents a series of computational experiments that validate the proposed formulations. Finally, concluding remarks are provided in Section~\ref{sec:conclusion}.

\section{The Alternative Relaxation Systems and Their Asymptotic Limits}\label{sec:relaxation}
\begin{definition}[Relaxation system]
A hyperbolic system of balance laws of the form
\begin{equation}\label{govern:eq1}
\partial_t U + \partial_x F(U) = -\frac{1}{\tau} R(U), \quad U \in \mathbb{R}^N,
\end{equation}
is called a \emph{relaxation system} in the sense of Whitham \cite{Whitham1969} and Liu \cite{Liu1979} if there exists a constant matrix $Q \in \mathbb{R}^{n \times N}$ with rank $n < N$ such that
\[
Q R(U) = 0 \quad \text{for all } U.
\]
This condition implies the existence of $n$ independent conserved quantities $v = Q U$. It is assumed that each value of $v$ uniquely determines a corresponding local equilibrium state $U = E(v)$ satisfying
\[
R(E(v)) = 0.
\]
The image of the mapping $E$ defines the manifold of local equilibria associated with the relaxation operator $R$. For the original system and the matrix $Q$, one obtains $n$ conservation laws
satisfied by any solution of \eqref{govern:eq1}:
\[
\partial_t (Q U) + \partial_x (Q F(U)) = 0.
\]
Under the local relaxation approximation $U = E(v)$, these equations form a closed system, referred to as the reduced or equilibrium system, for the conserved variables $v = Q U$. This system can be written as
\begin{equation}
\partial_t v + \partial_x G(v) = 0,
\end{equation}
where the flux function is given by $G(v) = Q F(E(v))$. When the relaxation time $\tau$ is much smaller than the characteristic time
scale associated with wave propagation, the source term becomes stiff, leading to a rapid convergence toward the equilibrium manifold.
\end{definition}

\begin{definition}[Integral invariance]
\label{def:integralInvariance}
A system of balance laws is said to possess an \emph{integral invariance} if there exists a scalar function $\Phi(U)$ and an associated flux $\Psi(U)$ such that
\begin{equation}
\partial_t \Phi(U) + \partial_x \Psi(U) = 0.
\end{equation}
Consequently, under appropriate boundary conditions, the quantity
\begin{equation}
\int \Phi(U(x,t))\,dx,
\end{equation}
is preserved in time, that is,
\begin{equation}
\int \Phi(U(x,t))\,dx = \int \Phi(U(x,0))\,dx,
\qquad \text{for all } t \ge 0.
\end{equation}
\end{definition}

Having established the general definition of relaxation systems and their equilibrium dynamics, we now focus on a specific relaxation approximation that falls within this class. The alternative Jin-Xin model proposed by Vasudeva Murthy serves as a prototypical example, allowing us to illustrate the Chapman-Enskog procedure and the emergence of integral invariants in a simple and transparent setting.

\subsection{Jin-Xin Relaxation Approximation:}
The Jin-Xin relaxation model~\cite{jin1995runge,Sjin_1995} provides a classical framework for approximating scalar conservation laws by means of a hyperbolic system with stiff relaxation. The model is given by
\begin{equation}\label{eq:xinJin}
\begin{cases}
\partial_t u_\tau + \partial_x v_\tau = 0, \\
\partial_t v_\tau + \partial_x u_\tau
= -\dfrac{1}{\tau}\,\big(v_\tau - f(u_\tau)\big),
\end{cases}
\end{equation}
where $\tau>0$ denotes the relaxation parameter. In the limit $\tau \to 0$, the system converges formally to the scalar conservation law
\[
\partial_t u + \partial_x f(u) = 0.
\]
Moreover, a Chapman-Enskog expansion shows that the relaxation system admits a parabolic-type approximation of the form
\begin{equation}
\partial_t u_\tau + \partial_x f(u_\tau) = \tau\,\partial_x\!\left(\eta(u_\tau)\,\partial_x u_\tau\right),
\end{equation}
where the diffusion coefficient $\eta(u) =1 -[f'(u)]^2$ satisfies Liu's sub-characteristic condition, which ensures the stability and consistency of the relaxation approximation. In addition, the Jin-Xin model admits a conservative formulation that leads to an integral invariance in the sense of
Definition~\ref{def:integralInvariance}.

We recall an alternative relaxation formulation of the Jin-Xin model proposed by Vasudeva Murthy~\cite{murthy2007alternative}, given by the system
\begin{equation}\label{modify:xinJin}
\begin{cases}
\partial_t u_\tau + \partial_x v_\tau = v_\tau - au_\tau, \\
\partial_t v_\tau + \partial_x u_\tau = -\dfrac{1}{\tau}\,(v_\tau - f(u_\tau)).
\end{cases}
\end{equation}
This system represents a linear relaxation-type model and constitutes a particular example of the general balance law of the form \eqref{eq:balance-law}. Following the classical Jin-Xin framework \cite{jin1995runge}, we consider the specific flux choice $f(u)=a u$ and investigate the asymptotic behavior of the model in the small relaxation limit $\tau \to 0$ by means of the Chapman-Enskog expansion.

The system is supplemented with the initial conditions
\begin{equation}
\begin{cases}
u_\tau(x,0) = u_0(x),\\
v_\tau(x,0) = f(u_0(x)) = a u_0(x) = v_0(x).
\end{cases}
\end{equation}
We assume that the relaxation variable $v_\tau$ admits an expansion in powers of $\tau$ of the form
\begin{equation}
v_\tau = v_\tau^{(0)} + \tau v_\tau^{(1)} + \mathcal{O}(\tau^2),
\end{equation}
where the first-order correction is given by
\begin{equation}
v_\tau^{(1)} = -\eta(u_\tau)\partial_x u_\tau, \qquad \eta(u) =1-a^2.
\end{equation}

Substituting this expansion into the first equation of
\eqref{modify:xinJin} yields, at leading order, the reduced parabolic-type
equation
\begin{equation}
\partial_t u_\tau + \partial_x v_\tau =\tau\,\partial_x\!\left(\eta(u_\tau)\,\partial_x u_\tau\right)
- \tau\,\eta(u_\tau)\,\partial_x u_\tau.
\end{equation}
This equation is parabolic provided that the Liu's sub-characteristic condition $\eta(u) \ge 0$ is satisfied, which guarantees the stability of the relaxation approximation and its consistency with the underlying hyperbolic system.

An important structural property of the alternative Jin-Xin relaxation model is the presence of an exact conservation law,
\begin{equation}\label{eqV}
\partial_t \big(u_\tau + \tau v_\tau\big)
+ \partial_x \big(v_\tau + \tau u_\tau\big) = 0,
\end{equation}
which implies an integral invariance in the sense of
Definition~\ref{def:integralInvariance}. This conservation property enhances both the physical interpretability and the numerical robustness of the alternative Jin-Xin relaxation model. The main advantage of equation \eqref{eqV} is that it satisfies integral invariance 
\[
\int \left(u_\tau + \tau v_\tau\right)(\cdot,t)
=
\int \left(u_0(\cdot)+\tau f(u_0(\cdot))\right),
\]
for all \(t\), whereas the integral invariant for \eqref{eq:xinJin} is
\[
\int u_\tau(\cdot,t)
=
\int u_0(\cdot).
\] Note that the integral invariant provides a fundamental conservation property of the relaxation system by accounting for the coupled contributions of $u_\tau$ and $v_\tau.$
It serves as an important criterion for assessing the consistency and conservation properties of numerical approximations. Moreover, as $\tau \rightarrow 0$, this invariant naturally reduces to the corresponding integral invariant of the limiting conservation law and is particularly useful to assess the conservation and asymptotic-preserving properties of numerical schemes.

Motivated by the alternative relaxation framework of Vasudeva Murthy introduced above, we now develop new relaxation approximations for several physically relevant systems. In particular, we construct the alternative relaxation models for the shallow water equations, the Broadwell model, and the compressible Euler equations with heat transfer, following the same structural principles. For each system, we derive the corresponding diffusive limit via the Chapman-Enskog expansion, obtain a general conservative formulation, and establish the associated integral invariance property in the sense of Definition~\ref{def:integralInvariance}.

\subsection{Shallow Water Equations:}
The shallow water equations \cite{SB,russo} are a fundamental model for the description of free-surface flows under the hydrostatic assumption. In one spatial dimension, a relaxation approximation of the shallow water system can be written as
\begin{equation}\label{orig:shallow}
\begin{cases}
\partial_t h_\tau + \partial_x (h_\tau u_\tau) = 0, \\[2mm]
\partial_t (h_\tau u_\tau)+ \partial_x \left(\dfrac{(h_\tau u_\tau)^2}{h_\tau}+ \dfrac{1}{2} h_\tau^2\right)= -\dfrac{1}{\tau}\left((h_\tau u_\tau) - \dfrac{1}{2} h_\tau^2\right),
\end{cases}
\end{equation}
where $h_\tau$ denotes the water depth and $u_\tau$ the depth-averaged velocity. In the relaxation limit $\tau \to 0$, the system yields a consistent approximation of the underlying conservation law. A Chapman-Enskog expansion yields a parabolic-type approximation whose well-posedness is ensured under an appropriate sub-characteristic condition. Moreover, the relaxation system admits a conservative formulation, leading to an integral invariance in the sense of Definition~\ref{def:integralInvariance}.

We begin with an alternative relaxation approximation of the one-dimensional shallow water equations constructed within the general balance-law framework \eqref{govern:eq1}. The relaxation system is given by
\begin{equation}\label{modify:shallow}
\begin{cases}
\partial_t h_\tau + \partial_x (h_\tau u_\tau)
= (h_\tau u_\tau) - \dfrac{1}{2} h_\tau^2, \\[2mm]
\partial_t (h_\tau u_\tau)
+ \partial_x \left( \dfrac{(h_\tau u_\tau)^2}{h_\tau}
+ \dfrac{1}{2} h_\tau^2 \right)
= -\dfrac{1}{\tau}
\left( (h_\tau u_\tau) - \dfrac{1}{2} h_\tau^2 \right),
\end{cases}
\end{equation}
supplemented with the initial conditions
\begin{equation}
\begin{cases}
h_\tau(x,0) = h_0(x),\\
(h_\tau u_\tau)(x,0) = (h_0 u_0)(x).
\end{cases}
\end{equation}
Here, the relaxation term enforces convergence toward the equilibrium relation $h_\tau u_\tau = \tfrac{1}{2} h_\tau^2$ as $\tau \to 0$.
To investigate the diffusive limit of the system, we apply the Chapman-Enskog expansion and assume that the momentum admits the asymptotic expansion
\begin{equation}
(h_\tau u_\tau)= (h_\tau u_\tau)^{(0)}+ \tau (h_\tau u_\tau)^{(1)}+ \mathcal{O}(\tau^2),
\end{equation}
with
\begin{equation}
(h_\tau u_\tau)^{(0)} = \frac{h_\tau^2}{2},\qquad(h_\tau u_\tau)^{(1)}= \partial_x h_\tau\left( \frac{h_\tau^2 - 4 h_\tau}{4} \right).
\end{equation}
Substituting this expansion into the mass balance equation of
\eqref{modify:shallow}, we obtain the reduced parabolic equation
\begin{equation}
\partial_t h_\tau+ \partial_x \left( \frac{h_\tau^2}{2} \right)=\tau\,\partial_x\!\left(\left( \frac{4 h_\tau - h_\tau^2}{4} \right)\partial_x h_\tau\right)-\tau \left( \frac{4 h_\tau - h_\tau^2}{4} \right)\partial_x h_\tau.
\end{equation}
The equation is parabolic under Liu’s sub-characteristic condition
$0 < h_\tau < 4$.
Furthermore, the relaxation system satisfies the conservative identity
\begin{equation}
\partial_t \big(h_\tau + \tau (h_\tau u_\tau)\big)+ \partial_x \left((h_\tau u_\tau)+ \tau \left(\frac{(h_\tau u_\tau)^2}{h_\tau}+ \frac{1}{2} h_\tau^2\right) \right) = 0.
\end{equation}
Consequently, the quantity $\int \big(h_\tau(x,t) + \tau (h_\tau u_\tau)(x,t)\big)\,dx$ is preserved in time and therefore constitutes an integral invariant in the sense of Definition~\ref{def:integralInvariance}.

\subsection{Broadwell System:}
The Broadwell model \cite{broadwell1964shock,PareschiRusso} is a classical discrete-velocity kinetic system that serves as a prototypical example for the analysis of nonlinear hyperbolic systems with relaxation. In one spatial dimension, the macroscopic formulation of the Broadwell relaxation model can be written as
\begin{equation}\label{eq:broadwell}
\begin{cases}
\partial_t \rho_\tau + \partial_x m_\tau = 0, \\[4pt]
\partial_t m_\tau + \partial_x z_\tau = 0, \\[6pt]
\partial_t z_\tau + \partial_x m_\tau
= \dfrac{1}{2\tau}
\left(\rho_\tau^2 + m_\tau^2 - 2\rho_\tau z_\tau\right),
\end{cases}
\end{equation}
where $\rho_\tau$ denotes the density, $m_\tau = \rho_\tau u_\tau$ represents the momentum with velocity $u_\tau$, $z_\tau$ corresponds to the momentum flux and $\tau>0$ is the relaxation parameter. In the stiff relaxation limit $\tau \to 0$, the system formally reduces to its hydrodynamic equilibrium and yields the associated system of conservation laws. The Chapman-Enskog expansion further reveals the presence of diffusive corrections, leading to a parabolic-type approximation whose well-posedness is ensured under a suitable sub-characteristic condition. Moreover, the Broadwell relaxation system admits a conservative structure that gives rise to an integral invariance in the sense of Definition~\ref{def:integralInvariance}.

We now construct an alternative relaxation formulation of the Broadwell system within the general framework introduced earlier. The corresponding relaxation model reads
\begin{equation}\label{eq:alternative_broadwell}
\begin{cases}
\partial_t \rho_\tau + \partial_x m_\tau = 0, \\[4pt]
\partial_t m_\tau + \partial_x z_\tau
= -\dfrac{1}{2}\left(\rho_\tau^2 + m_\tau^2 - 2\rho_\tau z_\tau\right), \\[6pt]
\partial_t z_\tau + \partial_x m_\tau
= \dfrac{1}{2\tau}
\left(\rho_\tau^2 + m_\tau^2 - 2\rho_\tau z_\tau\right),
\end{cases}
\end{equation}
supplemented with the initial conditions
\begin{equation}
\begin{cases}
\rho_\tau(x,0) = \rho_0(x),\\
m_\tau(x,0) = m_0(x),\\
z_\tau(x,0) = z_0(x).
\end{cases}
\end{equation}
The relaxation parameter $\tau>0$ controls the rate of approach toward equilibrium. In the limit $\tau \to 0$, the equilibrium relation is given by
\begin{equation}
z_\tau \longrightarrow
\frac{1}{2\rho_\tau}\left(\rho_\tau^2 + m_\tau^2\right).
\end{equation}
The characteristic structure of system \eqref{eq:alternative_broadwell} yields eigenvalues $\lambda = (1,0,-1)$. Moreover, the Liu's sub-characteristic condition
\begin{equation}
\left| \frac{m_\tau}{\rho_\tau} \right| < 1,
\end{equation}
ensures hyperbolicity and stability of the relaxation approximation.

To derive the corresponding diffusive limit, we apply the Chapman-Enskog expansion and assume that $z_\tau$ admits the expansion
\begin{equation}
z_\tau = z_\tau^{(0)} + \tau z_\tau^{(1)} + \mathcal{O}(\tau^2),
\end{equation}
where the leading-order equilibrium term is
\begin{equation}
z_\tau^{(0)}= \frac{1}{2\rho_\tau}\left(\rho_\tau^2 + m_\tau^2\right).
\end{equation}
By substituting the expansion into \eqref{eq:alternative_broadwell} and matching terms of equal powers of $\tau$, the first-order correction is obtained as
\begin{equation}
z_\tau^{(1)}= -\frac{1}{\rho_\tau}\left[\left(\frac{1}{2}-\frac{m_\tau^2}{2\rho_\tau^2}\right)\partial_x(m_\tau)
+\left(\frac{m_\tau^3}{2\rho_\tau^3}-\frac{m_\tau}{2\rho_\tau}\right)\partial_x(\rho_\tau)\right].
\end{equation}
Substituting the expressions for $z_\tau^{(0)}$ and $z_\tau^{(1)}$ into the momentum equation yields the reduced parabolic-type system
\begin{equation}\label{parabolic_modified:Broadwell}
\begin{cases}
\partial_t \rho_\tau + \partial_x m_\tau = 0, \\[6pt]
\partial_t m_\tau+ \partial_x\!\left(\dfrac{1}{2}\rho_\tau+ \dfrac{1}{2\rho_\tau} m_\tau^2\right)
&=\tau\,\partial_x\!\left[\dfrac{1}{\rho_\tau}\left(\left(\dfrac{1}{2}-\dfrac{m_\tau^2}{2\rho_\tau^2}\right)\partial_x(m_\tau)
+\left(\dfrac{m_\tau^3}{2\rho_\tau^3}-\dfrac{m_\tau}{2\rho_\tau}\right)\partial_x(\rho_\tau)\right)\right] \\[6pt]
\qquad\qquad
&- \dfrac{\tau}{\rho_\tau}\left(\left(\dfrac{1}{2}-\dfrac{m_\tau^2}{2\rho_\tau^2}\right)\partial_x(m_\tau)
+\left(\dfrac{m_\tau^3}{2\rho_\tau^3}-\dfrac{m_\tau}{2\rho_\tau}\right)\partial_x(\rho_\tau)\right).
\end{cases}
\end{equation}
The system \eqref{parabolic_modified:Broadwell} is of parabolic type under the Liu's sub-characteristic condition, with the diffusive correction acting as a stabilizing mechanism.

Finally, the proposed alternative relaxation model for the Broadwell system retains a conservative structure
\begin{equation}
\partial_t \rho_\tau + \partial_x m_\tau = 0, \qquad
\partial_t (m_\tau + \tau z_\tau) + \partial_x (z_\tau + \tau m_\tau) = 0,
\end{equation}
which implies an integral invariant in the sense of
Definition~\ref{def:integralInvariance}, with conserved variable
\[
U_\tau = (\rho_\tau,\; m_\tau + \tau z_\tau)^T.
\]

\subsection{ Euler Equations with Heat Transfer:}
We first recall the one-dimensional compressible Euler equations with heat transfer in relaxation form \cite{Xu1999HeatTransfer}, which serve as the reference model for the subsequent alternative relaxation formulation. Thermal interaction with an external heat bath is modeled through a relaxation mechanism acting on the total energy. In terms of the density $\rho_\tau$, momentum $m_\tau = \rho_\tau u_\tau$, and total energy density $w_\tau = \rho_\tau E_\tau$, the system takes the form
\begin{equation}\label{eq:Euler_heat_original}
\begin{cases}
\partial_t \rho_\tau + \partial_x m_\tau = 0, \\[4pt]
\partial_t m_\tau+ \partial_x \!\left[\left(1-\dfrac{\gamma-1}{2}\right)\dfrac{m_\tau^2}{\rho_\tau}
+ (\gamma -1)w_\tau\right] = 0, \\[6pt]
\partial_t w_\tau+ \partial_x \!\left[\dfrac{m_\tau}{\rho_\tau}\left(\gamma w_\tau
- \dfrac{\gamma-1}{2}\dfrac{m_\tau^2}{\rho_\tau}\right)\right]
=\dfrac{1}{\tau}\left(\rho_\tau T_0
- \dfrac{1}{c_v}\left(w_\tau - \dfrac{m_\tau^2}{2\rho_\tau}\right)\right).
\end{cases}
\end{equation}

Here $E_\tau = e_\tau + \tfrac{u_\tau^2}{2}$ denotes the total energy per unit mass, with $e_\tau$ the internal energy. For a $\gamma$-law gas, the thermodynamic relations are $p = (\gamma -1)\rho_\tau e_\tau$ and $e_\tau = \mathcal{R}T/(\gamma -1)$, where $\mathcal{R}$ is the specific gas constant. The parameter $\tau>0$ represents the thermal relaxation time, and $T_0$ denotes the prescribed temperature of the heat bath.

In the stiff relaxation limit $\tau \to 0$, the system formally reduces to the compressible Euler equations with an equilibrium temperature constraint. The model preserves the conservative structure of the mass and momentum equations and admits an exact conservation law leading to an integral invariance in the sense of Definition~\ref{def:integralInvariance}.

Within the general balance law framework \eqref{govern:eq1}, we now introduce another alternative relaxation formulation of the system,
\begin{equation}\label{Modified:Euler_with_heat_transfer}
\begin{cases}
\partial_t \rho_\tau + \partial_x m_\tau = 0, \\[4pt]
\partial_t m_\tau+ \partial_x \!\left[\left(1-\dfrac{\gamma-1}{2}\right)
\dfrac{m_\tau^2}{\rho_\tau}+ (\gamma -1)w_\tau\right]
= -\!\left(\rho_\tau T_0
- \dfrac{1}{c_v}\!\left(w_\tau - \dfrac{m_\tau^2}{2\rho_\tau}\right)\!\right), \\[6pt]
\partial_t w_\tau+ \partial_x \!\left[\dfrac{m_\tau}{\rho_\tau}
\!\left(\gamma w_\tau- \dfrac{\gamma-1}{2}
\dfrac{m_\tau^2}{\rho_\tau}\right)\!\right]
=\dfrac{1}{\tau}\!\left(\rho_\tau T_0
- \dfrac{1}{c_v}\!\left(w_\tau - \dfrac{m_\tau^2}{2\rho_\tau}\right)\!\right).
\end{cases}
\end{equation}
The proposed relaxation system is considered together with the prescribed initial data
\begin{equation}
\begin{cases}
\rho_\tau(x,0) = \rho_0(x),\\
m_\tau(x,0) = m_0(x),\\
w_\tau(x,0) = w_0(x).
\end{cases}
\end{equation}

For the hyperbolic part of the alternative relaxation system, the characteristic speeds are given by $\lambda = (u_\tau - c,\, u_\tau,\, u_\tau + c)$, where $c = \sqrt{\gamma p/\rho_\tau}$ denotes the sound speed. Equivalently, eigenvalues can be expressed as
\[
\lambda \approx \left(u_\tau - \sqrt{\gamma(\gamma -1)c_v T},\, u_\tau,\, u_\tau + \sqrt{\gamma(\gamma -1)c_v T}\right).
\]
The corresponding sub-characteristic condition for the relaxation approximation is $T_0 < \gamma T$, which ensures stability and consistency of the model.

In the equilibrium limit $\tau \to 0$, the internal energy relaxes toward the equilibrium value associated with the bath temperature $T_0$, yielding
\begin{equation}
w_\tau \longrightarrow \frac{m_\tau^2}{2\rho_\tau} + c_v \rho_\tau T_0.
\end{equation}

To derive the diffusive limit, we apply the Chapman-Enskog expansion to $w_\tau$ and assume
\begin{equation}
w_\tau = w_\tau^{(0)} + \tau w_\tau^{(1)} + \mathcal{O}(\tau^2),
\end{equation}
where
\begin{align}
w_\tau^{(0)} &= \frac{m_\tau^2}{2\rho_\tau} + c_v \rho_\tau T_0, \\[4pt]
w_\tau^{(1)} &=
- c_v^2 T_0 (\gamma -1)
\!\left(\partial_x m_{\tau}
- \frac{m_\tau}{\rho_\tau}\partial_x \rho_{\tau}\right).
\end{align}

Substituting these expressions into the governing equations yields the reduced parabolic-type system
\begin{equation}\label{parabolic_modified:Euler_with_heat_transfer}
\begin{cases}
\partial_t \rho_\tau + \partial_x m_\tau = 0, \\[6pt]
\partial_t m_\tau
+ \partial_x \!\left[(\gamma -1)c_v T_0\rho_\tau
+ \dfrac{m_\tau^2}{\rho_\tau}\right]
= \tau c_v^2 T_0 (\gamma -1)^2
\partial_x\!\left(\rho_\tau \partial_x u_{\tau}\right)
- \tau c_v T_0 (\gamma -1)\rho_\tau \partial_x u_{\tau}.
\end{cases}
\end{equation}

The system \eqref{parabolic_modified:Euler_with_heat_transfer} is parabolic under the sub-characteristic condition and represents the diffusive limit of the thermal relaxation process.

Finally, the proposed alternative relaxation model for Euler equations with heat transfer, preserves a conservative structure,
\begin{equation}\label{conservation:modified_Euler_with_heat_transfer}
\partial_t \rho_\tau + \partial_x m_\tau = 0, \qquad
\partial_t (m_\tau + \tau w_\tau)
+ \partial_x \!\left(m_\tau + \tau w_\tau\right) = 0,
\end{equation}
which implies an integral invariant in the sense of Definition~\ref{def:integralInvariance}, with conserved variable
\[
U_\tau = (\rho_\tau,\; m_\tau + \tau w_\tau)^T.
\]

\subsection*{Relaxation Approximation for Two-Dimensional Hyperbolic Balance Laws:}
In this subsection, we extend the proposed relaxation framework to two-dimensional hyperbolic balance laws. 
Consider the two-dimensional hyperbolic balance law of the form
\begin{equation}
U_t + F_1(U)_x + F_2(U)_y = S(U),
\end{equation}
where $U$ denotes the vector of conserved variables, $F_1(U)$ and $F_2(U)$ represent the flux functions in the $x$- and $y$-directions, respectively, and $S(U)$ is the source term. Motivated by the one-dimensional alternative relaxation formulation, we propose the following two-dimensional relaxation system.
\subsection{2D Jin-Xin relaxation model:}
A widely used relaxation framework is the Jin--Xin model, which in two dimensions is given by
\begin{equation}\label{Xin_Jin_2D_model}
\begin{cases}
\partial_t u_\tau + \partial_x v_\tau + \partial_y w_\tau
= 0,\\[2mm]
\partial_t v_\tau + \partial_x u_\tau
= -\dfrac{1}{\tau}\left(v_\tau-a u_\tau\right),\\[2mm]
\partial_t w_\tau + \partial_y u_\tau
= -\dfrac{1}{\tau}\left(w_\tau-b u_\tau\right).
\end{cases}
\end{equation}
Here, $u$ denotes the conserved variable, whereas $v$ and $w$ are auxiliary relaxation variables corresponding to the fluxes in the $x$- and $y$-directions, respectively. The constants $a$ and $b$ are prescribed relaxation parameters. The relaxation system is designed so that, as the relaxation parameter tends to zero, the auxiliary variables approach their equilibrium states and the original hyperbolic equation is recovered.

Indeed, in the limit $\tau\to0$, the stiff source terms enforce the equilibrium relations
\[
v_\tau=au_\tau,\quad w_\tau=bu_\tau.
\]
Substituting these equilibrium relations into the first equation of \eqref{modify:Xin_Jin_2D_model} yields the reduced conservation law
\begin{equation}
\partial_t u_\tau+a\,\partial_x u_\tau+b\,\partial_y u_\tau=0,
\end{equation}
which represents the asymptotic limit of the proposed relaxation system. 
Hence, the relaxation approximation is formally consistent with the underlying hyperbolic conservation law.

We introduce an alternative relaxation approximation of the two-dimensional Jin–Xin model, which is formulated as follows
\begin{equation}\label{modify:Xin_Jin_2D_model}
\begin{cases}
\partial_t u_\tau + \partial_x v_\tau + \partial_y w_\tau
= v_\tau + w_\tau -(a+b)u_\tau,\\[2mm]
\partial_t v_\tau + \partial_x u_\tau
= -\dfrac{1}{\tau}\left(v_\tau-a u_\tau\right),\\[2mm]
\partial_t w_\tau + \partial_y u_\tau
= -\dfrac{1}{\tau}\left(w_\tau-b u_\tau\right).
\end{cases}
\end{equation}
The system is supplemented with the initial conditions 
\begin{equation} 
\begin{cases} u_\tau(x,y,0) = u_0(x,y),\\ v_\tau(x,y,0) = a u_0(x,y) = v_0(x,y).\\ w_\tau(x,y,0) = b u_0(x,y) = w_0(x,y). 
\end{cases} 
\end{equation} 
We assume that the relaxation variable $v_\tau$ admits an expansion in powers of $\tau$ of the form 
\begin{equation} 
v_\tau = v_\tau^{(0)} + \tau v_\tau^{(1)} + \mathcal{O}(\tau^2),\, w_\tau = w_\tau^{(0)} + \tau w_\tau^{(1)} + \mathcal{O}(\tau^2). 
\end{equation} 
where the first-order correction is given by 
\begin{equation}
v_\tau^{(1)} = -\eta(u_\tau)\partial_x u_\tau, \quad \eta(u) =1-a^2,\qquad w_\tau^{(1)} = -\beta(u_\tau)\partial_y u_\tau, \quad \beta(u) =1-b^2, \end{equation} 
Substituting this expansion into the first equation of \eqref{modify:Xin_Jin_2D_model} yields, at leading order, the reduced parabolic-type equation \begin{equation} 
\partial_t u_\tau + \partial_x v_\tau =\tau\,\left\{\partial_x\!\left(\eta(u_\tau)\,\partial_x u_\tau\right) + \partial_y\!\left(\beta(u_\tau)\,\partial_y u_\tau\right)\right\} - \tau\,\left\{\eta(u_\tau)\,\partial_x u_\tau + \beta(u_\tau)\,\partial_y u_\tau \right\}. 
\end{equation} 
The above equation is parabolic provided that the Liu sub-characteristic conditions
\[
\eta(u_\tau)\ge0,\quad
\beta(u_\tau)\ge0,
\]
are satisfied. These conditions guarantee the positivity of the diffusion coefficients, thereby ensuring the stability of the relaxation approximation and its consistency with the underlying hyperbolic system.

An important structural property of the proposed alternative two-dimensional Jin--Xin relaxation model is that it satisfies the conservation law
\begin{equation} \label{conserve:2dxinjin} \partial_t \big(u_\tau + \tau v_\tau + \tau w_\tau\big) + \partial_x \big(v_\tau + \tau u_\tau\big) + \partial_y \big(w_\tau + \tau u_\tau\big) = 0. \end{equation}
Integrating \eqref{conserve:2dxinjin} over the spatial domain $\Omega$ gives \begin{equation*} 
\int_\Omega \partial_t \big(u_\tau + \tau v_\tau + \tau w_\tau\big) dx\, dy\, dt + \int_\Omega \partial_x \big(v_\tau + \tau u_\tau\big) dx\, dy\, dt \int_\Omega \partial_y \big(w_\tau + \tau u_\tau\big) dx\, dy\, dt=0. 
\end{equation*} 
Under appropriate boundary conditions, such as periodic or homogeneous zero-flux boundary conditions, the boundary integrals vanish. Consequently,
\begin{equation}
\frac{d}{dt}
\int_\Omega
\big(u_\tau(x,y,t)+\tau v_\tau(x,y,t)+\tau w_\tau(x,y,t)\big)\,dx\,dy
=0.
\end{equation}
\begin{equation} \int_\Omega \big(u_\tau(x,y,t) + \tau v_\tau(x,y,t) + \tau w_\tau (x,y,t)\big) dx\, dy = \int_\Omega \big(u_0(x,y) + \tau v_0(x,y) + \tau w_0(x,y)\big) dx\, dy. \end{equation} 
The integral invariance in the sense of Definition~\ref{def:integralInvariance}, the proposed relaxation model preserves the corresponding integral invariant. The preservation of this invariant provides a physically meaningful conserved quantity and contributes to the physical consistency and numerical robustness of the proposed alternative two-dimensional Jin--Xin relaxation model.

\section{Numerical Approximation and Validation}\label{sec:numerical}

This section describes the numerical methods employed to validate the alternative relaxation formulations developed in the previous section. Both one- and two-dimensional hyperbolic balance laws are considered. In one spatial dimension, the governing system is given by
\begin{equation}
\partial_t U+\partial_xF(U)=S(U),
\label{eq:hyperbolic}
\end{equation}
where $U$ denotes the vector of conserved variables, $F(U)$ is the flux function, and $S(U)$ represents the relaxation source term.
The corresponding two-dimensional system is
\begin{equation}
\partial_t U + \partial_x F_1 (U) + \partial_y F_2(U) = S(U),
\label{eq:hyperbolic2D}
\end{equation}
where $F_1(U)$ and $F_2(U)$ denote the fluxes in the $x$- and $y$-directions, respectively.

For the numerical approximation, we employ three second-order numerical schemes available in the literature: the recently proposed second-order finite-volume central schemes CS-EBT2~\cite{2026} and UCS2~\cite{Liotta}, together with the second-order IMEX Runge--Kutta (IMEX RK2) scheme~\cite{PareschiRusso}. 
Through a series of one- and two-dimensional benchmark problems, we investigate the accuracy, robustness, consistency with the diffusive limit, and preservation of the associated integral invariance properties.

\subsection{Second Order FV Central Schemes}
The one- and two-dimensional hyperbolic balance laws are discretized using the uniformly second-order implicit--explicit (IMEX) central schemes CS-EBT2~\cite{2026} and UCS2~\cite{Liotta}. Both schemes are based on the staggered-grid finite-volume framework of Nessyahu and Tadmor (NT)~\cite{NT}, in which the convective fluxes are approximated by a Riemann-solver-free central discretization, while the relaxation source term is treated semi-implicitly. Consequently, these methods achieve second-order accuracy in both space and time and provide stable and robust approximations for hyperbolic balance laws with stiff source terms. The formulations of the CS-EBT2 and UCS2 schemes are presented below.

For the one-dimensional formulation, let $\bar{U}_i(t)$ denote the cell average of $U$ over the interval $[x_{i-\frac12},x_{i+\frac12}]$,
\begin{equation}
\bar{U}_i(t) = \frac{1}{\Delta x} \int_{x_{i-\frac12}}^{x_{i+\frac12}} U(x,t)\,dx.
\end{equation}
Within each cell, a piecewise linear reconstruction is employed,
\begin{equation}
L_i(x,t) = \bar{U}_i(t) + (x-x_i)\frac{\bar{U}'_i}{\Delta x},
\qquad x\in[x_{i-\frac12},x_{i+\frac12}],
\end{equation}
where the slopes $\bar{U}'_i$ are computed using the MinMod limiter
\begin{equation}
\bar{U}'_i = \operatorname{MM}\!\left(\bar{U}_{i+1}-\bar{U}_i,\; \bar{U}_i-\bar{U}_{i-1}\right),
\end{equation}
with limiter
\begin{equation}\label{eq:minmod}
\operatorname{MM}(a,b) =
\begin{cases}
\operatorname{sgn}(a)\,\min(|a|,|b|), & \text{if } \operatorname{sgn}(a)=\operatorname{sgn}(b),\\
0, & \text{otherwise},
\end{cases}
\end{equation}
which provides a non-oscillatory reconstruction while retaining second-order accuracy in smooth regions.
Integrating the balance law over the space-time control volume 
$[x_i, x_{i+1}] \times [t^n, t^{n+1}]$, and using the piecewise linear polynomial reconstruction for $U(x,t^n)$ in $[x_i, x_{i+\frac12}]$ and $[x_{i+\frac12}, x_{i+1}],$ the cell average over the spatial cell $[x_i, x_{i+1}]$ can be expressed as
\begin{equation}\label{corrector1}
    \frac{1}{\Delta x} \int_{x_i}^{x_{i+1}} U(x,t^n)\, dx 
    = \frac{1}{2} \big(\bar{U}_i^n + \bar{U}_{i+1}^n \big) 
      + \frac{1}{8} \big(\bar{U}'_i - \bar{U}'_{i+1} \big).
\end{equation}
and the contribution of the flux term over the space-time cell is given by
\begin{equation}\label{corrector2}
\frac{1}{\Delta x}\int_{t^n}^{t^{n+1}}
\Bigl[F(U(x_{i+1},t))-F(U(x_i,t))\Bigr]\,dt.
\end{equation}
This temporal integral is approximated by the midpoint rule,
\begin{equation}\label{corrector3}
\int_{t^n}^{t^{n+1}}
\Bigl[F(U(x_{i+1},t))-F(U(x_i,t))\Bigr]\,dt
\;\approx\;
\Delta t\Bigl[
F(\bar{U}_{i+1}^{\,n+\frac12})
-
F(\bar{U}_i^{\,n+\frac12})
\Bigr],
\end{equation}
where the staggered half-time interface states
$\bar{U}_i^{\,n+\frac12}$ are obtained from the NT-type predictor step. 

\subsubsection{One Dimensional CS-EBT2 Scheme}
The space-time integral of the source term,
\begin{equation}
\mathcal{I}_S =
\int_{t^n}^{t^{n+1}}
\int_{x_{i}}^{x_{i+1}}
S(U(x,t))\,dx\,dt,
\end{equation}
is approximated using a semi-implicit backward Taylor expansion,
\begin{equation}
\mathcal{I}_S \approx
\frac{\Delta x\,\Delta t}{2} \Bigg[ 2S(\bar{U}_{i+1/2}^{n+1})
- \frac{\Delta t}{2}
\left(\frac{\partial S}{\partial U}(\bar{U}_i^n) + \frac{\partial S}{\partial U}(\bar{U}_{i+1}^n) \right)
\left(S(\bar{U}_{i+1/2}^{n+1})
- \frac{F_{i+1}^n - F_i^n}{\Delta x} \right)
\Bigg].
\end{equation}
The source term is handled using an implicit treatment.

Combining the above discretizations yields the CS-EBT2 predictor-corrector formulation. The predictor stage reads
\begin{align}
\bar{U}_i^{\,n+\frac12}
=\bar{U}_i^n+\frac{\Delta t}{2}\left( G(\bar{U}_i^{\,n+\frac12}) - \frac{F'_i}{\Delta x}
\right).
\end{align}
The corrector stage is given by
\begin{equation}
\begin{aligned}
\bar{U}_{i+\frac12}^{\,n+1}=&\frac12
\left(
\bar{U}_i^n+\bar{U}_{i+1}^n
\right)
+\frac18
\left(\bar{U'}_i-\bar{U'}_{j+1}
\right) -\lambda
\left[F(\bar{U}_{i+1}^{\,n+\frac12})-F(\bar{U}_i^{\,n+\frac12})
\right] \\
&
+\frac{\Delta t}{2}
\Bigg[
2S(\bar{U}_{i+\frac12}^{\,n+1})-\frac{\Delta t}{2}
\left(
\frac{\partial S}{\partial \bar{U}}(\bar{U}_i^n)
+
\frac{\partial S}{\partial \bar{U}}(\bar{U}_{i+1}^n)
\right)
\left(
S(\bar{U}_{i+\frac12}^{\,n+1})
-
\frac{F_{i+1}^n-F_i^n}{\Delta x}
\right)
\Bigg].
\end{aligned}
\end{equation}
where $\lambda=\Delta t/\Delta x$ and the MinMod limiter is used in the computation of $\bar{U}'_i$ and $F'_i$.

The CS-EBT2 scheme is second-order accurate in both space and time, employs an implicit treatment of stiff source terms, and preserves the conservative structure of the balance law. Under suitable boundary conditions, the corresponding integral invariants are also preserved. The authors in \cite{2026} demonstrated that the scheme possesses a larger stability region, allowing for the use of larger CFL numbers and, consequently, larger time steps. This enhanced stability is particularly advantageous for stiff relaxation problems. Moreover, it requires lower CPU time and offers a simpler and more efficient computational implementation.
\subsubsection{One Dimensional UCS2 Scheme}
The space-time integral of the source term,
\begin{equation}
\mathcal{I}_S =
\int_{t^n}^{t^{n+1}}
\int_{x_{i}}^{x_{i+1}}
S(U(x,t))\,dx\,dt,
\end{equation}
is approximated using a Radau quadrature rule,
\begin{equation}
\mathcal{I}_S \approx
\Delta x\,\Delta t
\left(
\frac{3}{8}S(\bar{U}_i^{\,n+\frac13})
+
\frac{3}{8}S(\bar{U}_{i+1}^{\,n+\frac13})
+
\frac{1}{4}S(\bar{U}_{i+\frac12}^{\,n+1})
\right).
\end{equation}
The first two evaluations rely on predictor states, while the last term is treated implicitly at the new time level, ensuring stability in the stiff relaxation regime.

Combining the above discretizations yields the UCS2 predictor-corrector formulation. The predictor stage reads
\begin{align}
\bar{U}_i^{\,n+\frac12}
&=
\bar{U}_i^n
+
\frac{\Delta t}{2}
\left(
S(\bar{U}_i^{\,n+\frac12})
-
\frac{F'_i}{\Delta x}
\right),\\
\bar{U}_i^{\,n+\frac13}
&=
\bar{U}_i^n
+
\frac{\Delta t}{3}
\left(
S(\bar{U}_i^{\,n+\frac13})
-
\frac{F'_i}{\Delta x}
\right).
\end{align}
The corrector stage is given by
\begin{align}
\bar{U}_{i+\frac12}^{\,n+1}
&=
\frac{1}{2}
(\bar{U}_i^n+\bar{U}_{i+1}^n)
+
\frac{1}{8}
(\bar{U}'_i-\bar{U}'_{i+1})
-
\lambda
\Bigl[
F(\bar{U}_{i+1}^{\,n+\frac12})
-
F(\bar{U}_i^{\,n+\frac12})
\Bigr]
\notag\\
&\quad
+
\Delta t
\left(
\frac{3}{8}S(\bar{U}_i^{\,n+\frac13})
+
\frac{3}{8}S(\bar{U}_{i+1}^{\,n+\frac13})
+
\frac{1}{4}S(\bar{U}_{i+\frac12}^{\,n+1})
\right),
\end{align}
where $\lambda=\Delta t/\Delta x$ and the MinMod limiter is used in the computation of $\bar{U}'_i$ and $F'_i$.

The UCS2 scheme is second-order accurate in both space and time, employs an implicit treatment of stiff source terms while some terms are evaluated using predictor steps, and preserves the conservative structure of the balance law. Under appropriate boundary conditions, the associated integral invariants are also preserved. Compared with the CS-EBT2 scheme, the UCS2 scheme has a relatively smaller stability region and may therefore require more restrictive CFL conditions for certain problems. Nevertheless, it provides uniform consistency across the considered regimes. Moreover, the UCS2 scheme involves a larger number of predictor stages, resulting in a higher computational cost and CPU time compared with the CS-EBT2 scheme.

The above one-dimensional finite-volume construction extends naturally
to the two-dimensional setting. Consider the hyperbolic balance law
\eqref{eq:hyperbolic2D} on $(x,y)\in\mathbb{R}^2$ and $t>0$, supplemented
with the initial condition
\begin{equation}
U(x,y,0)=U_0(x,y).
\end{equation}
For the numerical approximation, we introduce a Cartesian grid with mesh sizes
\[
\Delta x=x_{i+\frac12}-x_{i-\frac12},
\qquad
\Delta y=y_{j+\frac12}-y_{j-\frac12},
\]
and time step
\[
\Delta t=t^{n+1}-t^n.
\]
The computational domain is partitioned into uniform Cartesian cells
\[
I_{i,j}
=
[x_{i-\frac12},x_{i+\frac12}]
\times
[y_{j-\frac12},y_{j+\frac12}],
\qquad
i=1,\ldots,N_x,\quad
j=1,\ldots,N_y,
\]
with cell centers
\[
(x_i,y_j)
=
\left(
\frac{x_{i-\frac12}+x_{i+\frac12}}{2},
\frac{y_{j-\frac12}+y_{j+\frac12}}{2}
\right).
\]
The corresponding cell average at time $t^n$ is defined by
\begin{equation}
\bar{U}_{i,j}^n
=
\frac{1}{\Delta x\,\Delta y}
\int_{x_{i-\frac12}}^{x_{i+\frac12}}
\int_{y_{j-\frac12}}^{y_{j+\frac12}}
U(x,y,t^n)\,dy\,dx.
\end{equation}
A piecewise linear reconstruction is then employed within each Cartesian cell, with the slopes in the $x$- and $y$-directions computed using the MinMod limiter. The staggered-grid predictor states and numerical fluxes are constructed in both spatial directions following the one-dimensional formulation. Consequently, the two-dimensional finite-volume update involves the numerical fluxes across the four faces of each control volume, together with the corresponding space-time approximation of the relaxation source term. The resulting two-dimensional finite-volume formulation involves
numerical fluxes through the four faces of each Cartesian control volume together with the corresponding space-time approximation of the
relaxation source term. The CS-EBT2 and UCS2 schemes differ in the treatment of this source contribution, as described below.
\subsubsection{Two Dimensional CS-EBT2 Scheme}
The CS-EBT2 scheme for the two-dimensional problem \cite{2026} is obtained by applying the semi-implicit source discretization to the two-dimensional
finite-volume formulation described above. The resulting predictor and corrector steps are given by\\
\noindent\textbf{Predictor step:}\label{CSEBT2:scheme2d}  
The predictor state $\bar{U}_{i,j}^{\,n+\frac12}$ is obtained by solving
\begin{equation}
\bar{U}_{i,j}^{\,n+\frac12} = \bar{U}_{i,j}^{\,n} - \frac{\Delta t} {2} \left[ S\! \left( \bar{U}_{i,j}^{\,n+\frac12}\right) - \frac{(F_1)_{i,j}^{\dagger}}{\Delta x}
-\frac{(F_2)_{i,j}^{\ast}}{\Delta y} \right].
\end{equation}
\noindent\textbf{Corrector step:}  Then update to $\bar{U}_{i+\frac12,j+\frac12}^{\,n+1}$ by
\begin{equation}
\begin{aligned}
\bar{U}_{i+\frac12,j+\frac12}^{\,n+1}
&= \left\langle
\dfrac{1}{4}\left(\bar{U}_{i,\cdot}^n + \bar{U}_{i+1,\cdot}^n\right)
-
\dfrac{1}{8}\left(\bar{U}^{\dagger}_{i,\cdot} - \bar{U}^{\dagger}_{i+1,\cdot}\right)
-
\lambda\left[
F_1(\bar{U}_{i+1,\cdot}^{\,n+1/2})
-
F_1(\bar{U}_{i,\cdot}^{\,n+1/2})
\right]
\right.\\
&\qquad\left.
+\dfrac{\Delta t^2}{4}
\left(
\dfrac{\partial S}{\partial U}(\bar{U}_{i,\cdot}^n)
+
\dfrac{\partial S}{\partial U}(\bar{U}_{i+1,\cdot}^n)
\right)
\left(
\dfrac{(F_1)_{i+1,\cdot}^n-(F_1)_{i,\cdot}^n}{\Delta x}
\right)
\right\rangle_{j+1/2}\\
&\quad+
\left\langle
\dfrac{1}{4}\left(\bar{U}_{\cdot,j}^n + \bar{U}_{\cdot,j+1}^n\right)
-
\dfrac{1}{8}\left(\bar{U}^{\ast}_{\cdot,j} - \bar{U}^{\ast}_{\cdot,j+1}\right)
-
\mu\left[
F_2(\bar{U}_{\cdot,j+1}^{\,n+1/2})
-
F_2(\bar{U}_{\cdot,j}^{\,n+1/2})
\right]
\right.\\
&\qquad\left.
+\dfrac{\Delta t^2}{4}
\left(
\dfrac{\partial S}{\partial U}(\bar{U}_{\cdot,j}^n)
+
\dfrac{\partial S}{\partial U}(\bar{U}_{\cdot,j+1}^n)
\right)
\left(
\dfrac{(F_2)_{\cdot,j+1}^n-(F_2)_{\cdot,j}^n}{\Delta y}
\right)
\right\rangle_{i+1/2}\\
&\quad+
\dfrac{\Delta t}{2}
\left[
2S(\bar{U}_{i+\frac12,j+\frac12}^{\,n+1})
-\Delta t
\dfrac{\partial S}{\partial U}
(\bar{U}_{i+\frac12,j+\frac12}^{\,n})
S(\bar{U}_{i+\frac12,j+\frac12}^{\,n+1})
\right].
\end{aligned}
\end{equation}
\subsubsection{Two Dimensional UCS2 Scheme}
The two-dimensional UCS2 scheme follows the same finite-volume construction, with the source contribution approximated using the Radau quadrature formulation. The explicit formulation presented in this work is developed as an extension of the UCS2 scheme proposed in ~\cite{Liotta}; to the best of our knowledge, such a formulation is reported here for the first time in the literature. The resulting update is given by
\medskip

\noindent\textbf{Predictor step:}
\label{UCS2:scheme2d} 
Solve for $\bar{U}_{i,j}^{\,n+\frac12}$ and
$\bar{U}_{i,j}^{\,n+\frac13}$ from
\begin{equation}
\bar{U}_{i,j}^{\,n+\frac12}=\bar{U}_{i,j}^{\,n}
+\dfrac{\Delta t}{2}
\left[S\!\left(\bar{U}_{i,j}^{\,n+\frac12}\right)
-\dfrac{(F_1)^{\dagger}_{i,j}}{\Delta x}
-\dfrac{(F_2)^{\ast}_{i,j}}{\Delta y}
\right].
\end{equation}
\begin{equation}
\bar{U}_{i,j}^{\,n+\frac13} = \bar{U}_{i,j}^{\,n}
+ \dfrac{\Delta t}{3}
\left[S\!\left(\bar{U}_{i,j}^{\,n+\frac13}\right)
-\dfrac{(F_1)^{\dagger}_{i,j}}{\Delta x}
-\dfrac{(F_2)^{\ast}_{i,j}}{\Delta y}
\right].
\end{equation}

\noindent\textbf{Corrector step:} Then update to $\bar{U}_{i+\frac12,j+\frac12}^{\,n+1}$ by
\begin{equation}
\begin{aligned}
\bar{U}_{i+\frac12,j+\frac12}^{\,n+1}
&=\left\langle
\dfrac{1}{4}
\left(\bar{U}_{i,\cdot}^{\,n}+\bar{U}_{i+1,\cdot}^{\,n}\right)+\dfrac{1}{8}
\left(\bar{U}^{\dagger}_{i,\cdot}
-\bar{U}^{\dagger}_{i+1,\cdot}
\right)-\lambda\left[ F_1\!\left(\bar{U}_{i+1,\cdot}^{\,n+\frac12}\right)-F_1\!\left(\bar{U}_{i,\cdot}^{\,n+\frac12}\right)
\right]
\right.\\
&\qquad\left.
+\Delta t
\left(\dfrac{3}{16}
S\!\left(\bar{U}_{i,\cdot}^{\,n+\frac13}\right)
+\dfrac{3}{16}
S\!\left(\bar{U}_{i+1,\cdot}^{\,n+\frac13}\right)
\right)\right\rangle_{j+\frac12}
\\
&\quad+
\left\langle
\dfrac{1}{4}
\left(\bar{U}_{\cdot,j}^{\,n}+\bar{U}_{\cdot,j+1}^{\,n}\right)
+\dfrac{1}{8}
\left(\bar{U}^{\ast}_{\cdot,j}
-\bar{U}^{\ast}_{\cdot,j+1}\right)
-\mu
\left[F_2\!\left(\bar{U}_{\cdot,j+1}^{\,n+\frac12}\right)
-F_2\!\left(\bar{U}_{\cdot,j}^{\,n+\frac12}\right)
\right]\right.\\
&\qquad\left.
+\Delta t \left(
\dfrac{3}{16}
S\!\left(\bar{U}_{\cdot,j}^{\,n+\frac13}\right)+
\dfrac{3}{16}
S\!\left(\bar{U}_{\cdot,j+1}^{\,n+\frac13}\right)
\right)
\right\rangle_{i+\frac12}
+
\dfrac{\Delta t}{4}
S\!\left(
\bar{U}_{i+\frac12,j+\frac12}^{\,n+1}
\right).
\end{aligned}
\end{equation}
Here,
$\langle \bar{U}_{\cdot,j} \rangle_{i+\frac12}
:= \frac12(\bar{U}_{i,j}+\bar{U}_{i+1,j})$
and
$\langle \bar{U}_{i,\cdot} \rangle_{j+\frac12}
:= \frac12(\bar{U}_{i,j}+\bar{U}_{i,j+1})$
denote the staggered averages in the $x$- and $y$-directions, respectively. Moreover, we define
$\lambda=\Delta t/\Delta x$ and $\mu=\Delta t/\Delta y$.
Here,
$\bar{U}^{\dagger}=\partial \bar{U}/\partial x$
and
$\bar{U}^{\ast}=\partial \bar{U}/\partial y$.
The quantities $\bar{U}^{\dagger}_i$, $\bar{U}^{\ast}_j$, $(F_1)^{\dagger}_{i,j}$, and $(F_2)^{\ast}_{i,j}$ are evaluated using the MinMod limiter defined in \eqref{eq:minmod}.
\subsection{Second-Order IMEX Runge-Kutta Scheme}
For the numerical approximation of the relaxation balance laws, we employ the second-order implicit-explicit Runge-Kutta (IMEX RK2) scheme introduced by Pareschi and Russo~\cite{PareschiRusso} and further analyzed by Boscarino~\cite{Boscarino2009}. The method is particularly suitable for relaxation systems in which the convective and source terms evolve on different time scales. The convective terms are treated explicitly, while the relaxation source terms are treated implicitly at the Runge-Kutta stages. This explicit-implicit treatment alleviates the severe time-step restrictions associated with an implicit treatment of stiff source terms and provides suitable stability properties in the stiff-relaxation regime. In the present work, the same second-order IMEX RK2 time discretization is applied to both one- and two-dimensional balance laws. The spatial finite-volume operators and corresponding flux contributions are formulated according to the dimension of the problem, while the Runge-Kutta coefficients remain unchanged. The one- and two-dimensional formulations are described separately below.

\subsubsection{One dimensional IMEX Runge-Kutta Scheme}
Let the spatial domain be $\Omega=[a,b]\subset\mathbb{R}$, which is discretized into $N$ uniform cells
\[
I_i=[x_{i-\frac12},x_{i+\frac12}], \qquad i=1,\dots,N,
\]
with cell size $\Delta x=(b-a)/N$ and cell centers $x_i=\tfrac12(x_{i-\frac12}+x_{i+\frac12})$. The time interval $[0,T]$ is discretized into time levels $t^n=n\,\Delta t$, where $\Delta t$ denotes the time step.

The spatial discretization is based on a conservative finite-volume formulation for the balance law~\eqref{eq:hyperbolic}. In particular, the interface fluxes are computed using a second-order piecewise linear reconstruction with MinMod slope limiting, and the resulting semi-discrete operator is integrated in time by the IMEX scheme.

The IMEX RK2 method is characterized by the parameters
\begin{equation}
\gamma = 1 - \dfrac{1}{\sqrt{2}},
\qquad
c = \dfrac{1}{2\gamma},
\end{equation}
and by the explicit and implicit Butcher tableaux
\[
\begin{array}{c|cc}
0 & 0 & 0 \\
c & c & 0 \\
\hline
& 1-\gamma & \gamma
\end{array}
\qquad
\begin{array}{c|cc}
\gamma & \gamma & 0 \\
1-\gamma & 1-\gamma & \gamma \\
\hline
& 1-\gamma & \gamma
\end{array}.
\]

Let $\bar{U}_i^n$ denote the numerical approximation of the cell average over $I_i$ at grid point $x_i$ and time $t^n$. The first stage is given by
\begin{equation}
\bar{U}_i^{(1)} = \bar{U}_i^n - \gamma
\frac{\Delta t}{\Delta x} \bigl( F_{i+\frac12}(\bar{U}^n) - F_{i-\frac12}(\bar{U}^n) \bigr) + \gamma\Delta t\,
S\!\left(\bar{U}_i^{(1)}\right),
\end{equation}
where the flux contribution is evaluated explicitly and the source term implicitly.

Following the Boscarino reformulation, an auxiliary explicit state is introduced,
\begin{equation}
\bar{U}_i^{(e)} = \left(1-\frac{c}{\gamma}\right)U_i^n + \frac{c}{\gamma}\bar{U}_i^{(1)},
\end{equation}
which is used to evaluate the explicit flux in the second stage.

The second stage reads
\begin{equation}
\bar{U}_i^{n+1} = \bar{U}_i^{\ast} - \gamma \frac{\Delta t}{\Delta x}
\bigl( F_{i+\frac12}(\bar{U}^{(e)}) - F_{i-\frac12}(\bar{U}^{(e)}) \bigr)
+ \gamma\Delta t\, S\!\left(\bar{U}_i^{n+1}\right),
\end{equation}
where
\begin{equation}
\bar{U}_i^{\ast} = \left(1-\frac{1-\gamma}{\gamma}\right)\bar{U}_i^n
+ \frac{1-\gamma}{\gamma}\bar{U}_i^{(1)}.
\end{equation}

Since the numerical solution at time $t^{n+1}$ coincides with the final implicit stage, the scheme is stiffly accurate. Consequently, it remains stable in the limit of vanishing relaxation parameters and is consistent with the equilibrium limit, while preserving the conservative structure of the balance law.
\subsubsection{Two-Dimensional IMEX Runge--Kutta Scheme}
Let the two-dimensional spatial domain be
$\Omega=[a,b]\times[c,d]\subset\mathbb{R}^2$,
which is discretized into $N_x\times N_y$ uniform Cartesian cells
\[
I_{i,j}
=
[x_{i-\frac12},x_{i+\frac12}]
\times
[y_{j-\frac12},y_{j+\frac12}],
\qquad
i=1,\ldots,N_x,\quad
j=1,\ldots,N_y,
\]
with cell sizes $\Delta x=\dfrac{b-a}{N_x}$,
and $\Delta y=\dfrac{d-c}{N_y}$, and cell centers $x_i=\dfrac{1}{2}\left(x_{i-\frac12}+x_{i+\frac12}\right)$,
and $y_j=\dfrac{1}{2}\left(y_{j-\frac12}+y_{j+\frac12}\right)$.

The time interval $[0,T]$ is discretized as
$t^n=n\Delta t$, where $\Delta t$ denotes the time step. Let $\bar{U}_{i,j}^n$ denote the
cell average of $U$ over $I_{i,j}$ at time $t^n$.

The two-dimensional finite-volume spatial operator is defined by
\begin{equation}
\mathcal{L}_{i,j}(U)
=
-\dfrac{
F_{1,i+\frac12,j}(U)-F_{1,i-\frac12,j}(U)
}{\Delta x}
-\dfrac{
F_{2,i,j+\frac12}(U)-F_{2,i,j-\frac12}(U)
}{\Delta y},
\label{eq:2d_spatial_operator}
\end{equation}
where $F_1$ and $F_2$ denote the numerical fluxes in the $x$- and $y$-directions, respectively. The same second-order piecewise linear reconstruction with MinMod slope limiting is employed in both spatial directions.

The same IMEX RK2 coefficients
\begin{equation}
\gamma=1-\frac{1}{\sqrt{2}},
\qquad
c=\frac{1}{2\gamma},
\label{eq:2d_imex_coefficients}
\end{equation}
as in the one-dimensional formulation are used. The first stage is given by
\begin{equation}
\bar{U}_{i,j}^{(1)}
=
\bar{U}_{i,j}^{n}
+
c\Delta t\,\mathcal{L}_{i,j}(\bar{U}^{n})
+
\gamma\Delta t\,S\left(\bar{U}_{i,j}^{(1)}\right).
\label{eq:2d_first_stage}
\end{equation}

Following the Boscarino reformulation, the auxiliary state
\begin{equation}
\bar{U}_{i,j}^{(e)}
=
\left(1-\frac{c}{\gamma}\right)\bar{U}_{i,j}^{n}
+
\frac{c}{\gamma}\bar{U}_{i,j}^{(1)},
\label{eq:2d_explicit_state}
\end{equation}
is introduced for the explicit flux evaluation in the second stage. We
also define
\begin{equation}
\bar{U}_{i,j}^{*}
=
\left(1-\frac{1-\gamma}{\gamma}\right)\bar{U}_{i,j}^{n}
+
\frac{1-\gamma}{\gamma}\bar{U}_{i,j}^{(1)}.
\label{eq:2d_star_state}
\end{equation}
The second stage is then written as
\begin{equation}
\bar{U}_{i,j}^{n+1}
=
\bar{U}_{i,j}^{*}
+
\gamma\Delta t\,\mathcal{L}_{i,j}(\bar{U}^{(e)})
+
\gamma\Delta t\,S\left(\bar{U}_{i,j}^{n+1}\right).
\label{eq:2d_second_stage}
\end{equation}

For the relaxation models considered in this work, the implicit source update is solved locally in each computational cell. Thus, the spatial coupling arises only through the numerical fluxes, while the relaxation term is treated implicitly at the Runge-Kutta stages.

The two-dimensional formulation therefore differs from the one-dimensional scheme only through the spatial operator $\mathcal{L}_{i,j}$, which accounts for the numerical fluxes through the four faces of each Cartesian control volume. The convective terms in both spatial directions are treated explicitly, whereas the relaxation source term is treated implicitly. Since the final stage coincides with the numerical solution at $t^{n+1}$, the IMEX RK2 formulation retains the stiffly accurate property and is suitable for the stiff-relaxation regime.

\subsection{Numerical Experiments}
The numerical performance of the proposed relaxation models is investigated using the CS-EBT2, UCS2 and IMEX RK2 schemes across a series of benchmark problems. All simulations are performed on uniform grids, with the CS-EBT2, UCS2 and IMEX RK2 solutions evaluated either against the exact solution when available or against a high-resolution reference solution obtained with the IMEX RK2 scheme on a finer mesh. The time step is determined according to
\begin{equation}
\Delta t^n = \mathrm{CFL}\,\frac{\Delta x}{\Lambda_{\max}^n},
\end{equation}
where $\Lambda_{\max}^n$ denotes the maximum eigenvalue magnitude of the flux Jacobian $\partial F/\partial U$ at time $t^n$, and the CFL number is chosen to ensure numerical stability for each scheme.

The preservation of integral invariance is assessed for all conserved variables using the discrete $L^1$-norm, defined as $\|U - U_0\|_{L^1}$. Grid refinement studies are performed to examine whether this property is retained as the spatial resolution increases. Minimal variation in the $L^1$-norm across successive grid resolutions indicates that the integral invariance property is effectively maintained by the numerical schemes.

The results for each benchmark problem are presented in the accompanying figures and tables. Figures illustrate the evolution of solution variables under different relaxation parameters, and tables report the corresponding $L^1$-errors to quantify the preservation of the integral invariance property. These results demonstrate the behavior of the newly proposed relaxation models when discretized with numerical schemes, confirming that the models retain the expected conservation properties across the considered relaxation regimes and initial conditions.

For the two-dimensional benchmark problems, the proposed relaxation models are discretized on uniform Cartesian meshes using the CS-EBT2, UCS2, and IMEX RK2 schemes. The convective fluxes are evaluated in both the $x$- and $y$-directions, while the relaxation source terms are treated
according to the corresponding time discretization schemes. The time step is selected according to the multidimensional CFL condition. The numerical solutions are compared with the exact solution when available or with a
high-resolution IMEX RK2 reference solution otherwise. The preservation of integral invariance is assessed using the corresponding $L^1$-norm, and grid refinement studies are performed to examine its behavior for increasing spatial resolution.

\subsubsection{Jin-Xin Relaxation Model}
\begin{figure}[!ht]
     \centering
     \begin{minipage}[b]{0.48\linewidth}
         \centering
         \includegraphics[width=\linewidth]{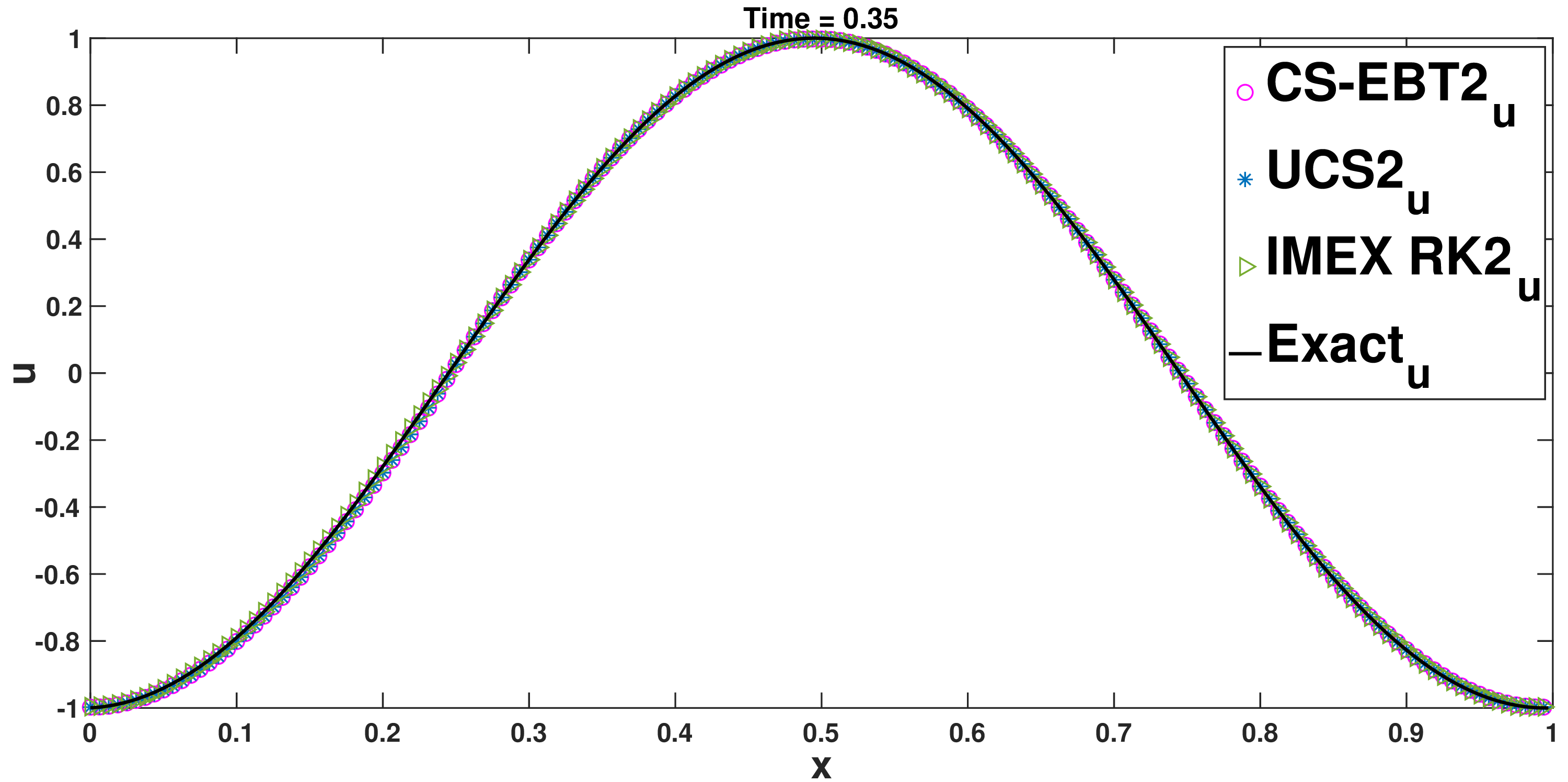}
     \end{minipage}
     \hfill
     \begin{minipage}[b]{0.48\linewidth}
         \centering
         \includegraphics[width=\linewidth]{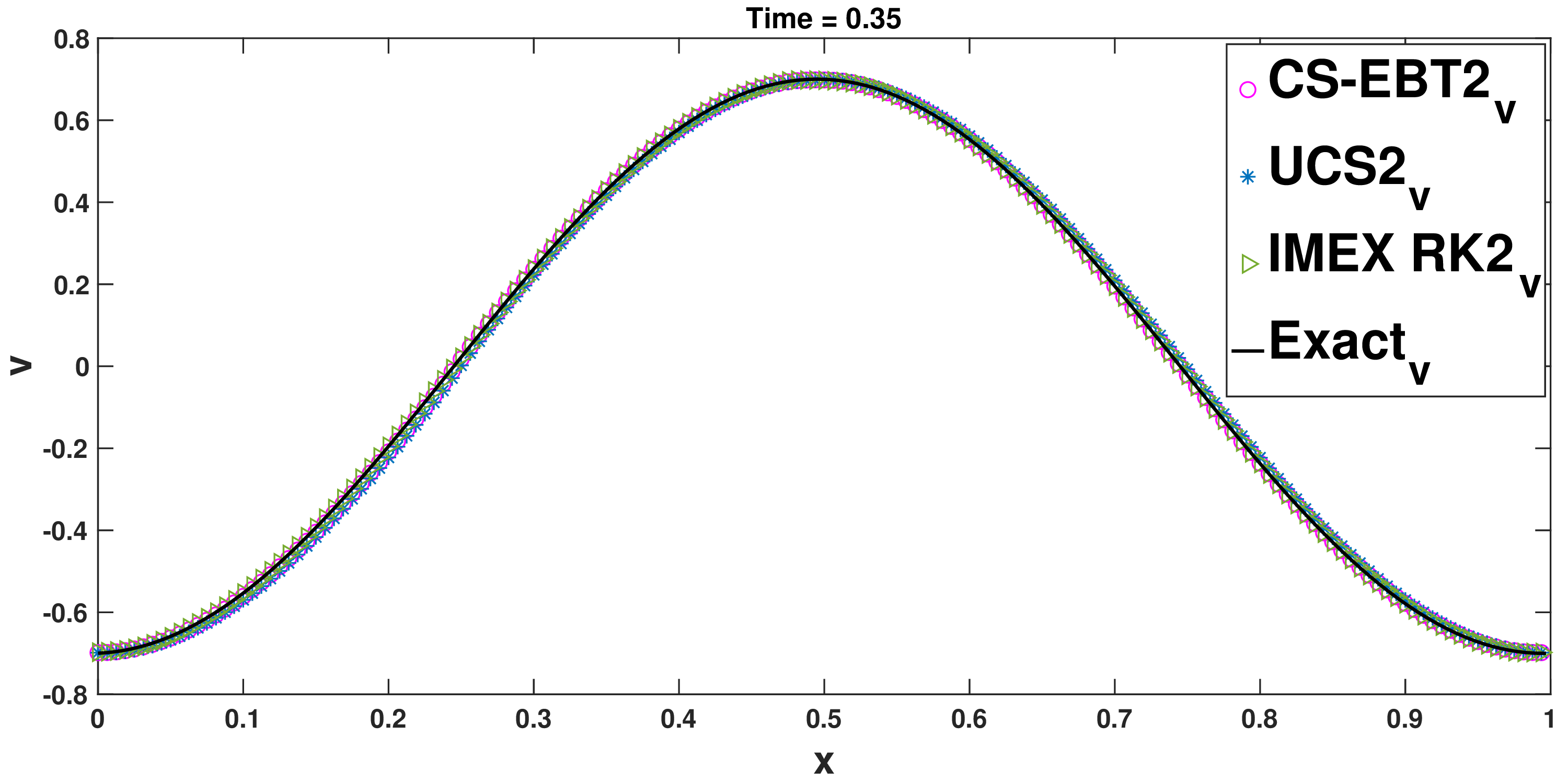}
     \end{minipage}
     \caption{Alternative Jin–Xin model with smooth (well-prepared) initial data. Comparison between the numerical and exact solutions for $u$ (left) and $v$ (right), computed with CFL $=0.9$, $\tau = 10^{-10}$, and $N = 320$.}
        \label{fig:example1a}
\end{figure}

We consider the alternative Jin-Xin relaxation model \eqref{modify:xinJin} on the spatial domain $[0,1]$, subject to suitable boundary conditions. Both smooth and non-smooth initial data are examined.\\
\noindent\textbf{Smooth case (well-prepared data).}
The initial conditions are prescribed by
\begin{equation}\label{smoothdata:XinJin}
\begin{aligned}
u(x,0) &= \sin(2\pi x), \\
v(x,0) &= au(x,0),
\end{aligned}
\end{equation}
with the parameter $a = 0.7$. The numerical computations are performed on a uniform grid with $N = 320$ points and a CFL number of $0.9$ with periodic boundary conditions. The stiff relaxation regime is enforced by setting the relaxation parameter $\tau = 10^{-10}$. The solution is advanced up to the final time $T = 0.35$.
Figure~\ref{fig:example1a} presents a comparison of the numerical solutions obtained using the CS-EBT2, UCS2 and IMEX RK2 schemes with the corresponding exact solution.
\par Tables~\ref{Tabu:xinjinerror} and~\ref{Tabv:xinjinerror} report the $L^1$-errors and the corresponding convergence rates obtained with the three numerical schemes for the variables $u$ and $v$, respectively, for the modified Jin--Xin relaxation model~\eqref{modify:xinJin}. The simulations are carried out using the smooth initial data~\eqref{smoothdata:XinJin}, with $\tau=10^{-10}$ and a CFL number of $0.9$.
\renewcommand{\arraystretch}{1.2}{
\begin{table}[ht!]
 \caption{Error and order of convergence for alternative Jin-Xin model with smooth initial condition \eqref{smoothdata:XinJin} for $u$ with $\tau = 10^{-10}$}
\centering
\begin{tabular}{*{7}{c}}
\toprule
\multirow{2}{*}{N} & 
\multicolumn{2}{c}{CS-EBT2} &
\multicolumn{2}{c}{UCS2} & 
\multicolumn{2}{c}{IMEX RK2} \\
\cmidrule(lr){2-3}
\cmidrule(lr){4-5}
\cmidrule(lr){6-7}
 & $L^{1}$-Error & Order 
 & $L^{1}$-Error & Order 
 & $L^{1}$-Error & Order \\
\midrule
32  &  1.1297e-03 & - &  1.0284e-03  & - &  2.3313e-02  & -\\
64 & 2.8779e-04  & 1.97  &  2.6116e-04 & 1.98  &  8.1999e-03 & 1.51\\
128 & 7.3145e-05 & 1.98  &  6.6603e-05 & 1.97  & 2.4093e-03 & 1.77\\
256 & 1.8458e-05 & 1.99 & 1.6853e-05  & 1.98  &  6.3227e-04 & 1.93\\
512 & 4.6325e-06 & 1.99  &  4.2377e-06  & 1.99  & 1.5796e-04 & 2.00\\
1024 & 1.1594e-06 & 2.00 & 1.0618e-06  &  2.00  & 3.7653e-05 & 2.07 \\  
\bottomrule
\end{tabular}
\label{Tabu:xinjinerror}
\end{table}}
 \renewcommand{\arraystretch}{1.2}{
\begin{table}[ht!]
 \caption{Error and order of convergence for alternative Jin-Xin model with smooth initial condition \eqref{smoothdata:XinJin} for $v$ with $\tau = 10^{-10}$}
\centering
\begin{tabular}{*{7}{c}}
\toprule
\multirow{2}{*}{N} & 
\multicolumn{2}{c}{CS-EBT2} &
\multicolumn{2}{c}{UCS2} & 
\multicolumn{2}{c}{IMEX RK2} \\
\cmidrule(lr){2-3}
\cmidrule(lr){4-5}
\cmidrule(lr){6-7}
 & $L^{1}$-Error & Order 
 & $L^{1}$-Error & Order 
 & $L^{1}$-Error & Order \\
\midrule
32 & 7.9081e-04    &  -  & 3.4926e-03  &  -   & 1.6319e-02  &   -\\
64  & 2.0145e-04   &  1.97   &  6.116e-04 &  2.51  &  5.7399e-03 &  1.51\\
128 & 5.1202e-05   &  1.98  &   1.7225e-04   & 1.83  & 1.6865e-03  & 1.77\\
256 & 1.2921e-05  &  1.99  &  5.2202e-05 & 1.72  & 4.4259e-04 & 1.93\\
512 & 3.2428e-06  &  1.99  &  1.755e-05 & 1.57   &  1.1057e-04 & 2.00\\
1024 & 8.1160e-07 & 2.00  & 4.1121e-06 & 2.09 &  2.6357e-05 & 2.07\\
\bottomrule
\end{tabular}
\label{Tabv:xinjinerror}
\end{table}}
\\
\noindent{\textbf{Non-Smooth case:}} The discontinuous profile is given by
\begin{eqnarray}\label{xinjin:discontinuos}
\begin{aligned}
&u(x,0)=
\begin{cases}
 2, \, \text{if}\,\,\,\, 0.25<x<0.5,\\
1, \, \text{otherwise},
\end{cases} \\
&v(x,0) = au(x,0),
\end{aligned}
\end{eqnarray}
with $a=0.7$. To assess the performance of the numerical schemes across different relaxation regimes, numerical simulations are conducted on the domain $[0,1]$ using $N=400$ uniformly spaced grid points. Transmissive boundary conditions are imposed for the CS-EBT2, UCS2, and IMEX RK2 schemes. In each case, the computed numerical solutions are compared with the corresponding exact solutions.

\begin{figure}[!ht]
     \centering
     \begin{minipage}[b]{0.48\linewidth}
         \centering
         \includegraphics[width=\linewidth]{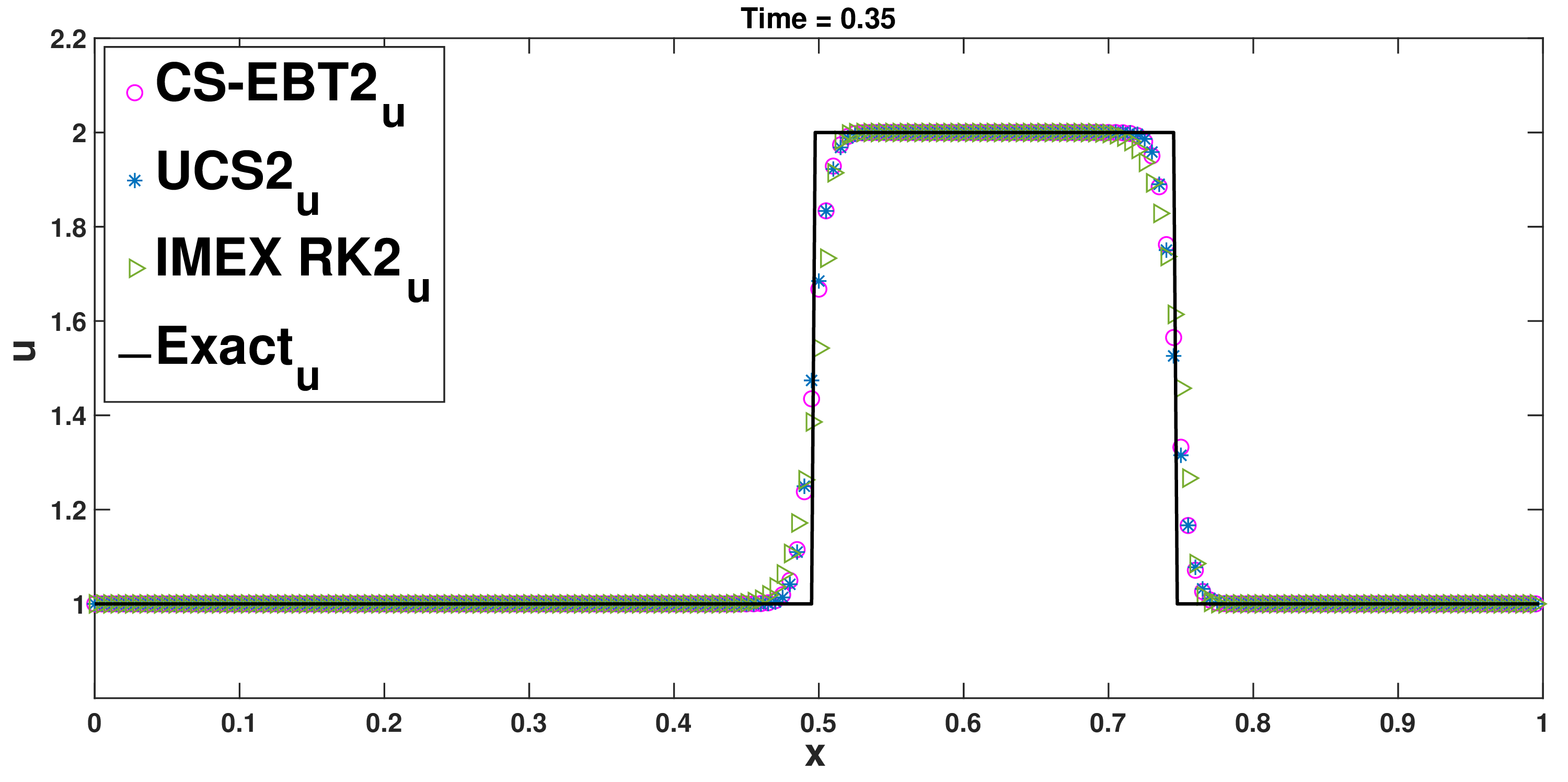}
     \end{minipage}
     \hfill
     \begin{minipage}[b]{0.48\linewidth}
         \centering
         \includegraphics[width=\linewidth]{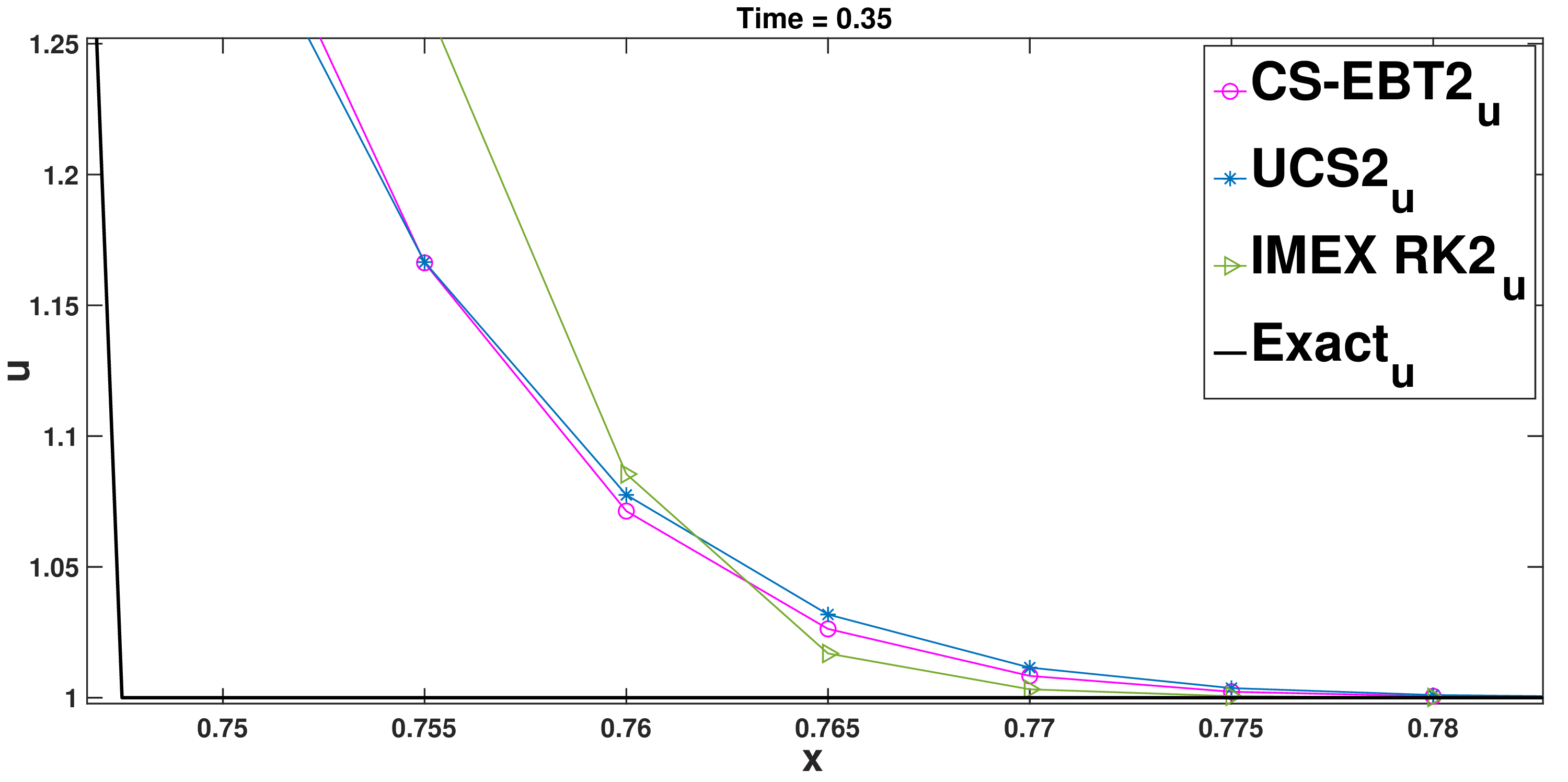}
     \end{minipage}
     \caption{Alternative Jin-Xin model with non-smooth case: comparison between numerical solution $u$(left) and $u$ zoomed version(right) and exact solution with CFL $0.9$, $\tau=10^{-10}$ and $N=400$.}
        \label{fig:example2a}
\end{figure}

\begin{figure}[!ht]
     \centering
     \begin{minipage}[b]{0.48\linewidth}
         \centering
         \includegraphics[width=\linewidth]{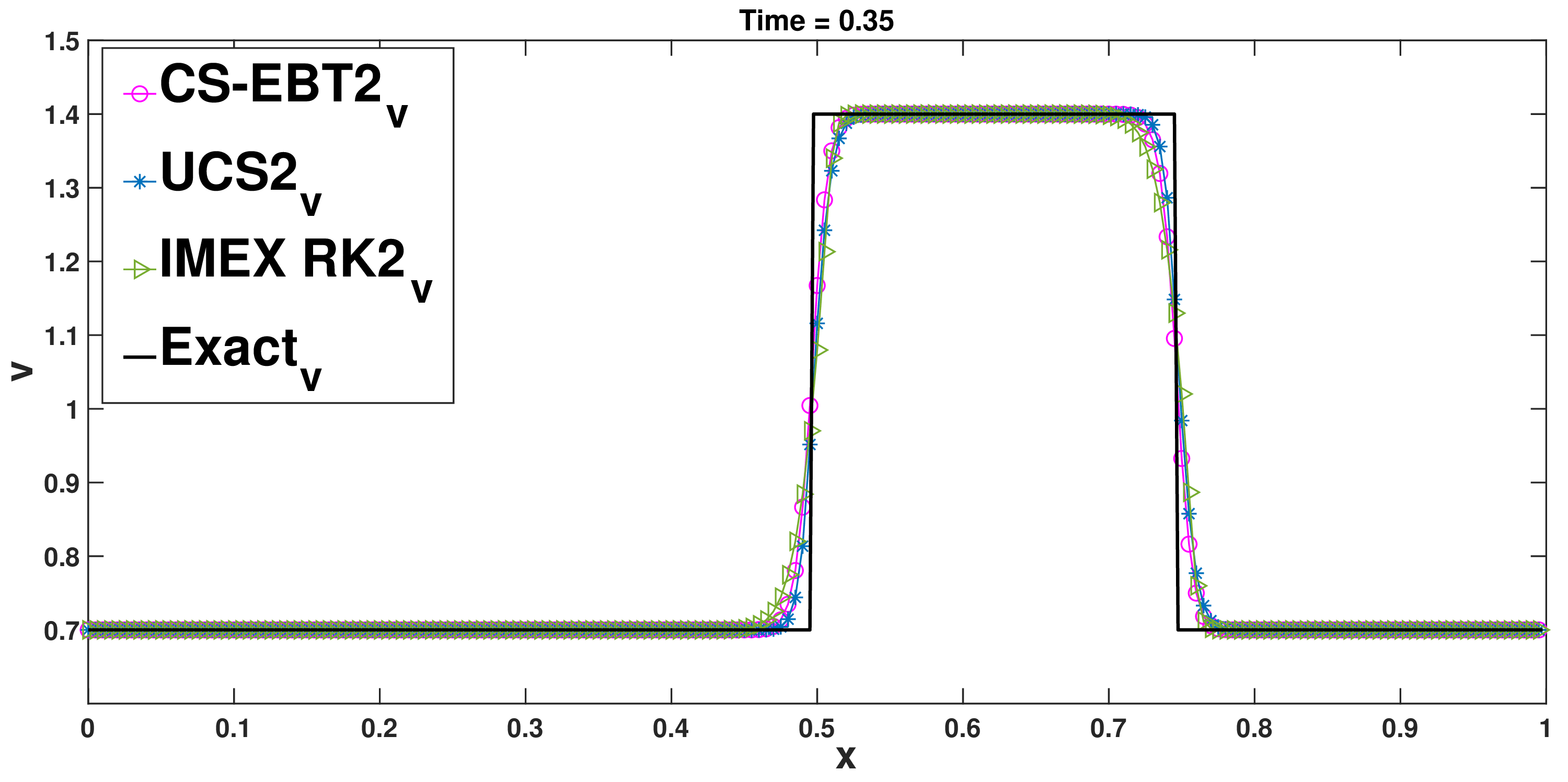}
     \end{minipage}
     \hfill
     \begin{minipage}[b]{0.48\linewidth}
         \centering
         \includegraphics[width=\linewidth]{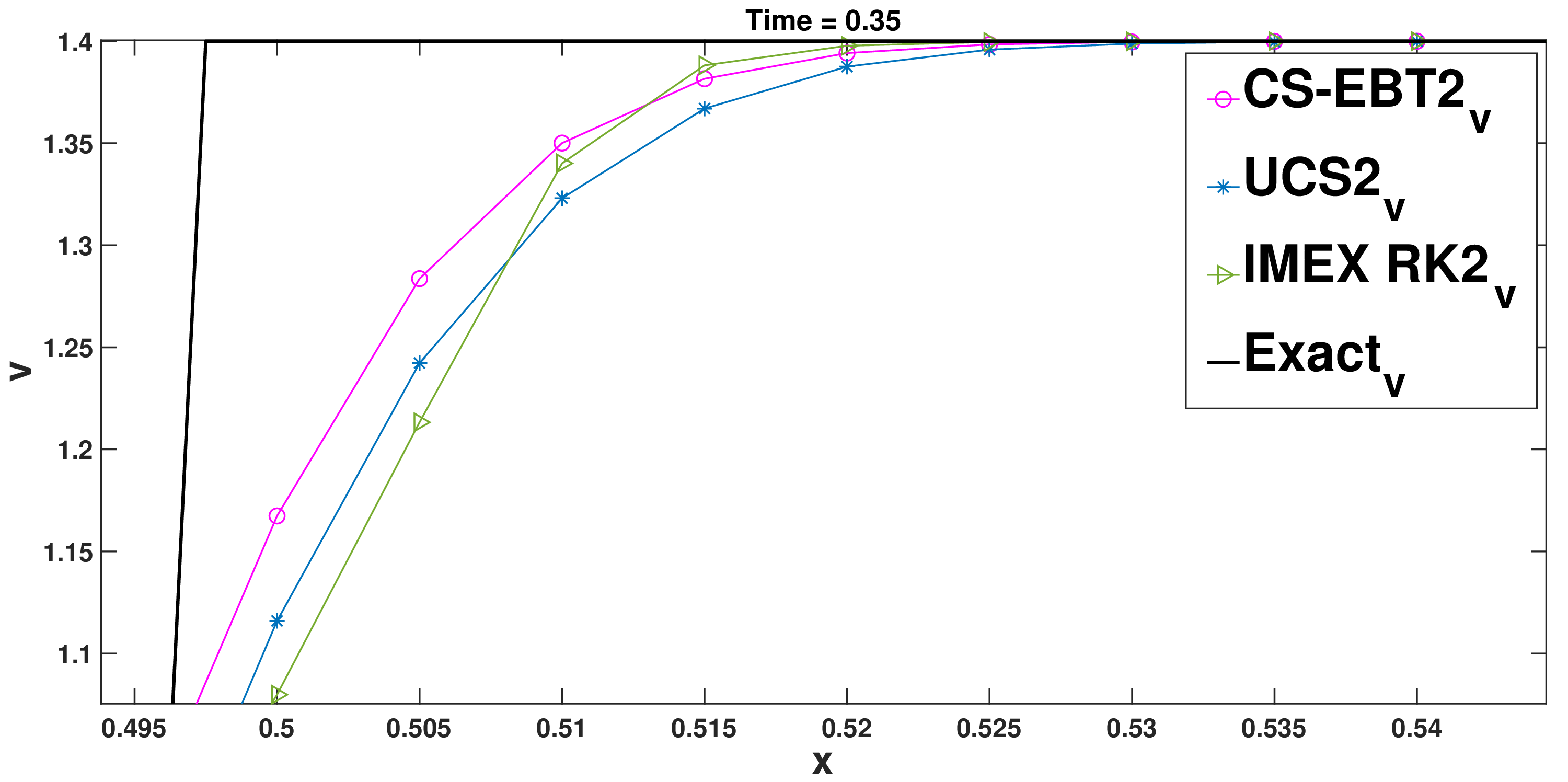}
     \end{minipage}
     \caption{Alternative Jin-Xin model with non-smooth case: comparison between numerical solution $v$(left) and $v$ zoomed version(right) and exact solution with CFL $0.9$, $\tau=10^{-10}$ and $N=400$.}
        \label{fig:example21a}
\end{figure}

\begin{figure}[!ht]
     \centering
     \begin{minipage}[b]{0.48\linewidth}
         \centering
         \includegraphics[width=\linewidth]{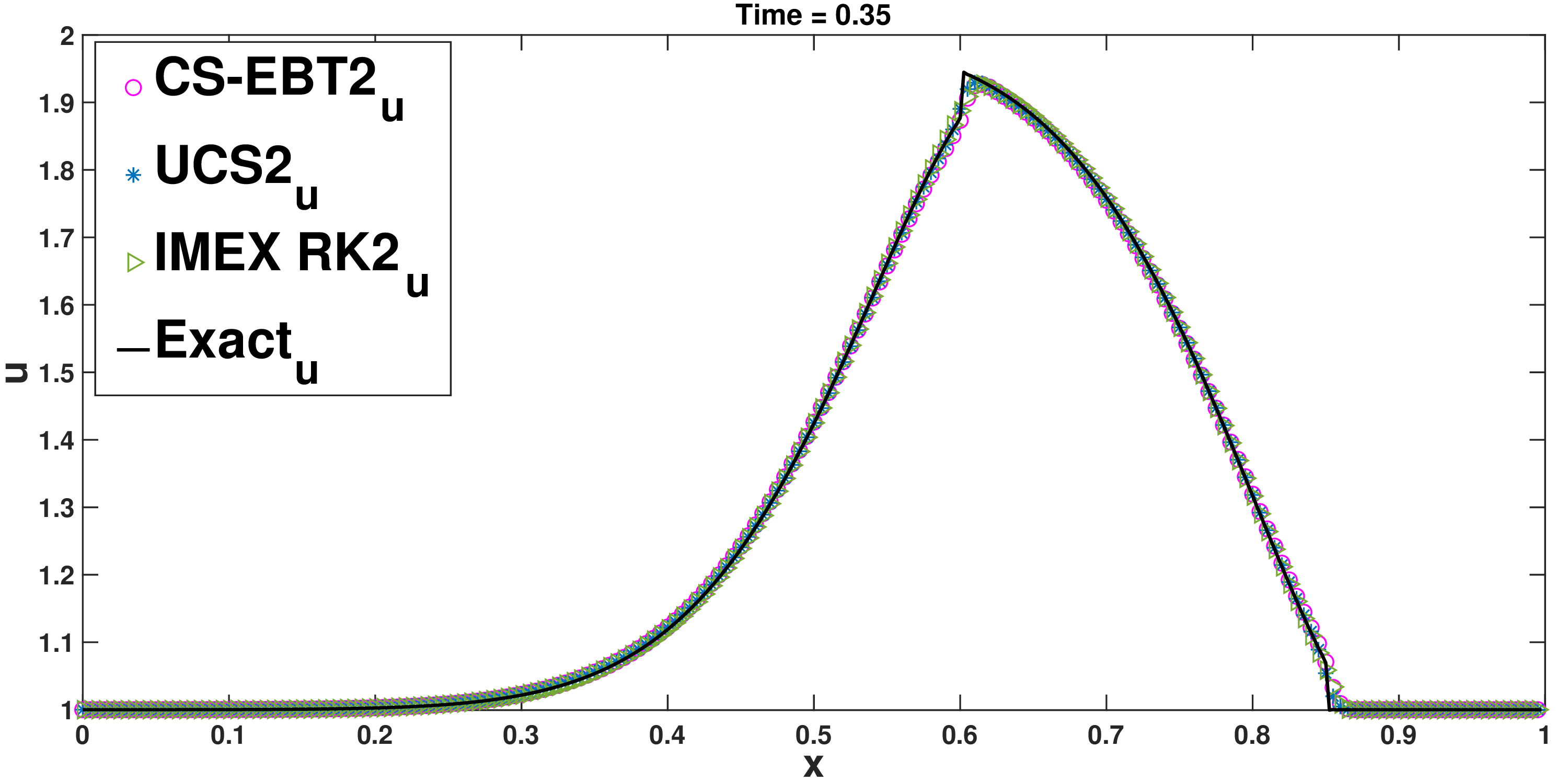}
     \end{minipage}
     \hfill
     \begin{minipage}[b]{0.48\linewidth}
         \centering
         \includegraphics[width=\linewidth]{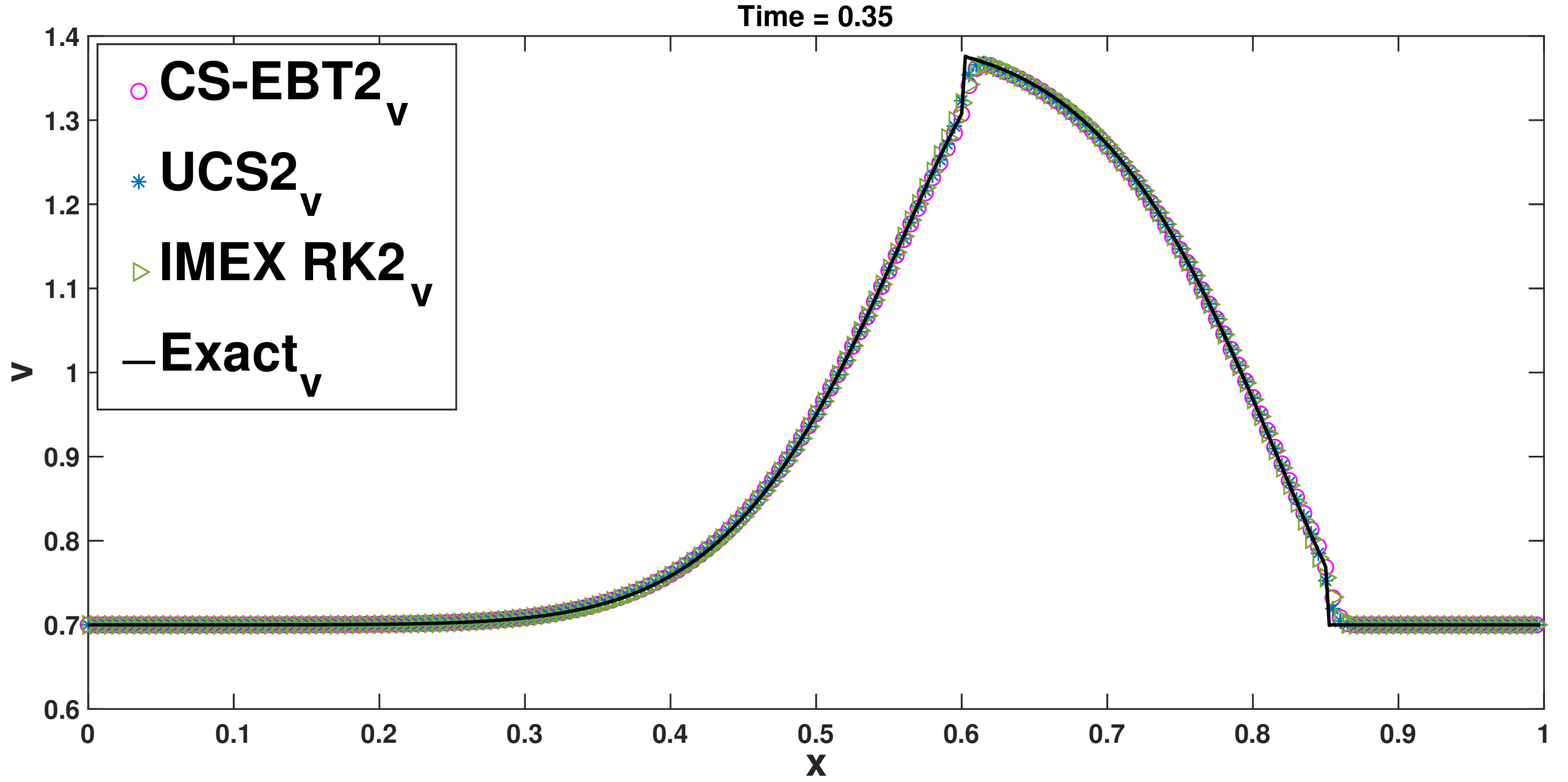}
     \end{minipage}
     \caption{Alternative Jin-Xin model with non-smooth case: comparison between numerical solution $u$(left) and $v$(right) and exact solution with CFL $0.9$, $\tau=0.02$ and $N=400$.}
        \label{fig:example2b}
\end{figure}

\begin{figure}[!ht]
     \centering
     \begin{minipage}[b]{0.48\linewidth}
         \centering
         \includegraphics[width=\linewidth]{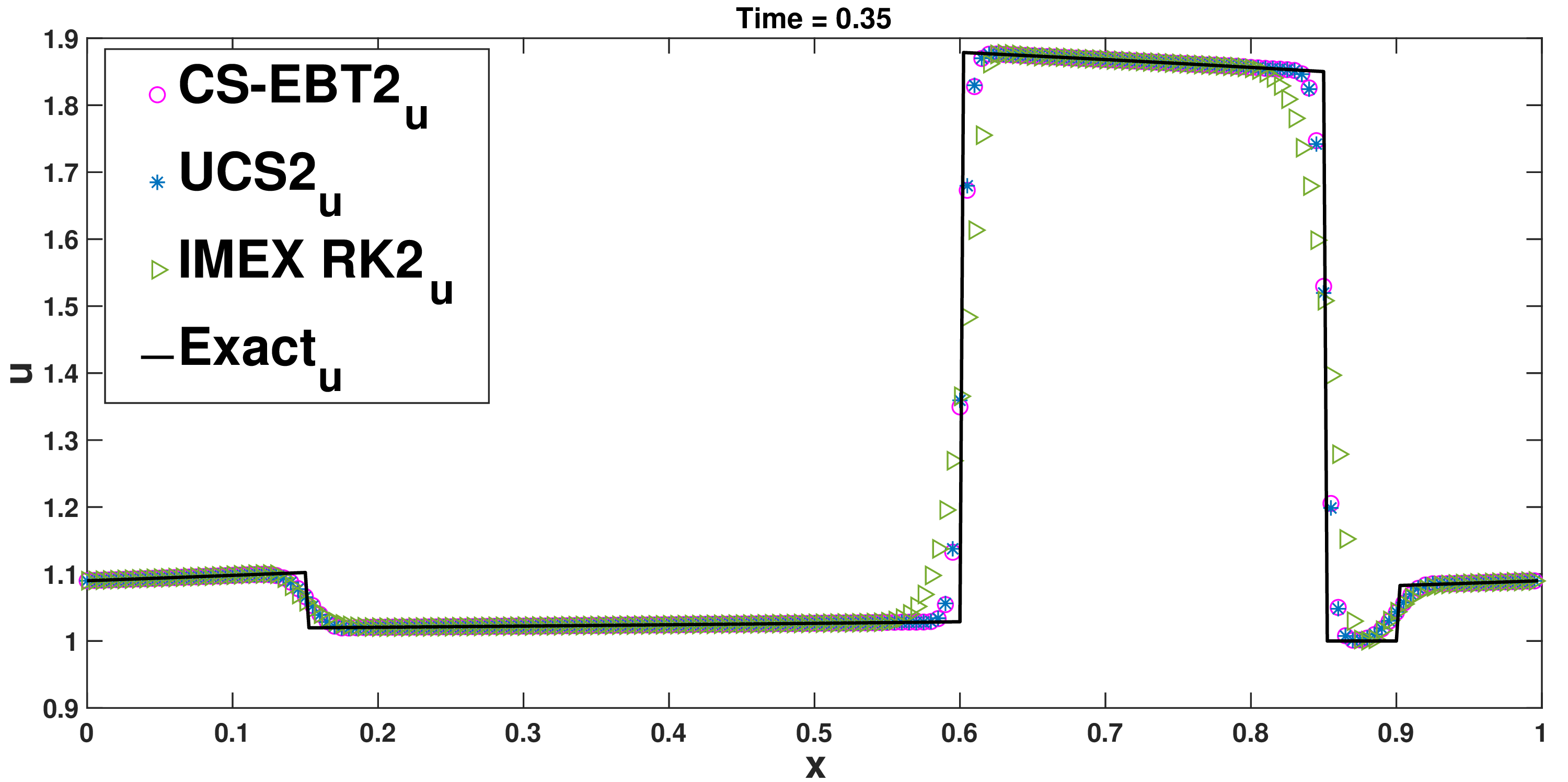}
     \end{minipage}
     \hfill
     \begin{minipage}[b]{0.48\linewidth}
         \centering
         \includegraphics[width=\linewidth]{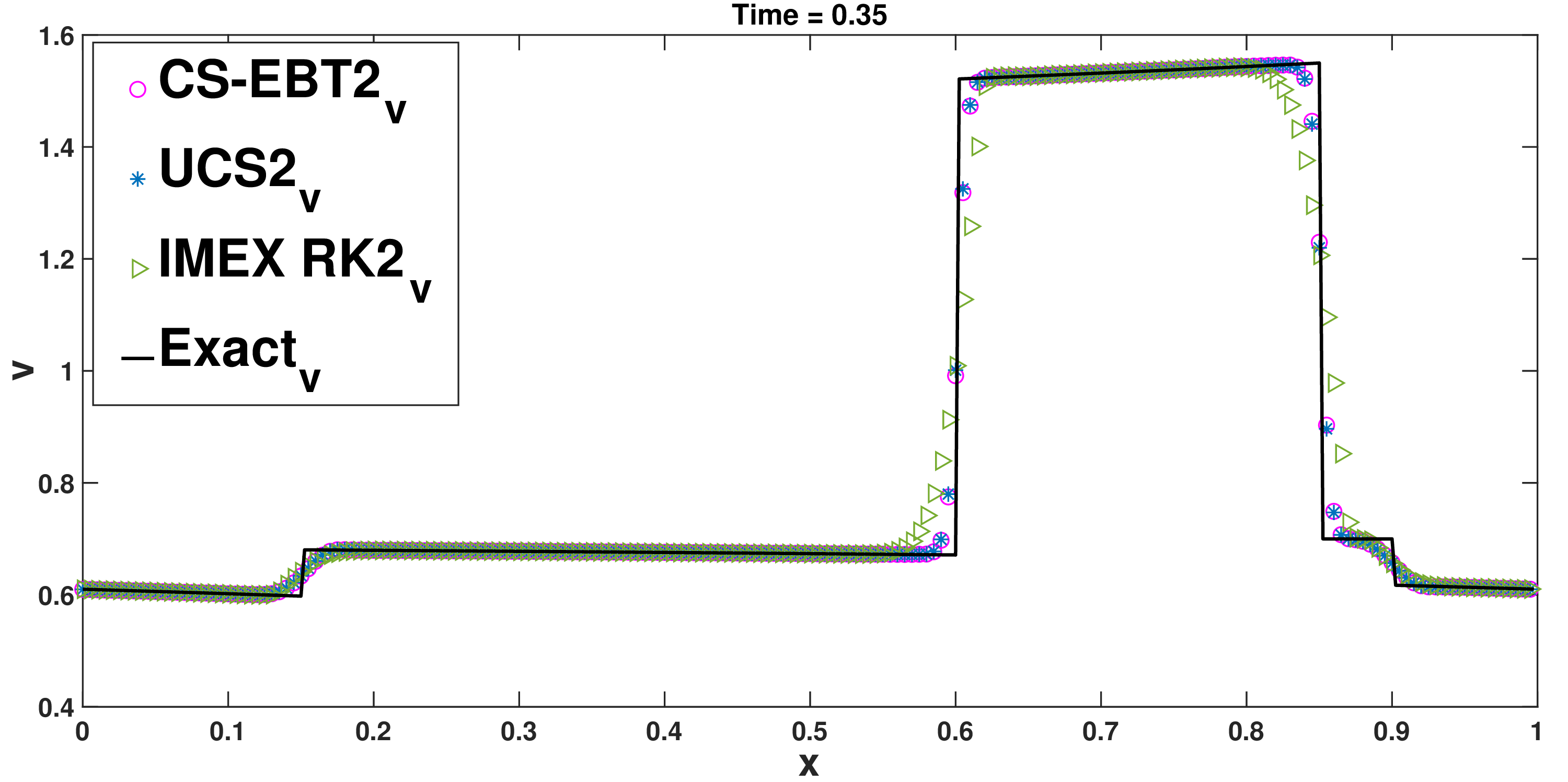}
     \end{minipage}
     \caption{Alternative Jin-Xin model with non-smooth case: comparison between numerical solution $u$(left) and $v$(right) and exact solution with CFL $0.9$, $\tau=1$ and $N=400$.}
        \label{fig:example2c}
\end{figure}

Figures~\ref{fig:example2a}-\ref{fig:example2c} show the numerical results at $T=0.35$ for non-smooth initial data with $\tau=10^{-10}$, $0.02$, and $1$, respectively. Table~\ref{Tabxinjin} provides a numerical verification of the integral invariance property for the CS-EBT2, UCS2, and IMEX RK2 schemes. The results demonstrate that the integral invariant is preserved by all three numerical methods.
 \renewcommand{\arraystretch}{1.2}{
\begin{table}[ht!]
\centering
\caption{Integral invariance for alternative Jin-Xin model with initial condition
\eqref{xinjin:discontinuos} for $u + \tau v$ with $\tau = 10^{-10}$ and CFL $0.9$}
\begin{tabular}{*{5}{c}}
\toprule
\multirow{2}{*}{$N$}
& \multicolumn{1}{c}{CS-EBT2}
& \multicolumn{1}{c}{UCS2}
& \multicolumn{1}{c}{IMEX RK2} \\
\cmidrule(lr){2-2} 
\cmidrule(lr){3-3}
\cmidrule(lr){4-4}
& Error $u + \tau v$ & Error $u + \tau v$ & Error $u + \tau v$ \\
\midrule

40 & 3.99999995e-01 & 4.00000024e-01  & 3.42169701e-06  \\

80 & 4.00000000e-01 & 4.00000000e-01  & 3.66917607e-12 \\

160 & 4.00000000e-01 & 4.00000000e-01  & 1.11022302e-16 \\

320 & 4.00000000e-01 & 4.00000000e-01 & 1.11022302e-16 \\

640 & 4.00000000e-01 & 4.00000000e-01  & 1.11022302e-16 \\

1280 & 4.00000000e-01 & 4.00000000e-01 & 1.11022302e-16\\
\bottomrule
\end{tabular}
\label{Tabxinjin}
\end{table}}

\subsubsection{Shallow Water Equations}
\begin{figure}[ht]
    \begin{minipage}[b]{0.48\linewidth}
        \includegraphics[width=\linewidth]{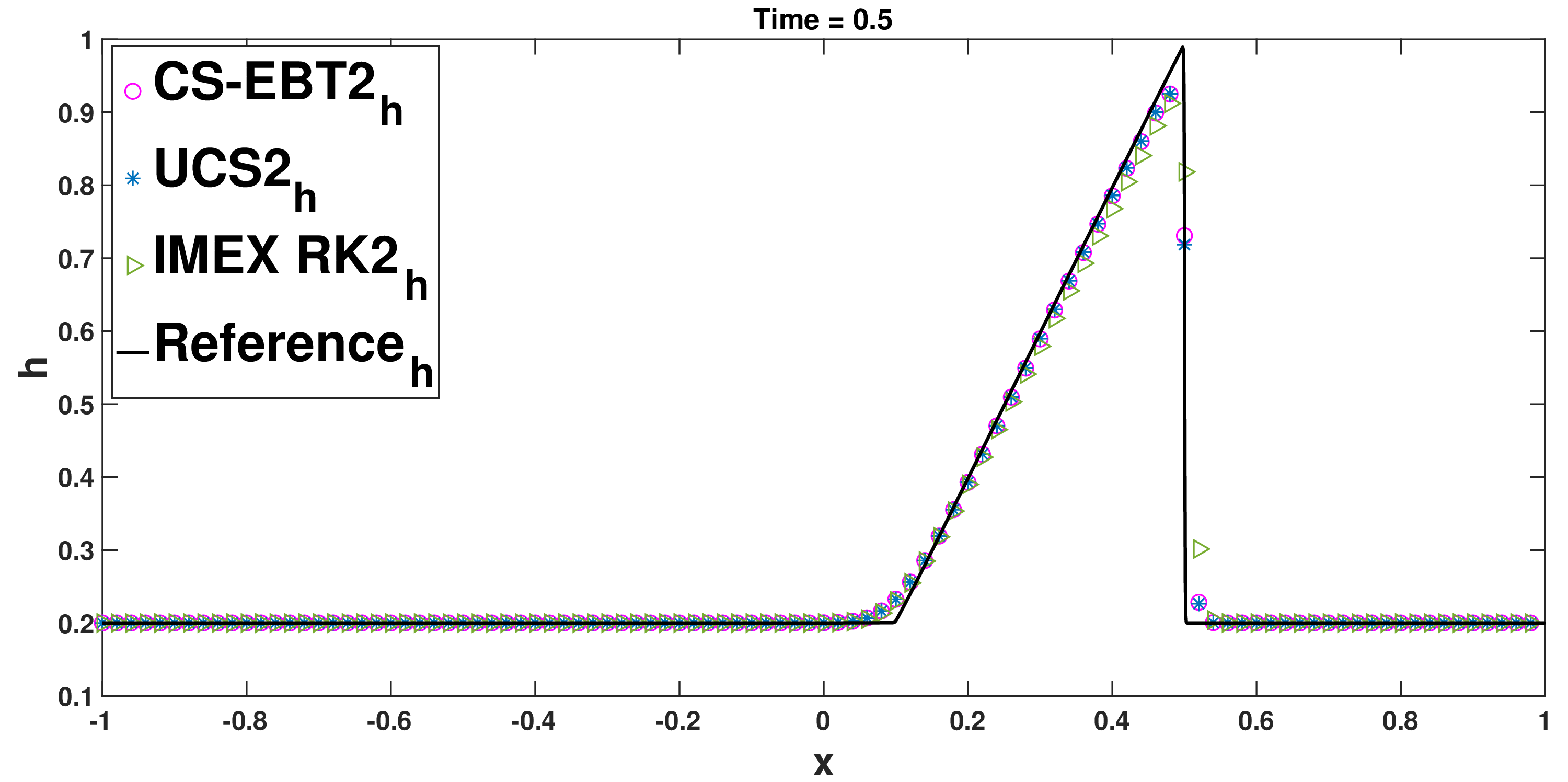}
    \end{minipage}
    \hfill
    \begin{minipage}[b]{0.48\linewidth}
        \includegraphics[width=\linewidth]{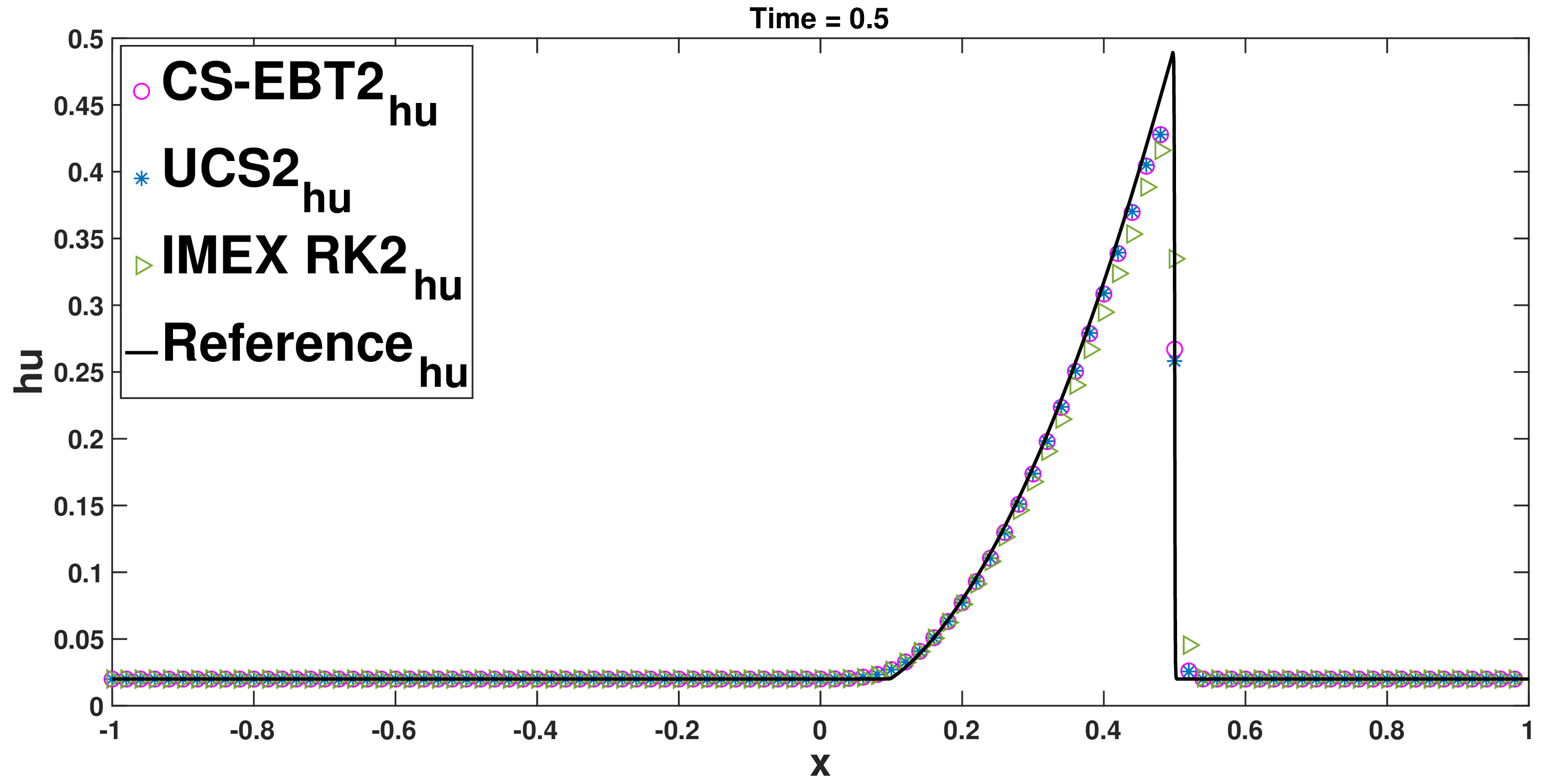}
    \end{minipage}
       \caption{Alternative shallow water model with non-smooth case: comparison between numerical solution $h$(left), $hu$(right) and the reference solution (IMEX RK2) with $\tau = 10^{-10}$, CFL $0.9$ and $N=200$.}
    \label{shallowWater}
\end{figure}

The alternative relaxation model for the shallow water equations \eqref{modify:shallow}, together with the corresponding smooth and non-smooth initial conditions, is given by\\
\textbf{Smooth case:}
\begin{eqnarray}\label{shallow:continuous}
\begin{cases}
\begin{aligned}
&h(x,0)= 1+0.2 \sin(8\pi x),\\
&hu(x,0) = \dfrac{1}{2}\;h^2(x,0).
\end{aligned}
\end{cases}
\end{eqnarray}
The numerical experiments are performed on the spatial domain $[-1,1]$ with a CFL number of $0.9$. To investigate the stiff relaxation regime, the relaxation parameter is fixed at $\tau = 10^{-10}$. All simulations are carried out up to the final time $T = 0.35$ subject to periodic boundary conditions.
The numerical results reported in Tables~\ref{Tabh:shallowerror} and~\ref{Tabhu:shallowerror} demonstrate the accuracy of the three numerical schemes for the modified shallow water relaxation model~\eqref{modify:shallow}. The tables list the $L^1$-errors and the corresponding convergence orders for the variables $h$ and $hu$, respectively, computed using the smooth initial condition~\eqref{shallow:continuous} with $\tau = 10^{-10}$ and a CFL number of $0.9$.\\
\renewcommand{\arraystretch}{1.2}{
\begin{table}[ht!]
 \caption{Error and order of convergence for modify shallow water model with smooth initial condition \eqref{shallow:continuous} for $h$ with $\tau = 10^{-10}$}
\centering
\begin{tabular}{*{7}{c}}
\toprule
\multirow{2}{*}{N} & 
\multicolumn{2}{c}{CS-EBT2} &
\multicolumn{2}{c}{UCS2} & 
\multicolumn{2}{c}{IMEX RK2} \\
\cmidrule(lr){2-3}
\cmidrule(lr){4-5}
\cmidrule(lr){6-7}
 & $L^{1}$-Error & Order 
 & $L^{1}$-Error & Order 
 & $L^{1}$-Error & Order \\
\midrule
64 & 1.6823e-03 &  -  & 1.6282e-03 & - & 2.0274e-03  & -\\

128 & 4.8501e-04  & 1.79 & 4.7006e-04 & 1.79  & 5.7192e-04 & 1.83\\

256  & 1.2023e-04  & 2.01  & 1.1671e-04 &  2.01 & 1.4209e-04 & 2.01\\

512  & 3.1253e-05 & 1.94  &  2.9062e-05 & 2.01 & 3.2779e-05 &  2.12\\

1024 & 8.3624e-06 & 1.90  & 6.3715e-06  &  2.19  & 7.4911e-06  & 2.13  \\ 
\bottomrule
\end{tabular}
\label{Tabh:shallowerror}
\end{table}}
  \renewcommand{\arraystretch}{1.2}{
\begin{table}[ht!]
 \caption{Error and order of convergence for modify shallow water model with smooth initial condition \eqref{shallow:continuous} for $hu$ with $\tau = 10^{-10}$}
\centering
\begin{tabular}{*{7}{c}}
\toprule
\multirow{2}{*}{N} & 
\multicolumn{2}{c}{CS-EBT2} &
\multicolumn{2}{c}{UCS2} & 
\multicolumn{2}{c}{IMEX RK2} \\
\cmidrule(lr){2-3}
\cmidrule(lr){4-5}
\cmidrule(lr){6-7}
 & $L^{1}$-Error & Order 
 & $L^{1}$-Error & Order 
 & $L^{1}$-Error & Order \\
\midrule
64  & 1.6783e-03 &  - & 1.6248e-03  & -  & 2.0212e-03  & -\\

128 & 4.8018e-04 &  1.81 &  4.6648e-04  & 1.80  & 5.6771e-04 & 1.83\\

256 & 1.1883e-04  &  2.01 & 1.1567e-04 & 2.01 &  1.3978e-04 & 2.02\\

512  & 3.0786e-05  & 1.95 & 2.8761e-05  & 2.01 &  3.2366e-05 & 2.11\\

1024 & 8.2011e-06 & 1.91 &  6.2862e-06 & 2.19  & 7.4153e-06 & 2.13
\\
\bottomrule
\end{tabular}
\label{Tabhu:shallowerror}
\end{table}}
 \renewcommand{\arraystretch}{1.2}{
\begin{table}[ht!]
 \caption{Integral invariance for alternative shallow water model with initial condition \eqref{shallow:discontinuous} for $h + \tau hu$ with $\tau = 10^{-10}$ and CFL $0.9$}
 \centering 
\begin{tabular}{*{5}{c}}
\toprule
\multirow{2}{*}{$N$}
& \multicolumn{1}{c}{CS-EBT2}
& \multicolumn{1}{c}{UCS2}
& \multicolumn{1}{c}{IMEX RK2} \\
\cmidrule(lr){2-2} 
\cmidrule(lr){3-3}
\cmidrule(lr){4-4}
& Error $h + \tau hu$ & Error $h + \tau hu$ & Error $h + \tau hu$ \\
\midrule

20 & 1.8999970e-01 & 1.89995963e-01  & 2.00000107e-02 \\

40 & 1.89992548e-01 &  1.89980593e-01  & 2.00000117e-02 \\

80 & 1.89994731e-01 & 1.89968542e-01  & 2.00000160e-02 \\

160 & 1.89996961e-01 & 1.89947282e-01  & 2.00000160e-02 \\

320 & 1.89998389e-01 & 1.89904793e-01 & 2.00000245e-02 \\

640 & 1.89999173e-01 & 1.89819442e-01  & 2.00000595e-02 \\

1280 & 1.89999582e-01 & 1.89653760e-01 & 2.00000595e-02 \\
\bottomrule
\end{tabular}
\label{Tabshallowwater}
\end{table}}
\newline
\noindent  \textbf{Non-Smooth case:}
\begin{eqnarray}\label{shallow:discontinuous}
\begin{cases}
\begin{aligned}
&h(x,0)=
\begin{cases}
 1, \, \text{if}\,\,\,\, 0<x<0.2,\\
0.2, \, \text{otherwise},
\end{cases} \\
&hu(x,0) = -\dfrac{1}{2}\;h^2(x,0),
\end{aligned}
\end{cases}
\end{eqnarray}
where $hu$ be the non-local-equilibrium initial data.
The numerical simulations are performed under transmissive boundary conditions. In the stiff relaxation regime, computations are carried out using $N = 200$ grid points, a CFL number of $0.9$, relaxation parameter $\tau = 10^{-10}$, and final time $T = 0.5$. The numerical solutions for the height $h$ and momentum $hu$, obtained using the CS-EBT2, UCS2 and IMEX RK2 schemes on the $N = 200$ grid, are compared in Figure~\ref{shallowWater} with a high-resolution reference solution computed using the IMEX RK2 scheme on a finer mesh with $N = 3200$ grid points.
Table~\ref{Tabshallowwater} provides numerical verification of the integral invariance property of the alternative shallow water model, confirming that the CS-EBT2, UCS2 and IMEX RK2 schemes preserve this property.
\begin{remark}
The use of different CFL numbers is motivated by the stability properties of the NT scheme in the presence of relaxation terms.  While the NT scheme is subject to a restrictive CFL condition for homogeneous hyperbolic systems, the stiff relaxation limit ($\tau \to 0$) enforces rapid convergence toward equilibrium, thereby permitting the use of larger CFL numbers without loss of stability. For larger values of $\tau$, the equilibrium effect becomes less pronounced, and a more restrictive CFL condition  required to maintain numerical stability and accuracy.    
\end{remark}
\subsubsection{Broadwell System}
For alternative Broadwell model~\eqref{eq:alternative_broadwell}, the initial conditions are prescribed in terms of the conservative variables $\rho$, $m$, and $z$. Periodic boundary conditions are imposed for the smooth initial data, while transmissive boundary conditions are used for the non-smooth initial data.\\
\noindent
\textbf{Smooth case:} The initial data are given by
\begin{equation}\label{broadwell:continuous}
(\rho(x,0),m(x,0),z(x,0))=
\begin{cases}
\left(\rho,\rho u,\dfrac{1}{2}\rho(1+u^2)\right), & 0\le x\le l,
\end{cases}
\end{equation}
where
\[
\rho(x,0)=1+0.3\sin\!\left(\frac{2\pi x}{l}\right), \quad
u(x,0)=0.5+0.1\sin\!\left(\frac{2\pi x}{l}\right),
\]
with $l=1$. Periodic boundary conditions are imposed at both ends of the computational domain $[-1,1]$. The numerical solution is advanced to the final time $T=0.5$ using a CFL number of $0.9$ and a relaxation parameter $\tau=10^{-8}$ for all the schemes i.e., CS-EBT2, UCS2 and IMEX RK2 schemes. Tables~\ref{Tabrho:broadwellerror}, \ref{Tabm:broadwellerror}, and \ref{Tabz:broadwellerror} present the $L^1$-error norms and the corresponding orders of convergence for the conservative variables $\rho$, $m$, and $z$, respectively. The numerical results demonstrate that the CS-EBT2, UCS2 and IMEX RK2 schemes attain the expected second-order order of accuracy for the alternative Broadwell model.
\renewcommand{\arraystretch}{1.2}{
\begin{table}
 \caption{Error and order of convergence for modify shallow water model with smooth initial condition \eqref{broadwell:continuous} for $\rho$ with $\tau = 10^{-8}$}
\centering
\begin{tabular}{*{7}{c}}
\toprule
\multirow{2}{*}{N} & 
\multicolumn{2}{c}{CS-EBT2} &
\multicolumn{2}{c}{UCS2} & 
\multicolumn{2}{c}{IMEX RK2} \\
\cmidrule(lr){2-3}
\cmidrule(lr){4-5}
\cmidrule(lr){6-7}
 & $L^{1}$-Error & Order 
 & $L^{1}$-Error & Order 
 & $L^{1}$-Error & Order \\
\midrule
64 & 2.9622e-04 &  -  & 2.8028e-04 & - & 1.1424e-02  & -\\

128 & 7.0987e-05  & 2.06 & 6.4033e-05 & 2.13  & 4.1643e-03 & 1.46\\

256  & 1.7335e-05 & 2.03 & 1.5033e-05 &  2.09 &  1.1791e-03 & 1.82\\

512  & 4.2730e-06 & 2.02 & 3.6286e-06  & 2.05 & 2.4880e-04 &  2.24\\

1024 &1.0608e-06 & 2.01 & 8.9413e-07  &  2.02 & 4.8338e-05 & 2.36\\ 
\bottomrule
\end{tabular}
\label{Tabrho:broadwellerror}
\end{table}}
  \renewcommand{\arraystretch}{1.2}{
\begin{table}
 \caption{Error and order of convergence for modify shallow water model with smooth initial condition \eqref{broadwell:continuous} for $m$ with $\tau = 10^{-8}$}
\centering
\begin{tabular}{*{7}{c}}
\toprule
\multirow{2}{*}{N} & 
\multicolumn{2}{c}{CS-EBT2} &
\multicolumn{2}{c}{UCS2} & 
\multicolumn{2}{c}{IMEX RK2} \\
\cmidrule(lr){2-3}
\cmidrule(lr){4-5}
\cmidrule(lr){6-7}
 & $L^{1}$-Error & Order 
 & $L^{1}$-Error & Order 
 & $L^{1}$-Error & Order \\
\midrule
64  & 2.5231e-04   & - & 2.3758e-04  & -  & 1.0198e-02  & -\\

128 & 6.0060e-05   &  2.07   &  5.4864e-05   & 2.11  & 3.7103e-03 & 1.46\\

256 & 1.4635e-05  &   2.04  & 1.3011e-05 & 2.08 &  1.0581e-03 & 1.81\\

512  & 3.6132e-06  & 2.02 & 3.1562e-06  & 2.04  &  2.2388e-04 & 2.24\\

1024 & 8.990e-07 & 2.01 & 7.8616e-07 & 2.01  & 4.3592e-05 & 2.36\\
\bottomrule
\end{tabular}
\label{Tabm:broadwellerror}
\end{table}}
\renewcommand{\arraystretch}{1.2}{
\begin{table}
 \caption{Error and order of convergence for modify shallow water model with smooth initial condition \eqref{broadwell:continuous} for $z$ with $\tau = 10^{-8}$}
\centering
\begin{tabular}{*{7}{c}}
\toprule
\multirow{2}{*}{N} & 
\multicolumn{2}{c}{CS-EBT2} &
\multicolumn{2}{c}{UCS2} & 
\multicolumn{2}{c}{IMEX RK2} \\
\cmidrule(lr){2-3}
\cmidrule(lr){4-5}
\cmidrule(lr){6-7}
 & $L^{1}$-Error & Order 
 & $L^{1}$-Error & Order 
 & $L^{1}$-Error & Order \\
\midrule
64  & 2.3012e-04 &  -  & 2.1650e-04  & - & 9.3592e-03  & -\\

128 & 5.4853e-05 &  2.07 &  4.9425e-05  & 2.13 & 3.4164e-03 & 1.45\\

256 & 1.3395e-05  &  2.03  & 1.1627e-05 & 2.09&  9.7147e-04 & 1.81\\

512  & 3.3051e-06  & 2.02 & 2.8088e-06  & 2.05  &  2.0541e-04 & 2.24\\

1024 & 8.2198e-07  & 2.01 & 6.9520e-07 & 2.01  & 3.9940e-05 & 2.36\\
\bottomrule
\end{tabular}
\label{Tabz:broadwellerror}
\end{table}}
\newline
\noindent \textbf{Non-smooth case:}
\begin{eqnarray}\label{broadwell:discontinuous}
    \begin{aligned}
    (\rho(x,0), m(x,0), z(x,0))=
        \begin{cases}
            (2,1,1), \, \text{if}\,\,\,\,x\leq0.2,\\
            (1,0.13962,1), \, \text{if}\,\,\,\,x>0.2,
            \end{cases} \\
    \end{aligned}
\end{eqnarray}
which introduces a discontinuity at $x = 0.2$. 
\renewcommand{\arraystretch}{1.2}{
\begin{table}
 \caption{Integral invariance for alternative Broadwell model with initial condition \eqref{broadwell:discontinuous} for $\rho$ and $m + \tau z$ with $\tau = 10^{-8}$ and CFL $0.9$}
\centering
\begin{tabular}{*{7}{c}}
\toprule
\multirow{2}{*}{$N$}
& \multicolumn{2}{c}{CS-EBT2}
& \multicolumn{2}{c}{UCS2}
& \multicolumn{2}{c}{IMEX RK2} \\
\cmidrule(lr){2-3} 
\cmidrule(lr){4-5}
\cmidrule(lr){6-7}
& Error $\rho$ & Error $m + \tau z$ & Error $\rho$ & Error $m + \tau z$ &  Error $\rho$ & Error $m + \tau z$ \\
\midrule

20 & 2.15094988e-01 & 1.87153063e-01 & 2.15094937e-01  & 1.85062784e-01 & 8.33333349e-02  & 8.33350379e-02\\

40 & 2.15094989e-01 & 1.86205478e-01 & 2.15094872e-01  & 1.85062365e-01 & 8.33333548e-02  & 8.33333556e-02\\

80 & 2.15094992e-01 & 1.85648639e-01 & 2.15094750e-01  & 1.850615590e-01 & 8.33333349e-02  & 8.33333357e-02\\

160 & 2.15094991e-01 & 1.85360041e-01 & 2.150945030e-01 & 1.85059977e-01 & 8.33333349e-02 & 8.33333357e-02\\

320 & 2.15094989e-01 & 1.85213185e-01 & 2.15094014e-01 &1.85056821e-01 & 8.33333349e-02 &8.33333357e-02\\

640 & 2.15094992e-01 & 1.85138455e-01 & 2.15093023e-01  &1.85050458e-01 & 8.33333349e-02  &8.33333357e-02\\

1280 & 2.15094991e-01 & 1.85100929e-01 & 2.15091043e-0 & 1.85037747e-01 & 8.33333349e-02 & 8.33333357e-02\\
\bottomrule
\end{tabular}
\label{Tabbroadwell}
\end{table}} 
\begin{figure}[!ht]
    \begin{minipage}[b]{0.32\linewidth}
        \includegraphics[width=\linewidth]{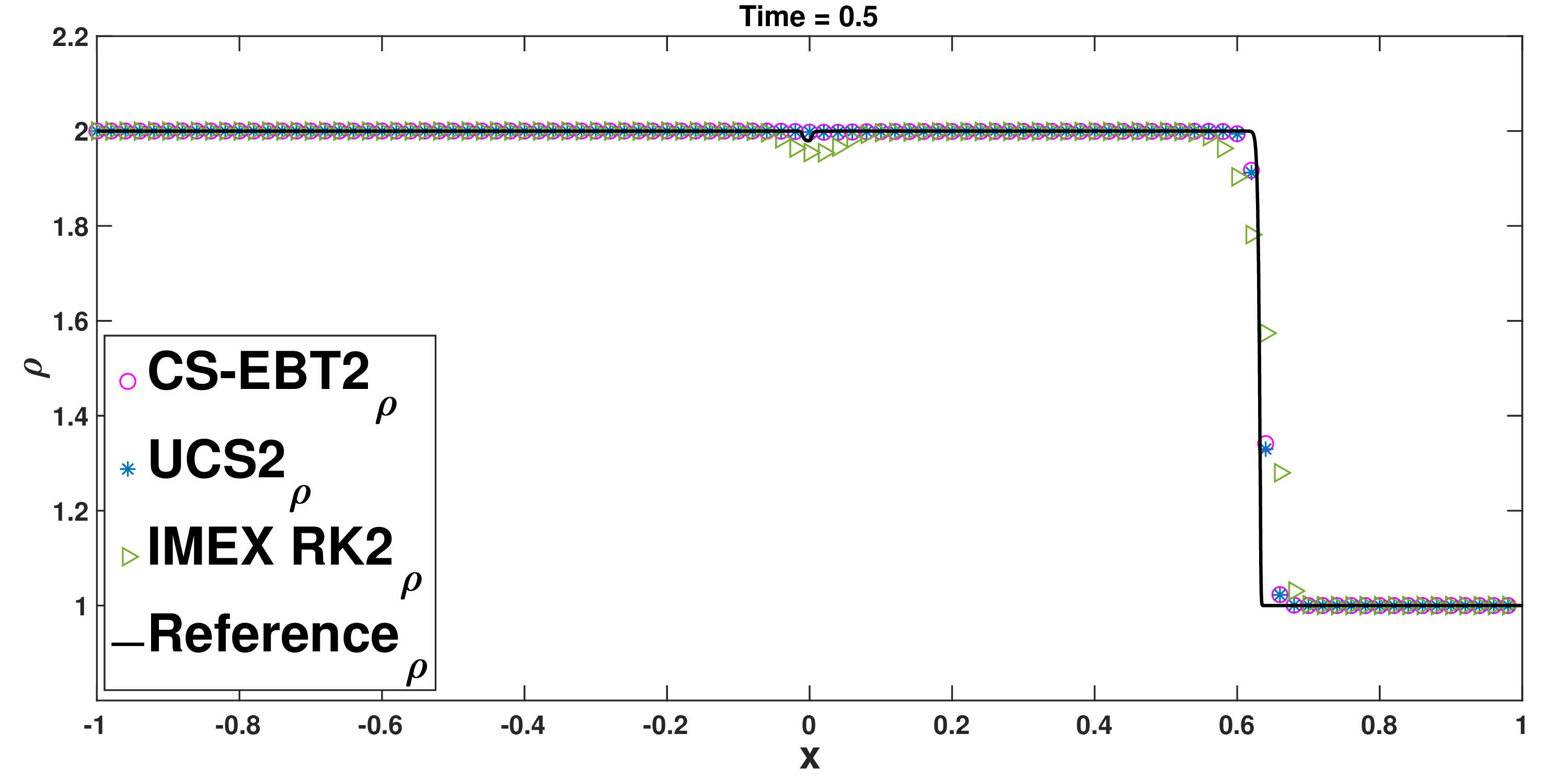}
    \end{minipage}
    \hfill
    \begin{minipage}[b]{0.32\linewidth}
        \includegraphics[width=\linewidth]{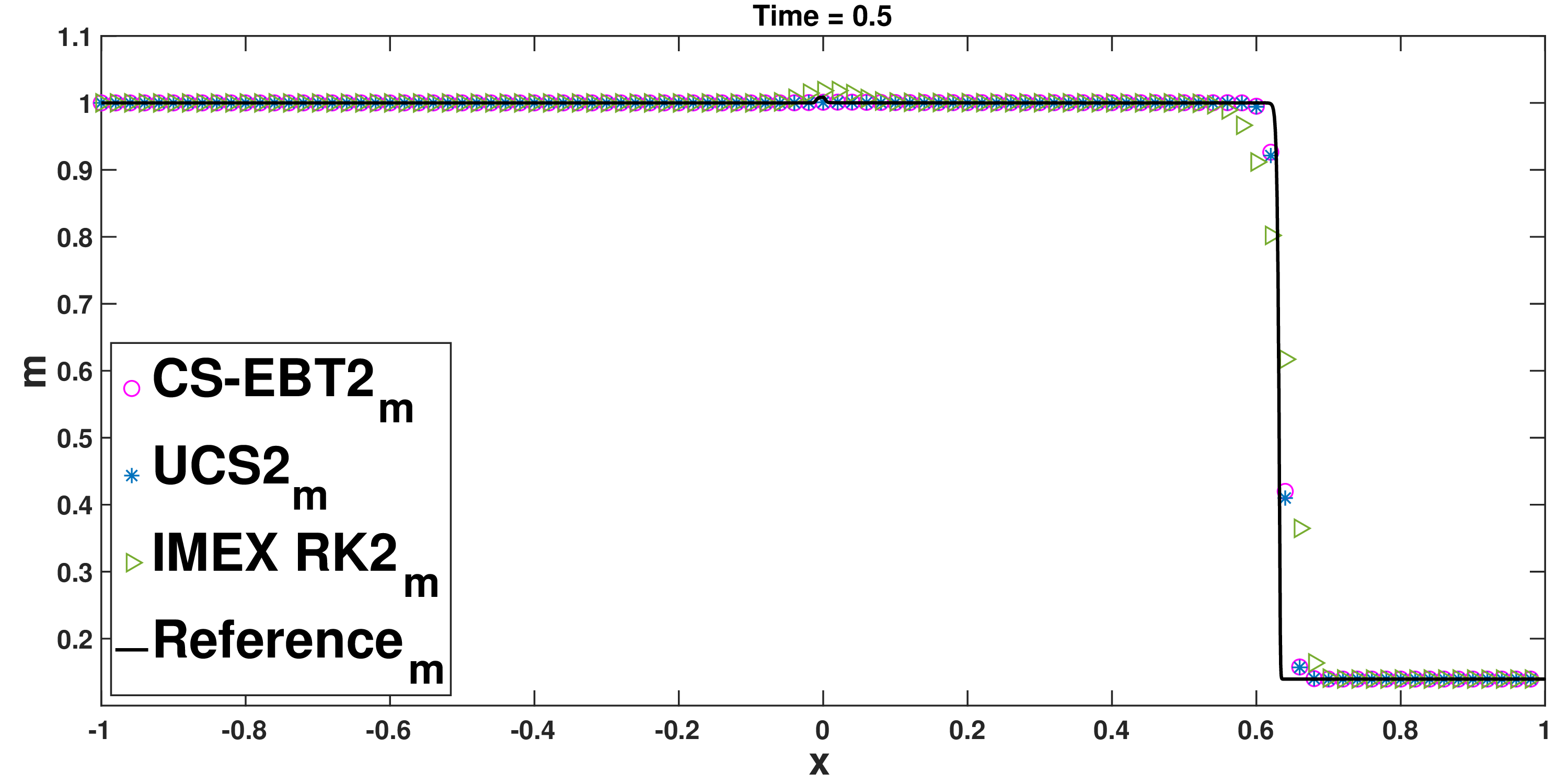}
    \end{minipage}
    \hfill
    \begin{minipage}[b]{0.32\linewidth}
        \includegraphics[width=\linewidth]{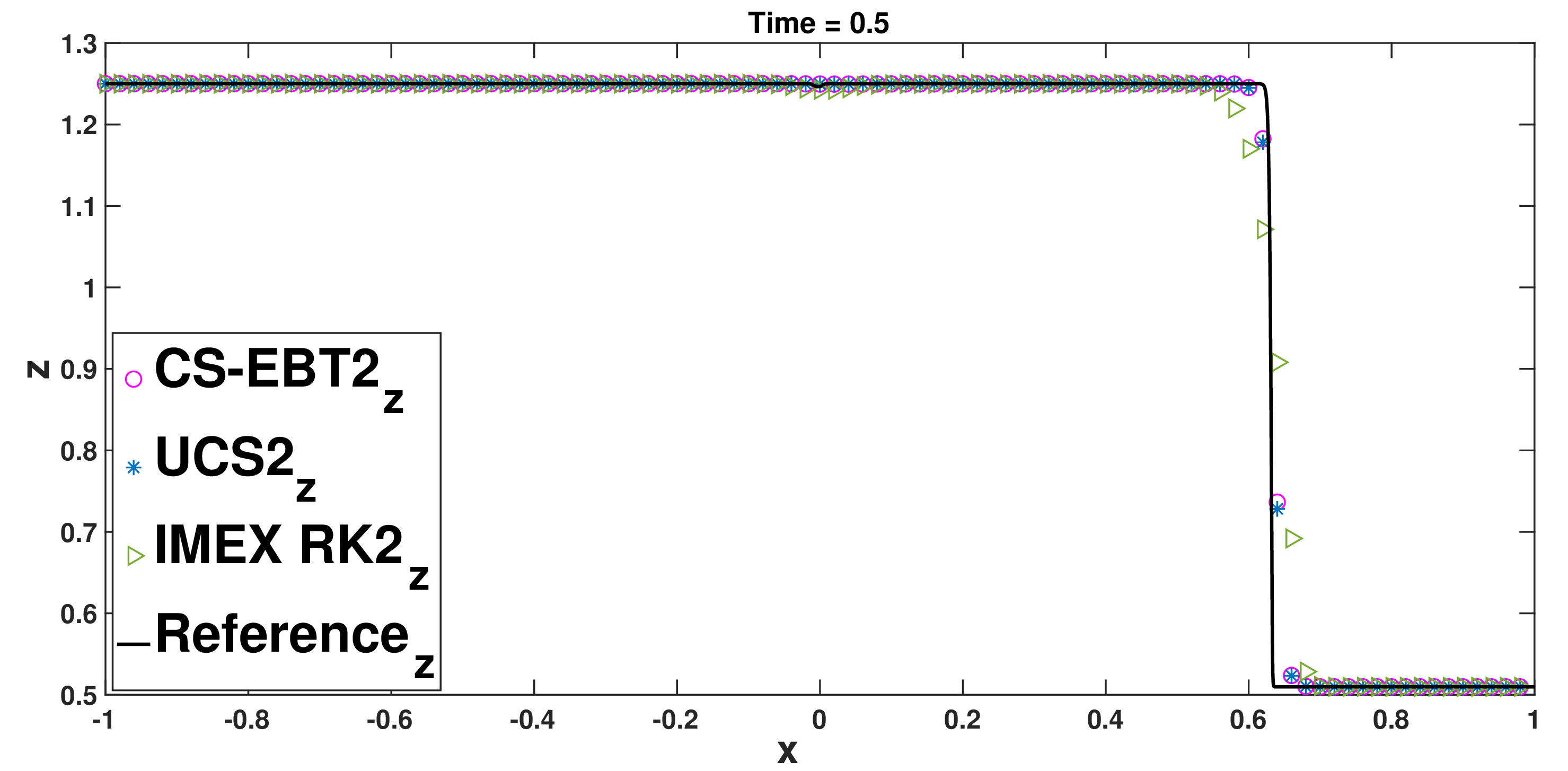}
    \end{minipage}
    \hfill
    \begin{minipage}[b]{0.32\linewidth}
        \includegraphics[width=\linewidth]{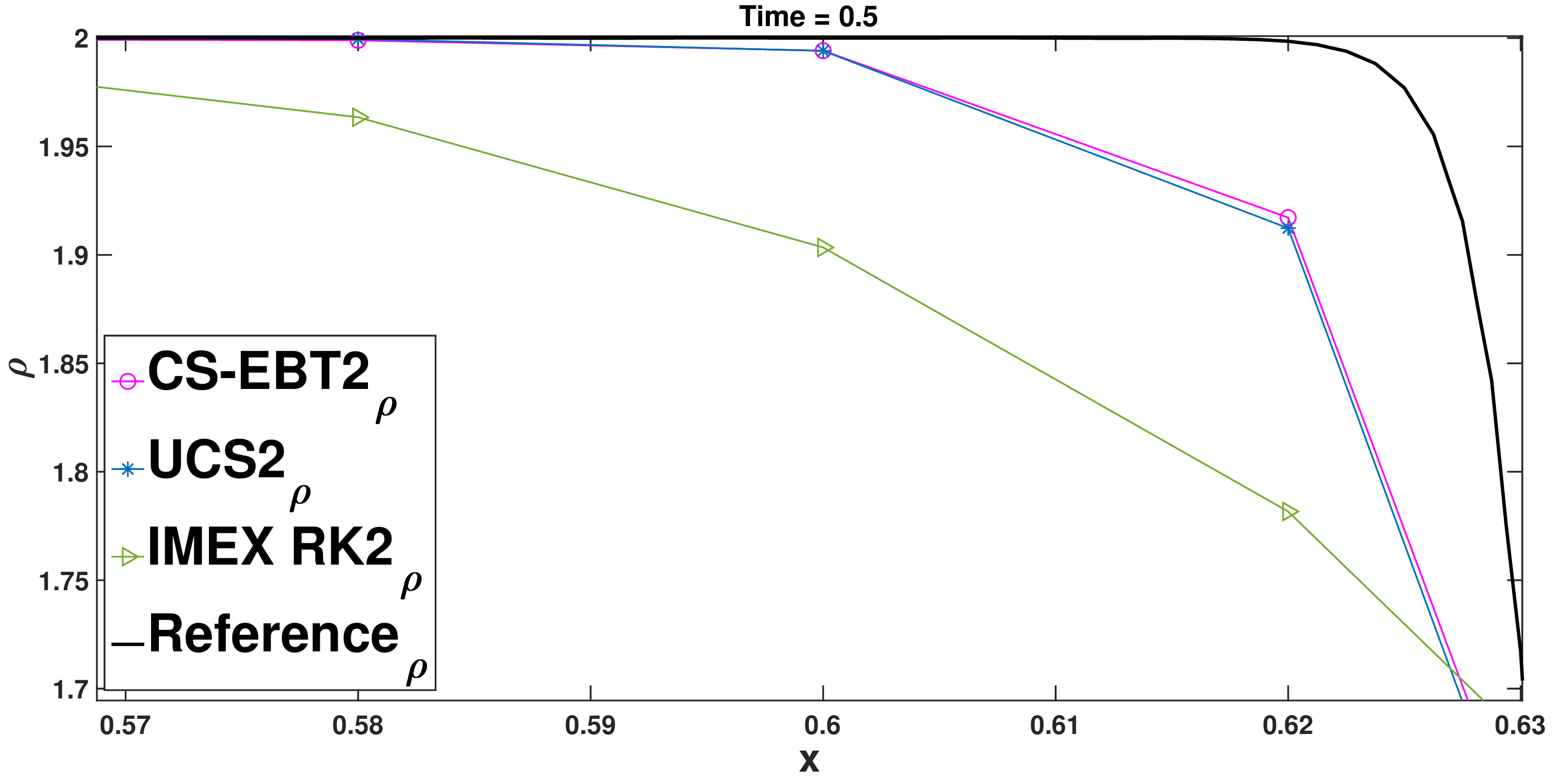}
    \end{minipage}
    \hfill
    \begin{minipage}[b]{0.32\linewidth}
        \includegraphics[width=\linewidth]{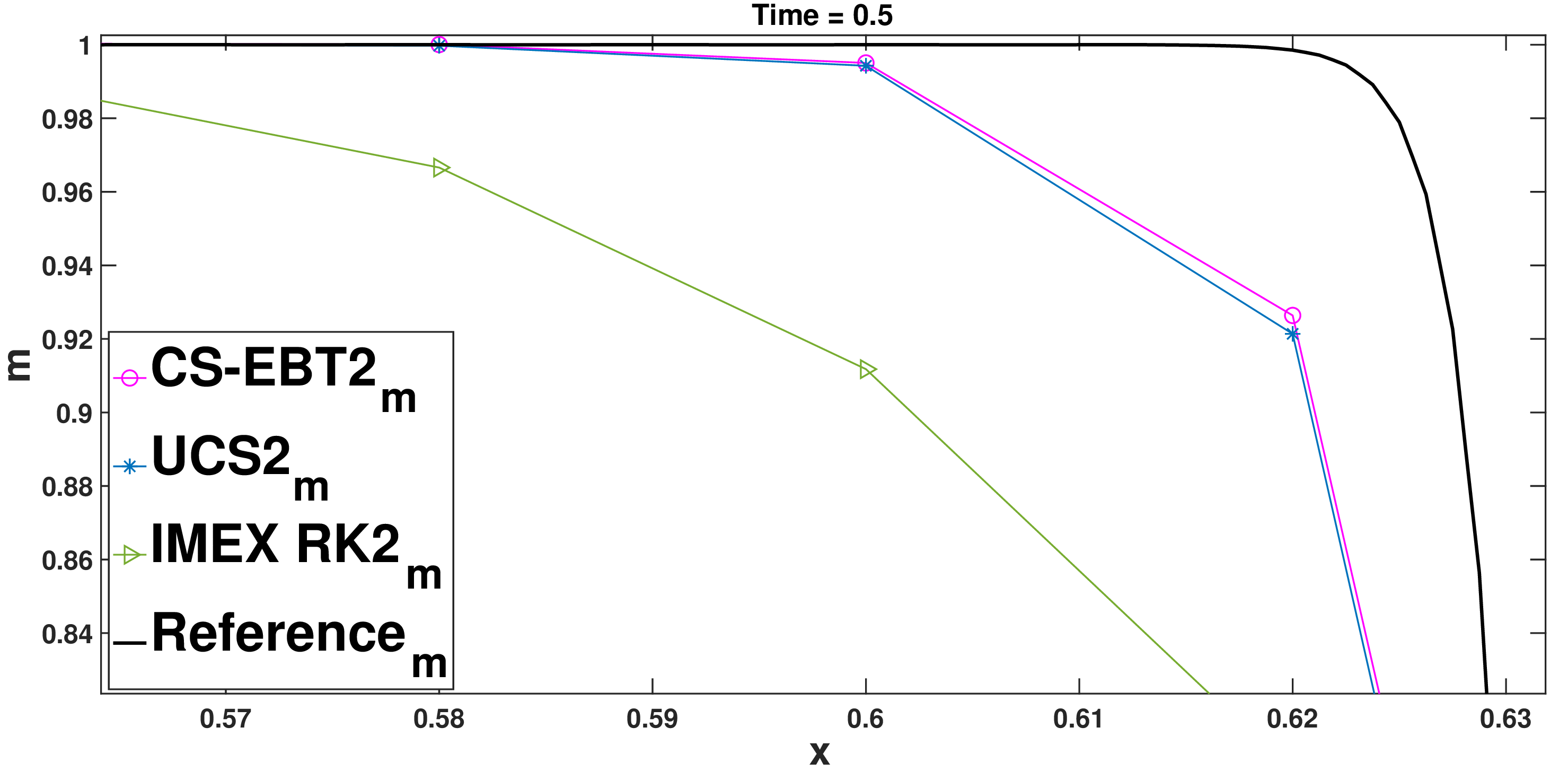}
    \end{minipage}
    \hfill
    \begin{minipage}[b]{0.32\linewidth}
        \includegraphics[width=\linewidth]{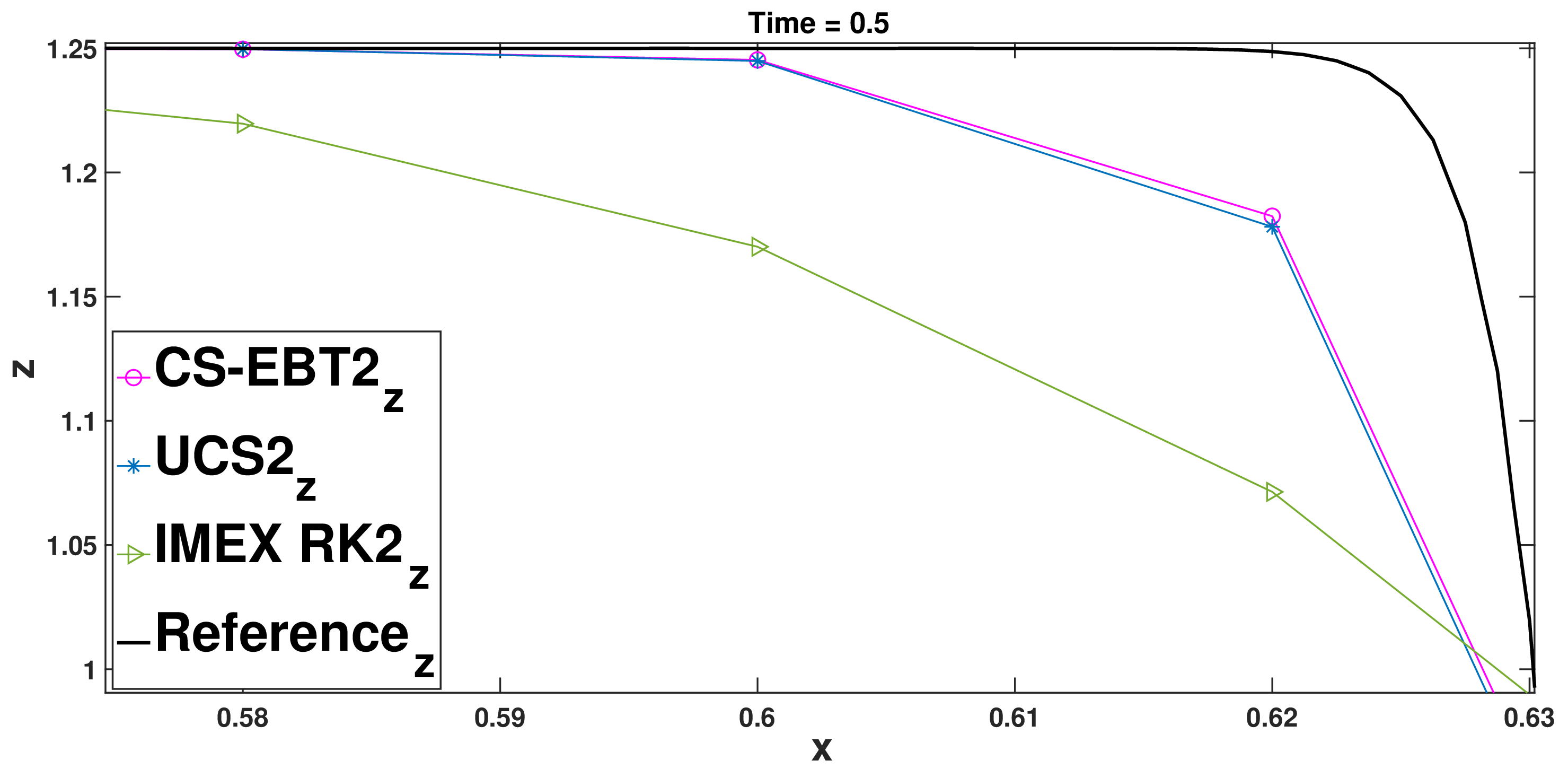}
    \end{minipage}
    \caption{Alternative Broadwell model with non-smooth case: comparison between Numerical solution $\rho$, $m$, $z$(first row) and zoomed version of $\rho$, $m$, $z$(second row) and the reference solution (IMEX RK2) with $\tau = 10^{-8}$, CFL $0.9$ and $N=200$.}
    \label{(3a)}
\end{figure}
\begin{figure}[!ht]
    \begin{minipage}[b]{0.48\linewidth}
        \includegraphics[width=\linewidth]{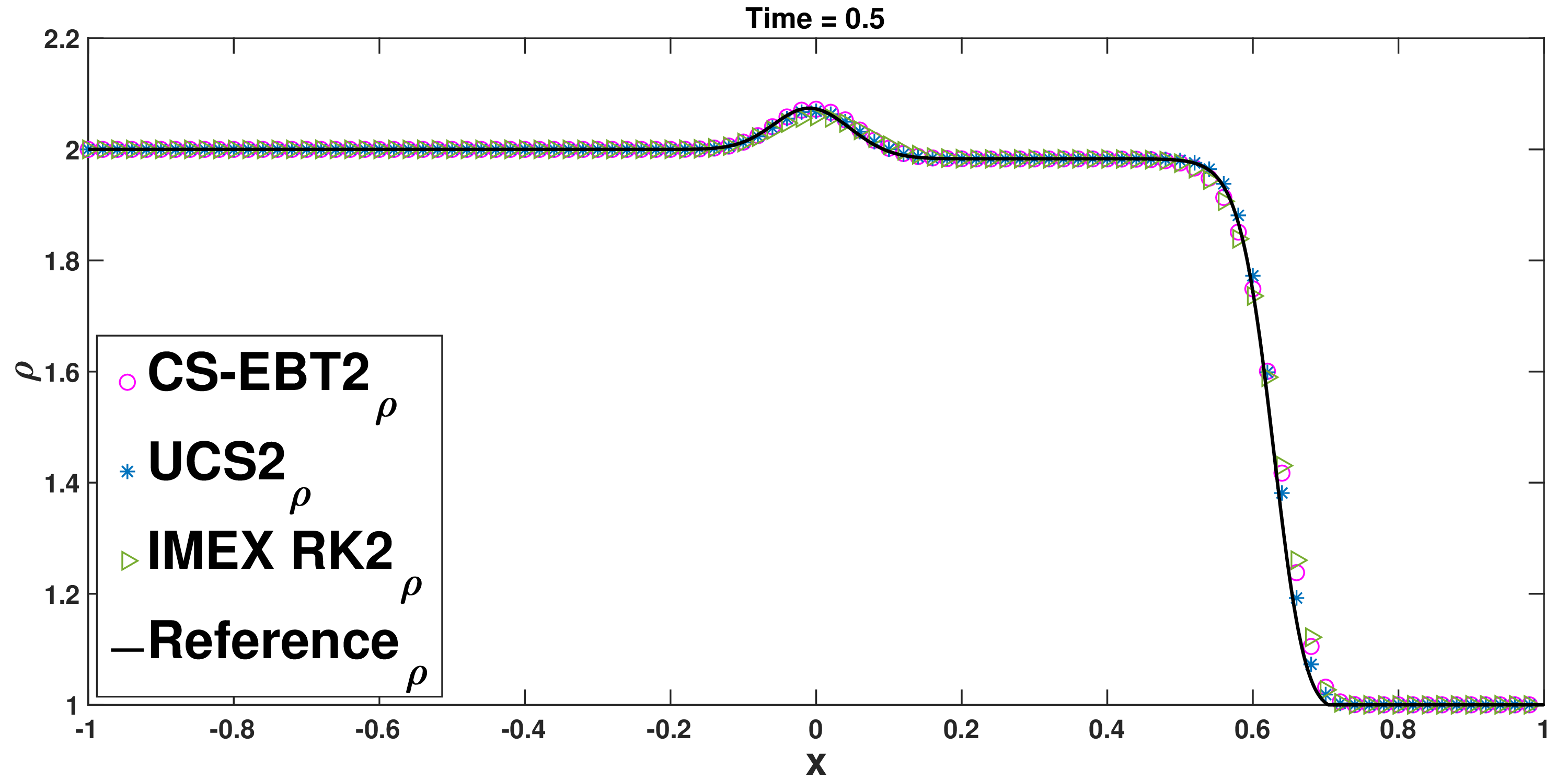}
    \end{minipage}
    \hfill
    \begin{minipage}[b]{0.48\linewidth}
        \includegraphics[width=\linewidth]{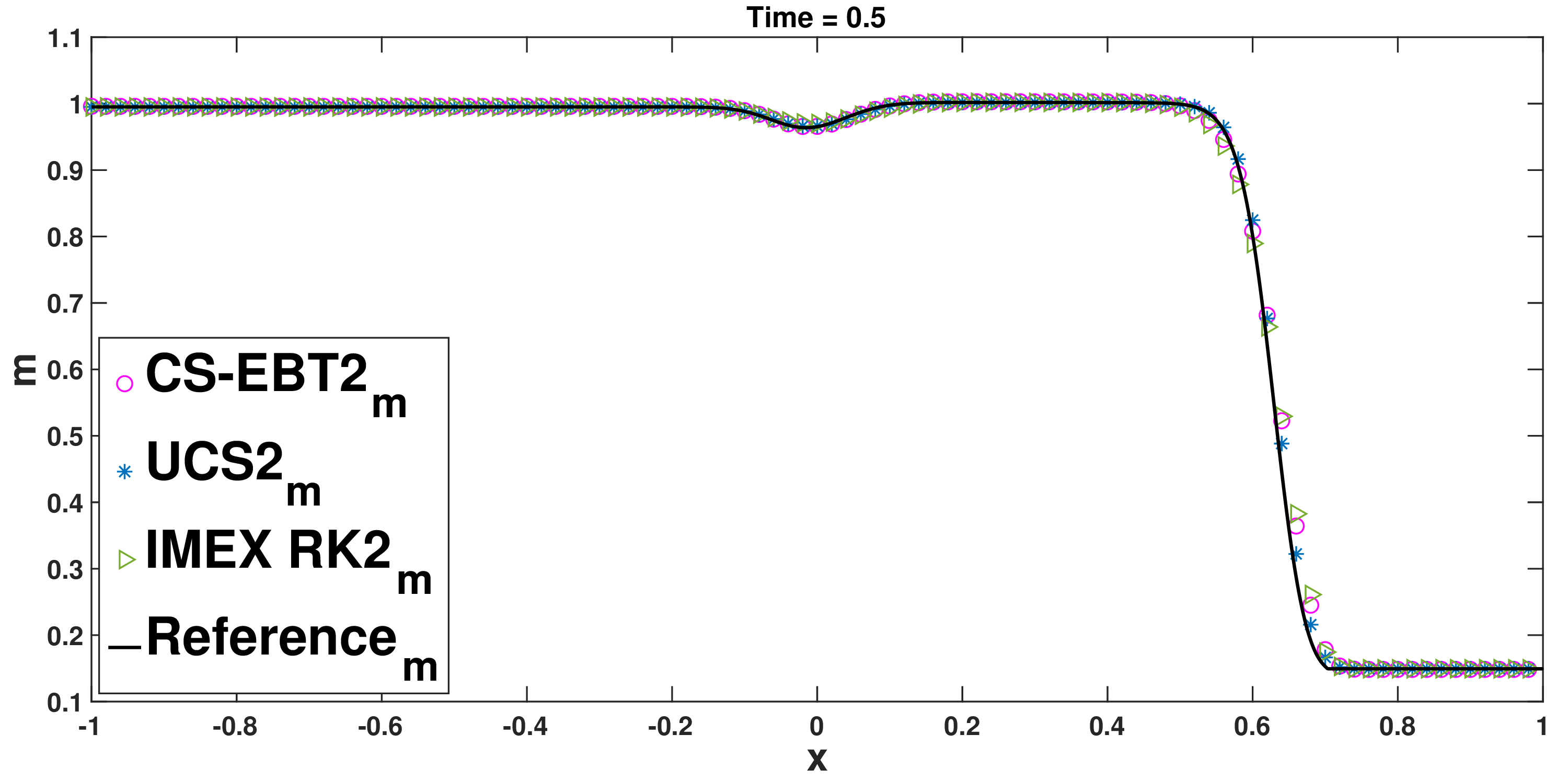}
    \end{minipage}

    \vspace{0.5cm} 

    \centering
    \includegraphics[width=0.5\linewidth]{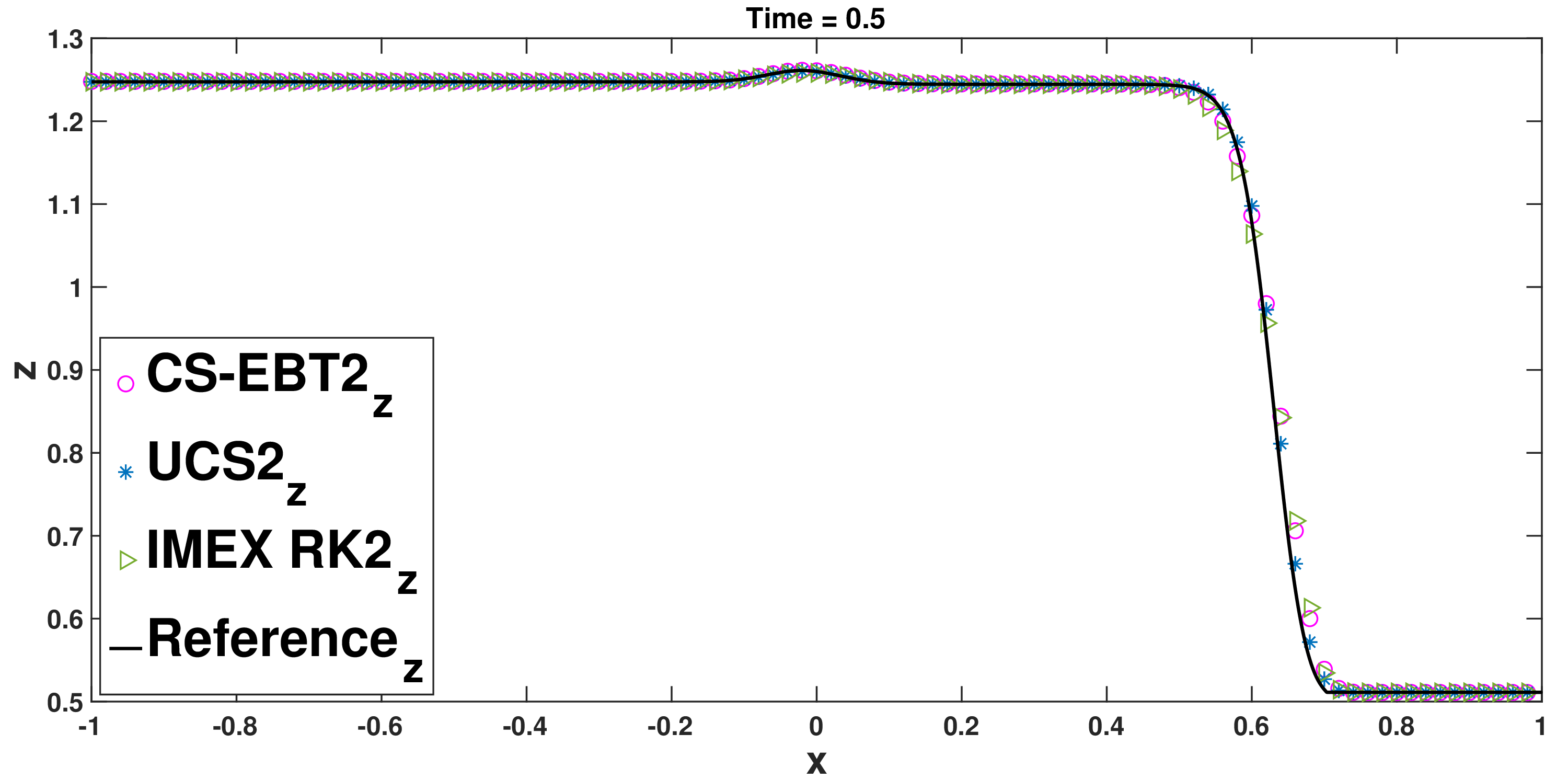}
    
    \caption{Alternative Broadwell model with non-smooth case: comparison between numerical solution $\rho$ (left), $m$ (right) and $z$ (center) and the reference solution (IMEX RK2) with $\tau = 0.02$, CFL $0.9$ and $N=200$.}
    \label{(3b)}
\end{figure}

\begin{figure}[!ht]
    \begin{minipage}[b]{0.48\linewidth}
        \includegraphics[width=\linewidth]{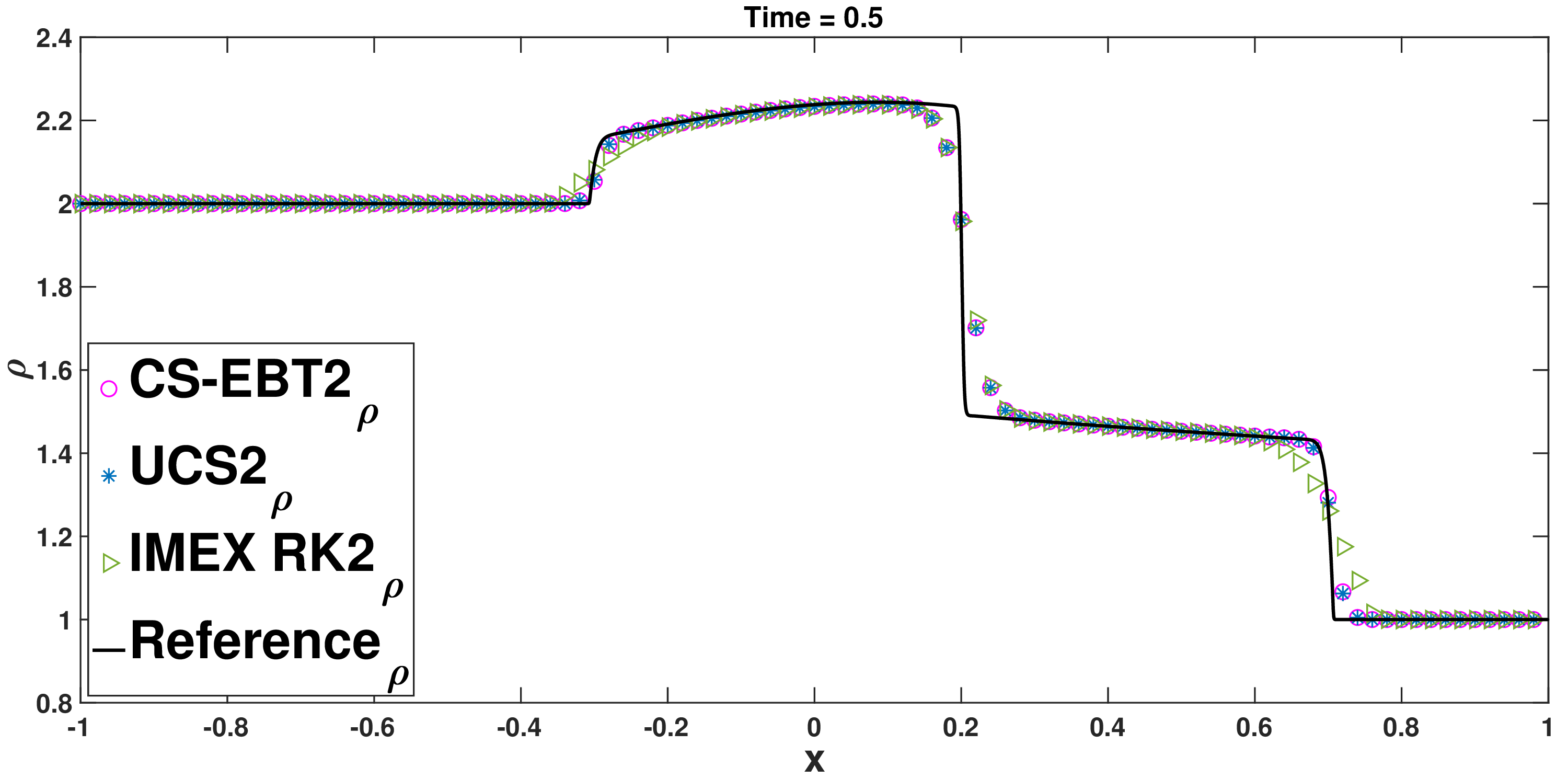}
    \end{minipage}
    \hfill
    \begin{minipage}[b]{0.48\linewidth}
        \includegraphics[width=\linewidth]{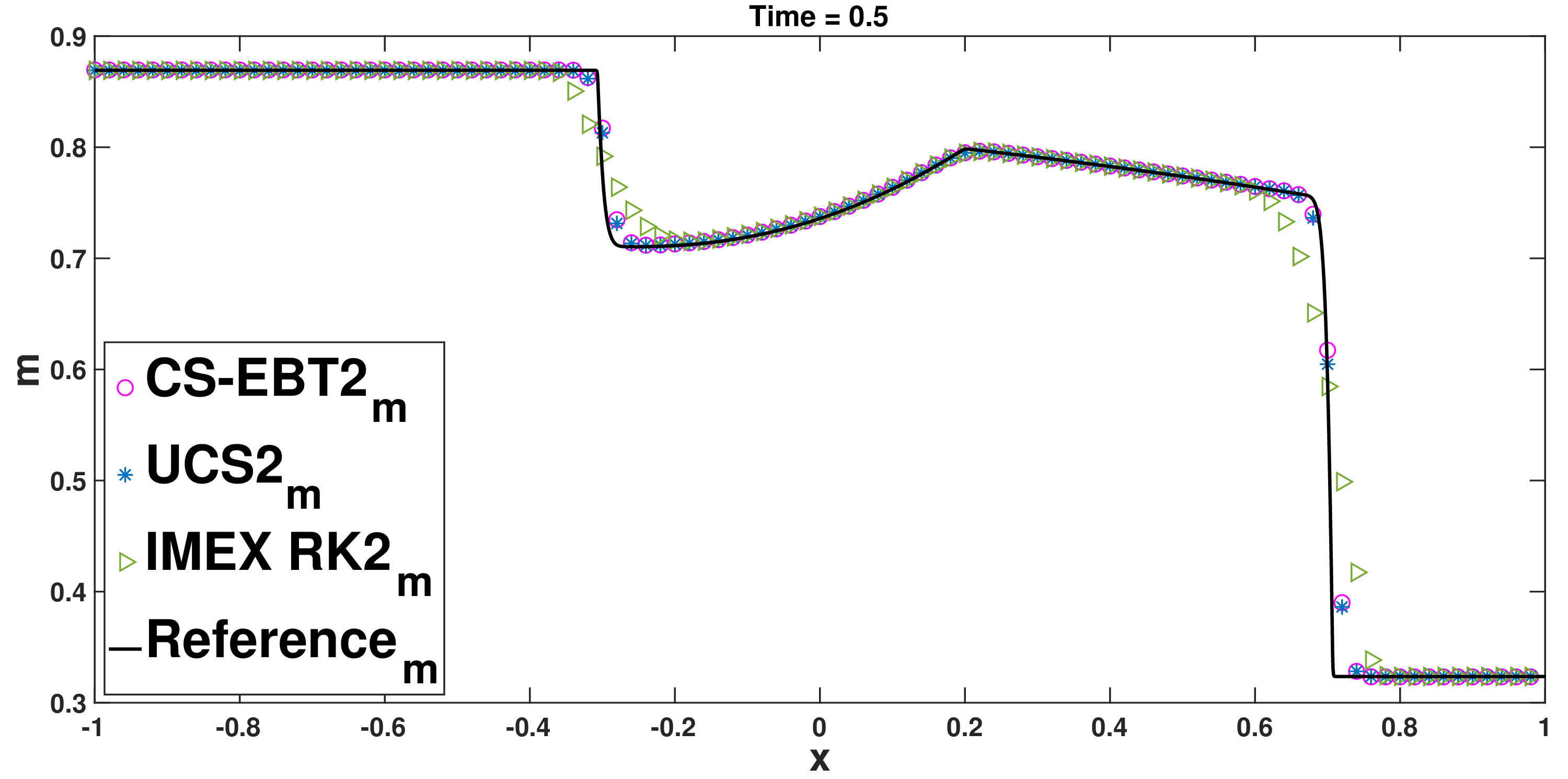}
    \end{minipage}

    \vspace{0.5cm} 

    \centering
    \includegraphics[width=0.5\linewidth]{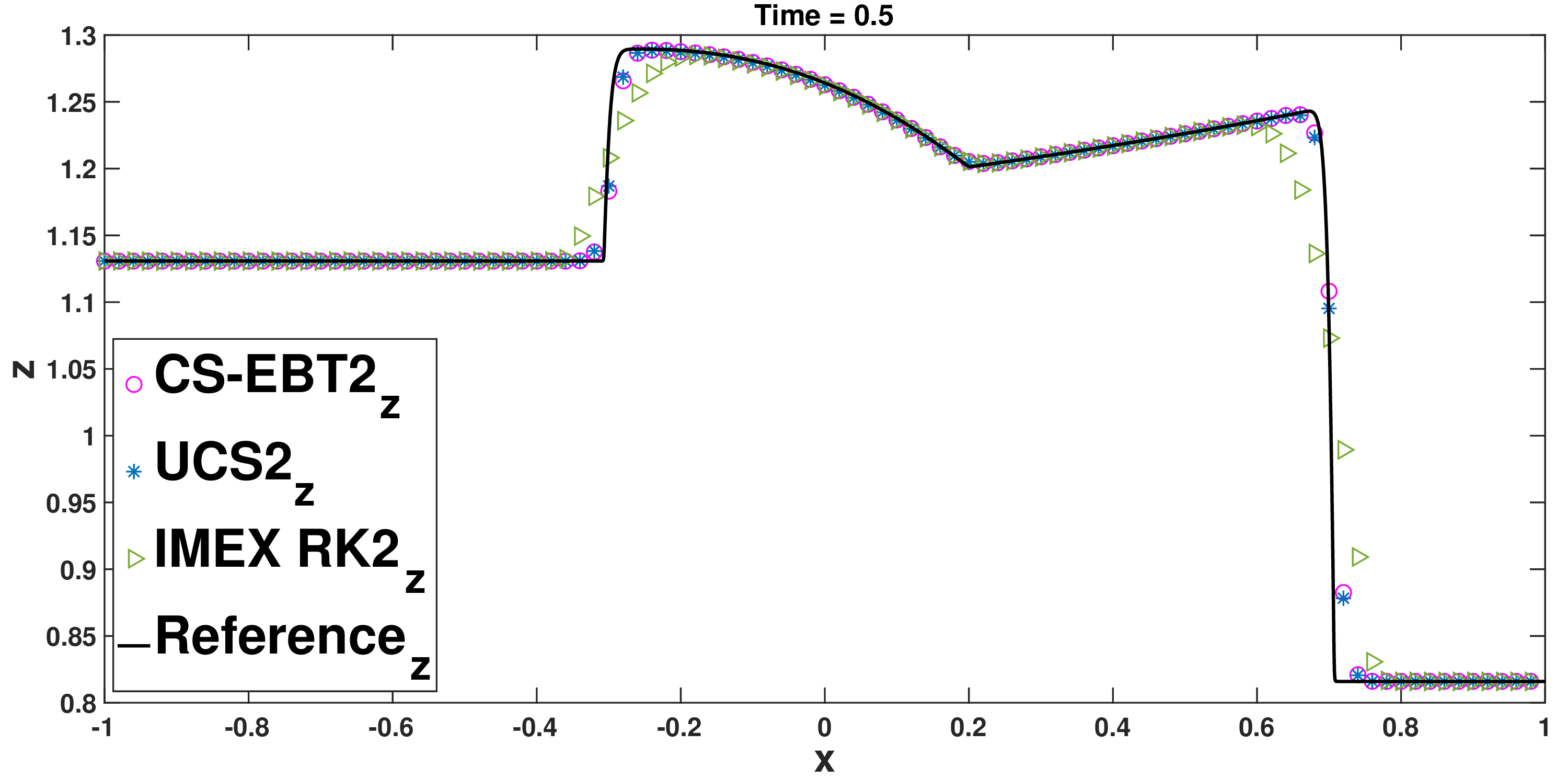}
    
    \caption{Alternative Broadwell model with non-smooth case: comparison between numerical solution $\rho$ (left), $m$ (right) and $z$ (center) and the reference solution (IMEX RK2) with $\tau = 1$, CFL $0.9$ and $N=200$.}
    \label{(3c)}
\end{figure}

Transmissive boundary conditions are applied at both ends of the computational domain $[-1,1]$ for the CS-EBT2, UCS2 and IMEX RK2 schemes. The solution is evolved to the final time $T=0.5$ using a CFL number of $0.9$ and a relaxation parameter $\tau=10^{-8}$, with the corresponding results shown in Fig.~\ref{(3a)}. The numerical solutions for $\tau=0.02$ and $\tau=1$, obtained using the same CFL number and final time, are presented in Figs.~\ref{(3b)} and~\ref{(3c)}, respectively.
For this model, the integral invariance property is examined numerically, and the corresponding results are reported in Table~\ref{Tabbroadwell}.
\begin{figure}[!ht]
    \begin{minipage}[b]{0.32\linewidth}
        \includegraphics[width=\linewidth]{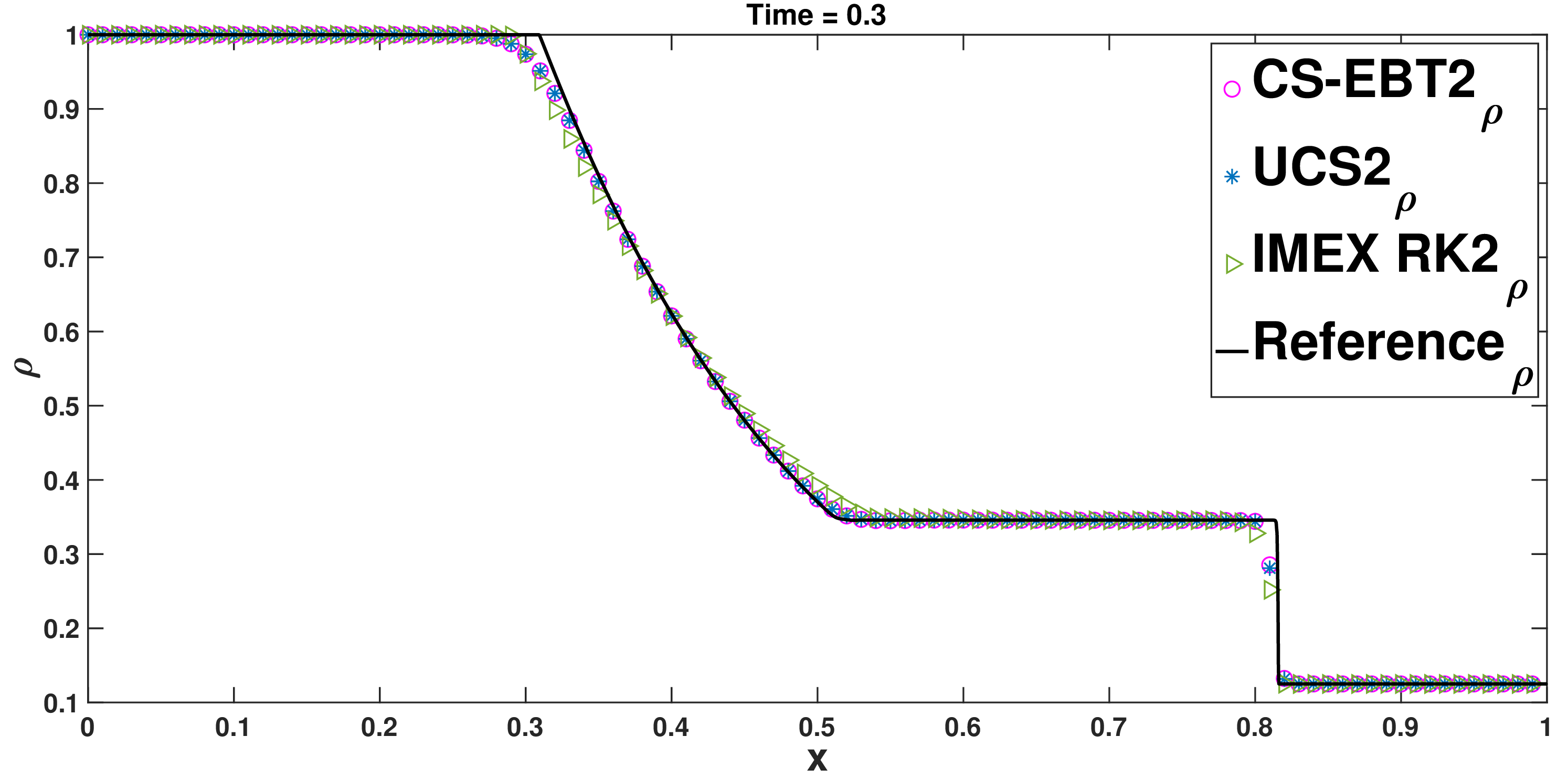}
    \end{minipage}
    \hfill
    \begin{minipage}[b]{0.32\linewidth}
        \includegraphics[width=\linewidth]{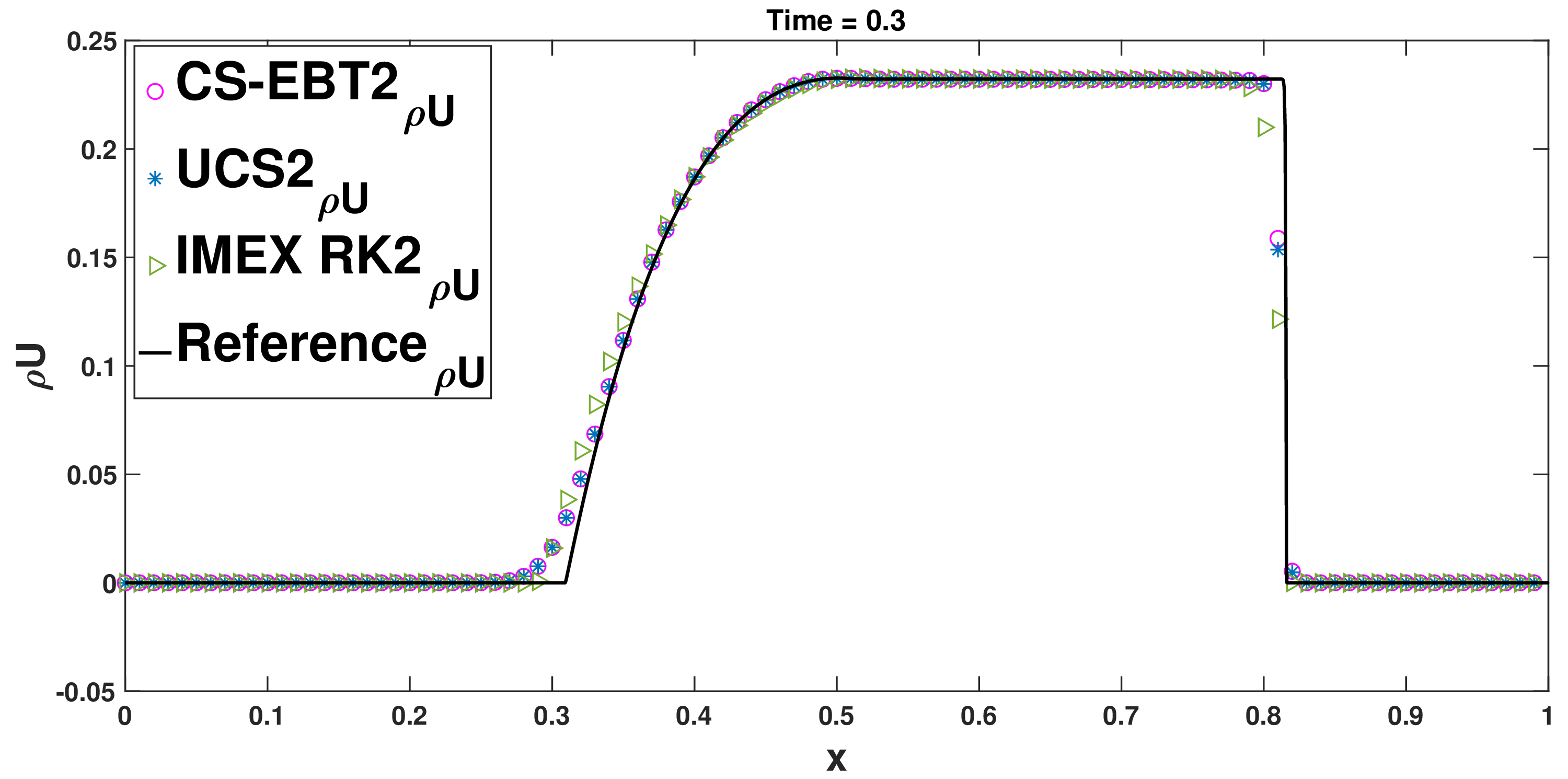}
    \end{minipage}
    \hfill
    \begin{minipage}[b]{0.32\linewidth}
    \includegraphics[width=\linewidth]{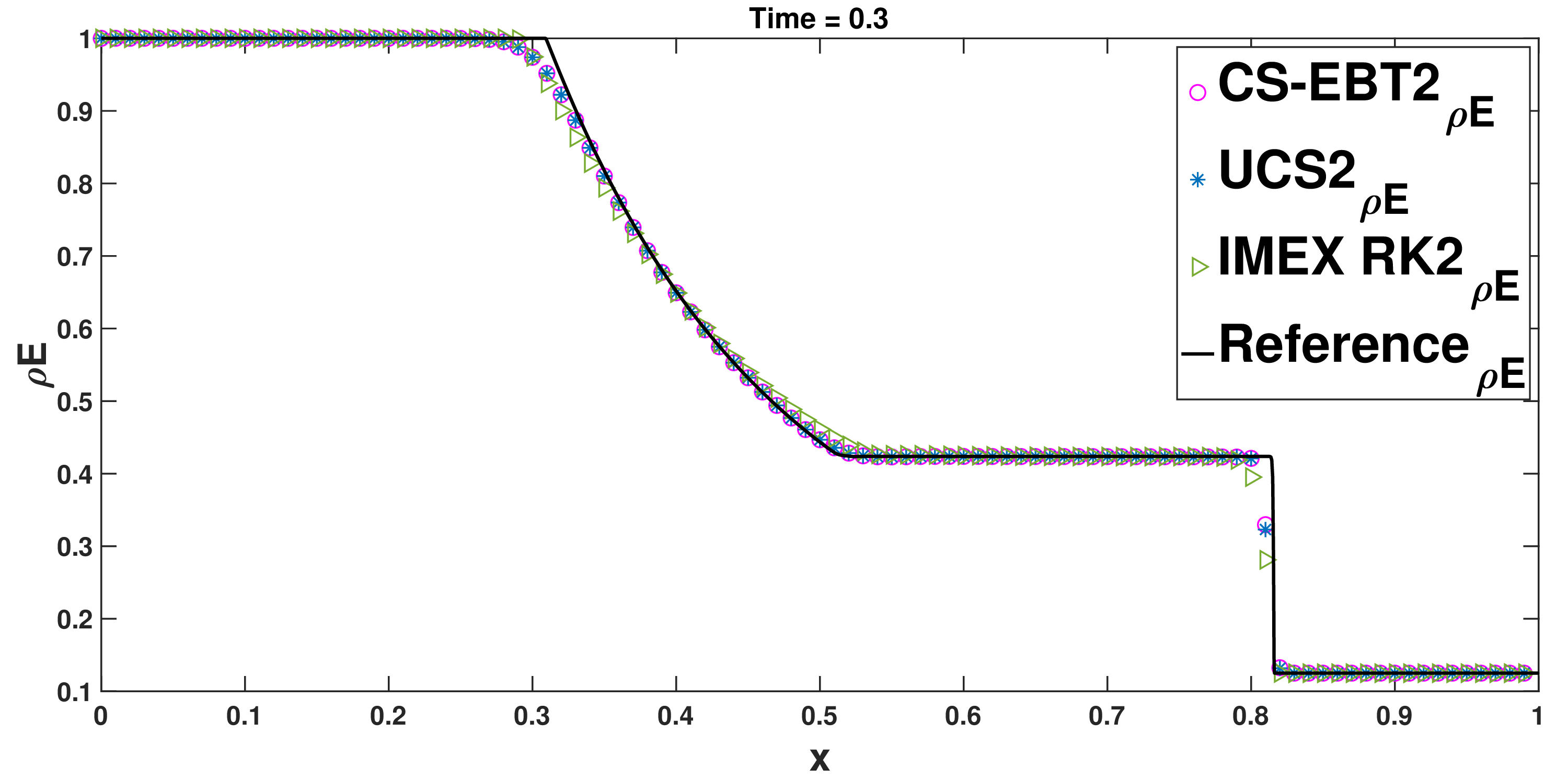}
    \end{minipage}
    \hfill
    \begin{minipage}[b]{0.32\linewidth}
        \includegraphics[width=\linewidth]{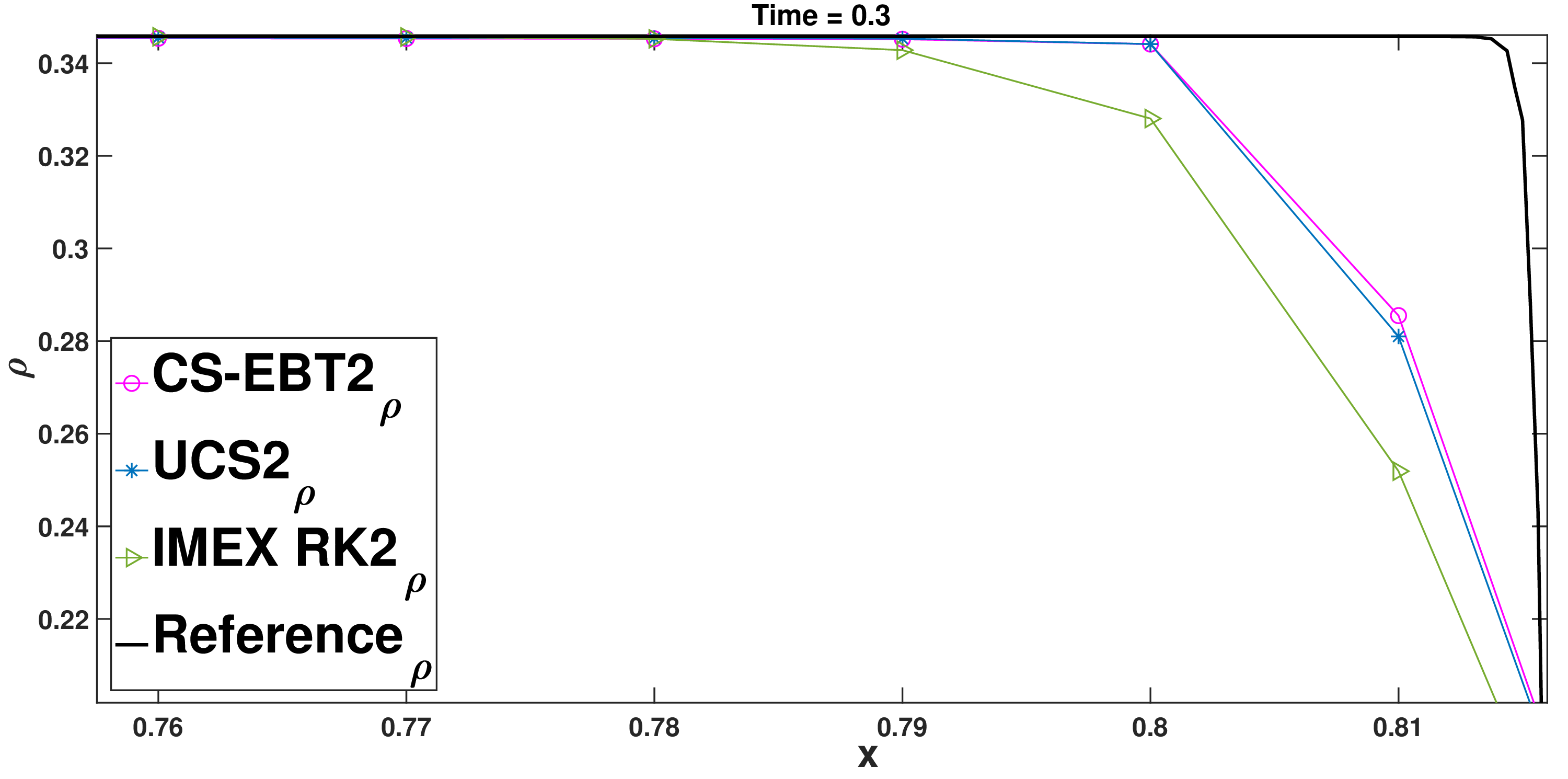}
    \end{minipage}
    \hfill
    \begin{minipage}[b]{0.32\linewidth}
        \includegraphics[width=\linewidth]{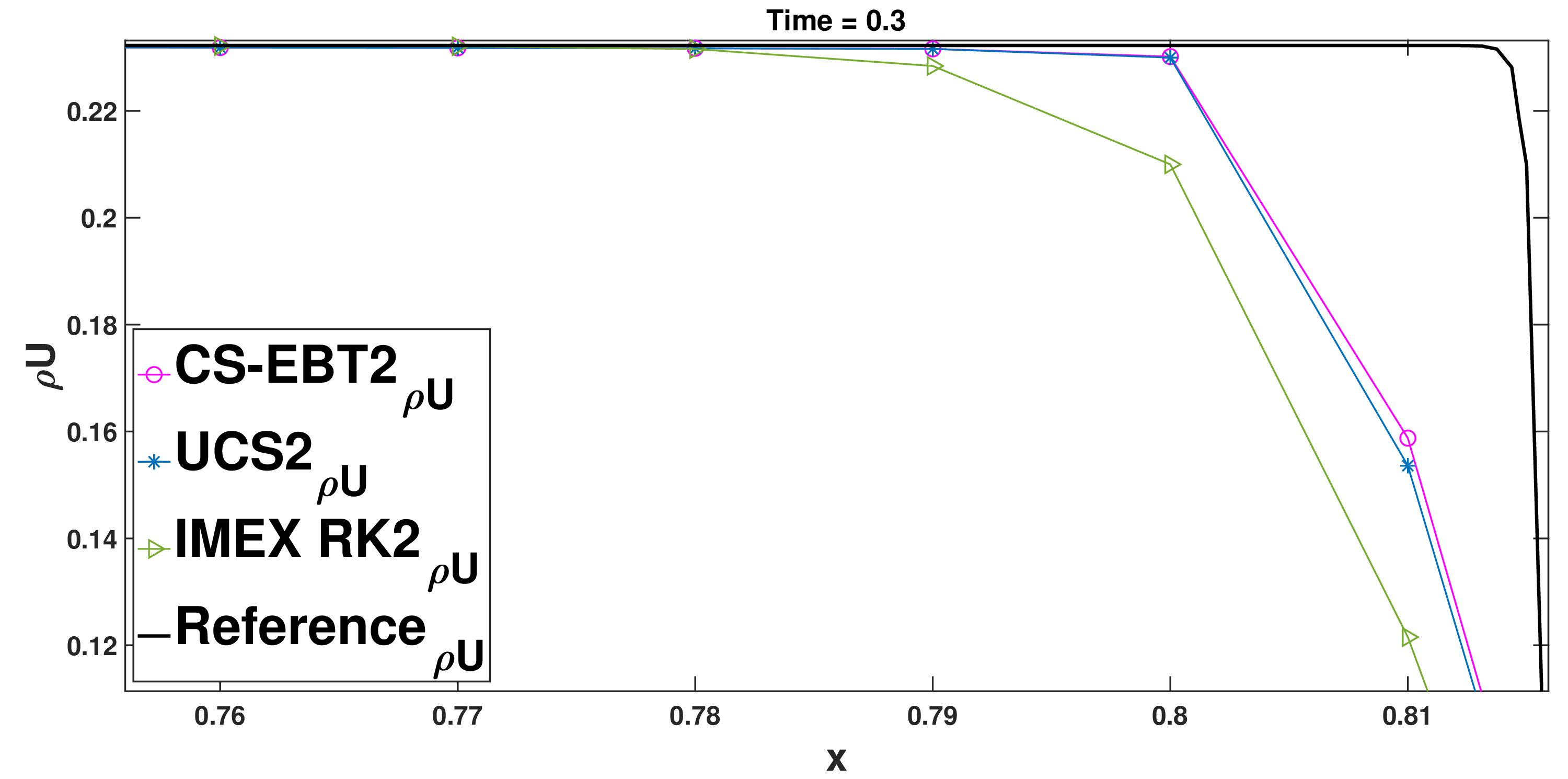}
    \end{minipage}
    \hfill
    \begin{minipage}[b]{0.32\linewidth}
    \includegraphics[width=\linewidth]{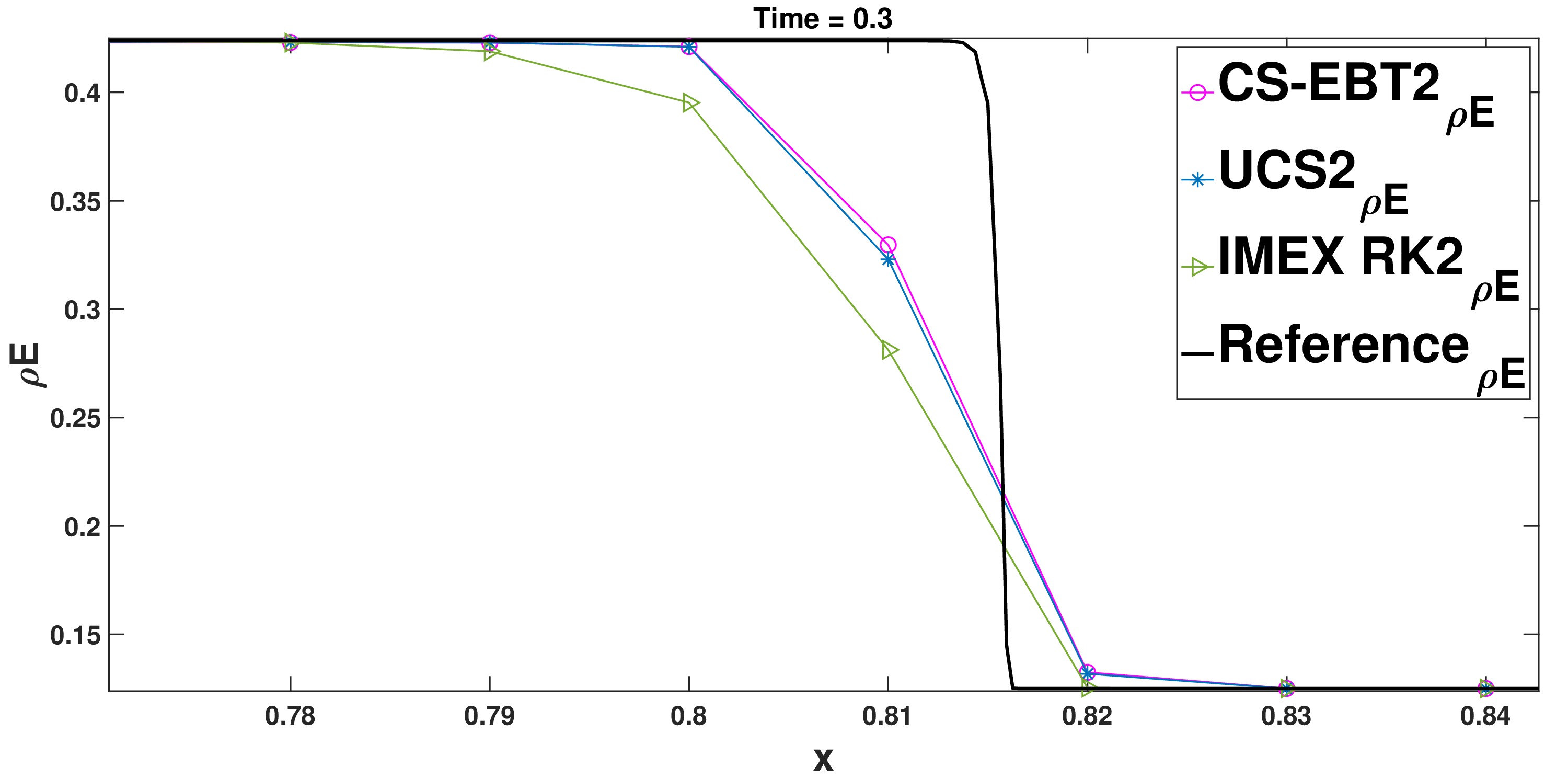}
    \end{minipage}   
    \caption{Alternative Euler equations with heat transfer model with non-smooth case: comparison between numerical solution $\rho$, $\rho u$ and $\rho E$ (first row), zoomed version of $\rho$, $\rho u$ and $\rho E$ (second row) and the reference solution (IMEX RK2) with $\tau = 10^{-8}$, CFL $0.9$ and $N = 200$.}
    \label{(4a)}
\end{figure}
\subsubsection{Euler Equations with Heat Transfer}
We consider the alternative relaxation Euler equations with heat transfer
\eqref{Modified:Euler_with_heat_transfer}. The initial data for the conserved
variables density $(\rho)$, velocity $(u)$, and total energy $(E)$ are
prescribed as follows
\begin{eqnarray}\label{discontinuous_euler}
    \begin{aligned}
    (\rho(x,0), u(x,0), E(x,0))=
        \begin{cases}
            (1,0,1), \, \text{if}\,\,\,x\leq0.5,\\
            (1/8,0,1) \, \text{if}\,\,\,x>0.5.
            \end{cases} \\
    \end{aligned}
\end{eqnarray}

\begin{figure}[!ht]
    \begin{minipage}[b]{0.48\linewidth}
        \includegraphics[width=\linewidth]{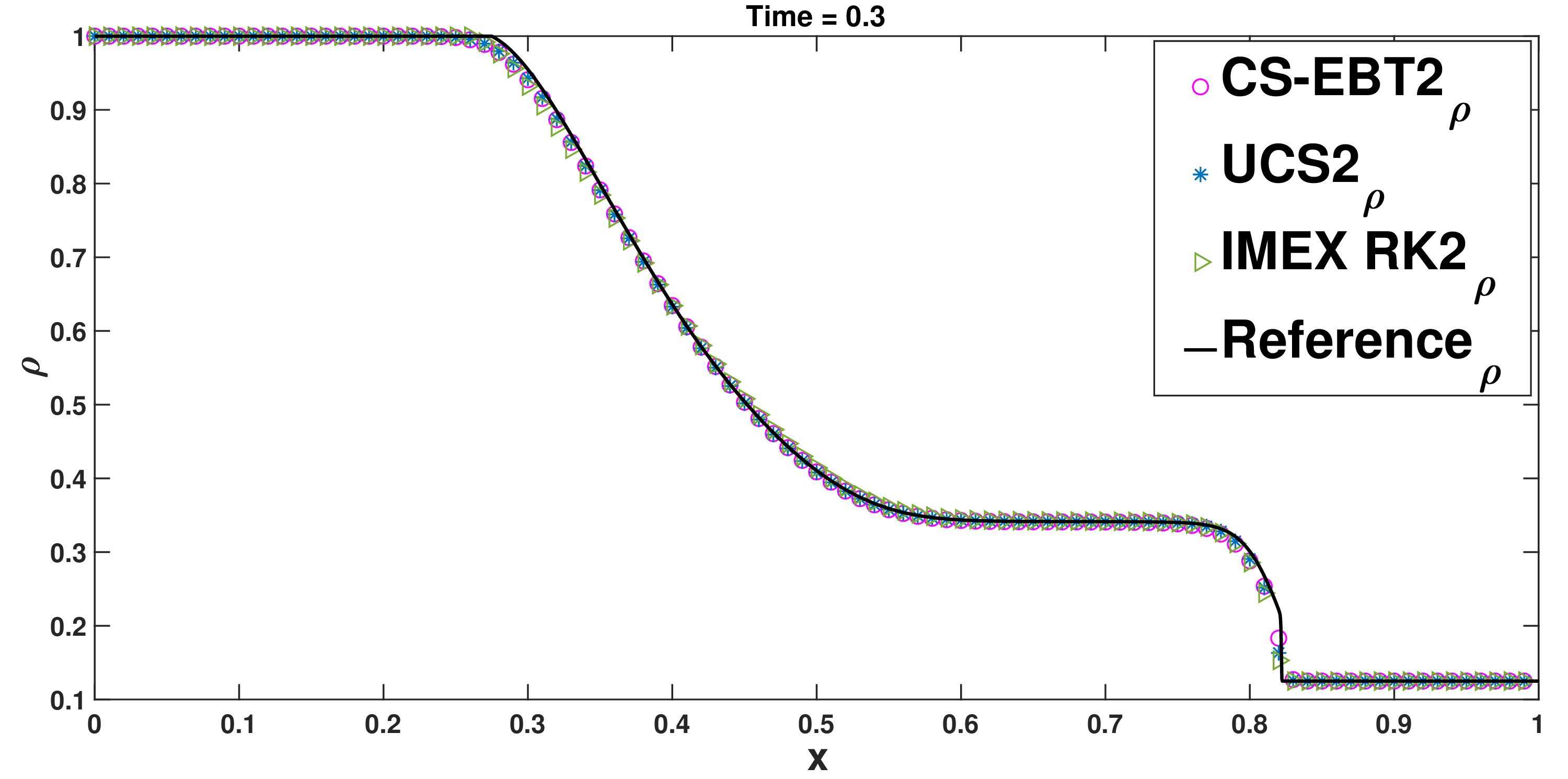}
    \end{minipage}
    \hfill
    \begin{minipage}[b]{0.48\linewidth}
        \includegraphics[width=\linewidth]{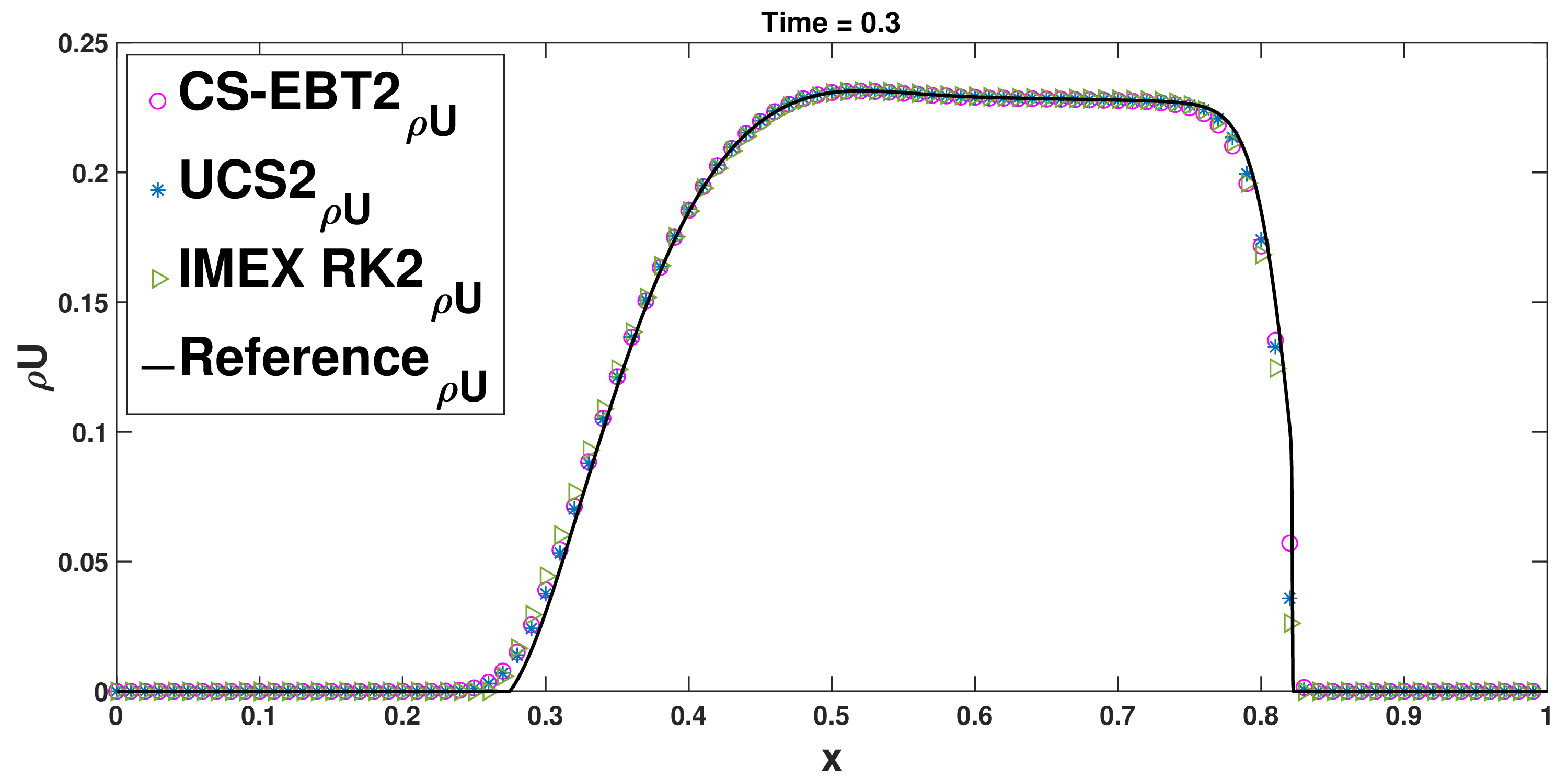}
    \end{minipage}
    \vspace{0.5cm} 
    \centering
    \includegraphics[width=0.5\linewidth]{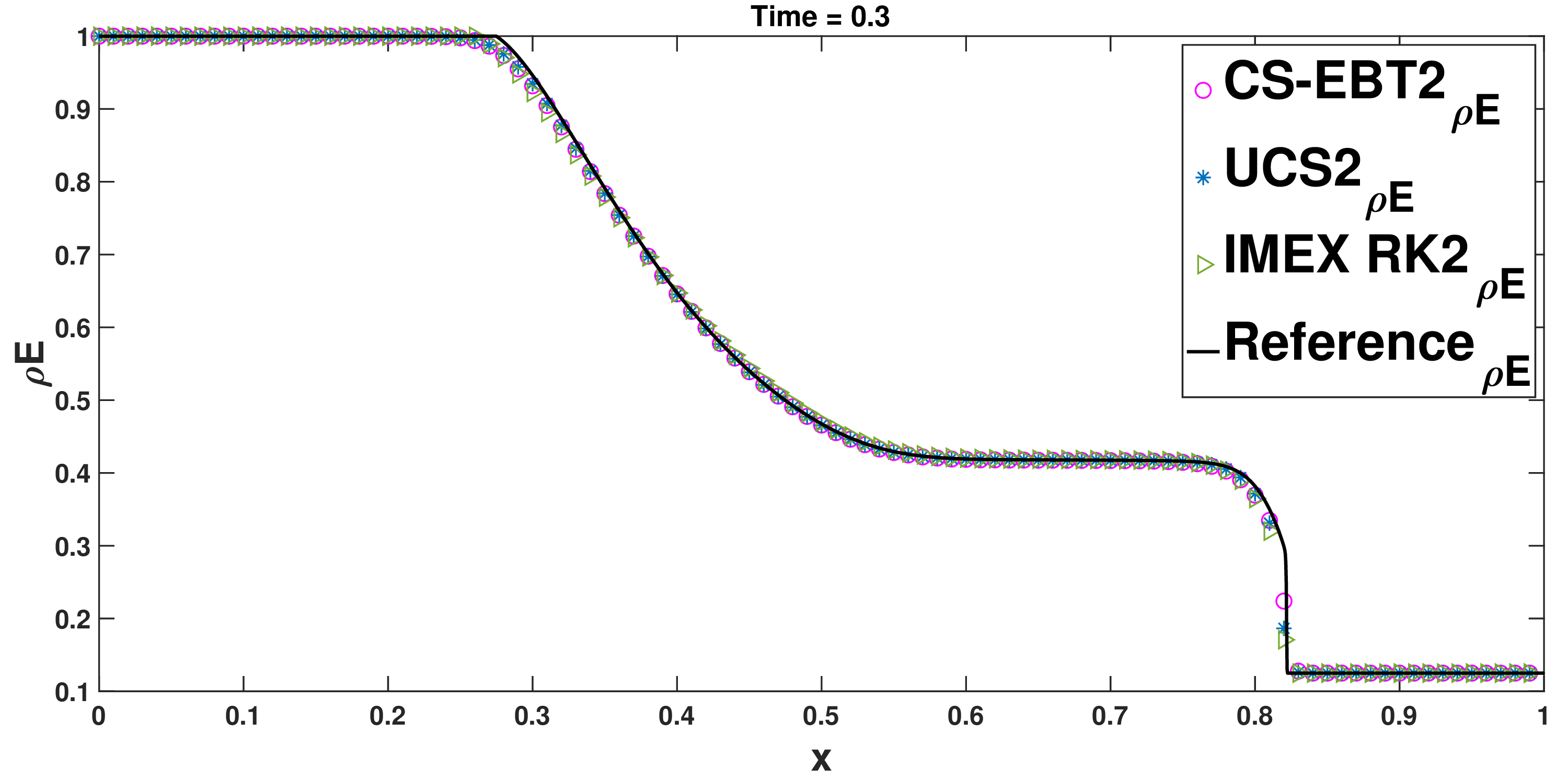}
    \caption{Alternative Euler equations with heat transfer model with non-smooth case: comparison between numerical solution $\rho$(left), $\rho u$(right) and $\rho E$(center) and the reference solution (IMEX RK2) with $\tau = 0.02$, CFL $0.9$ and $N = 200$.}
    \label{(euler4b)}
\end{figure}
\begin{figure}[!ht]
    \begin{minipage}[b]{0.48\linewidth}
        \includegraphics[width=\linewidth]{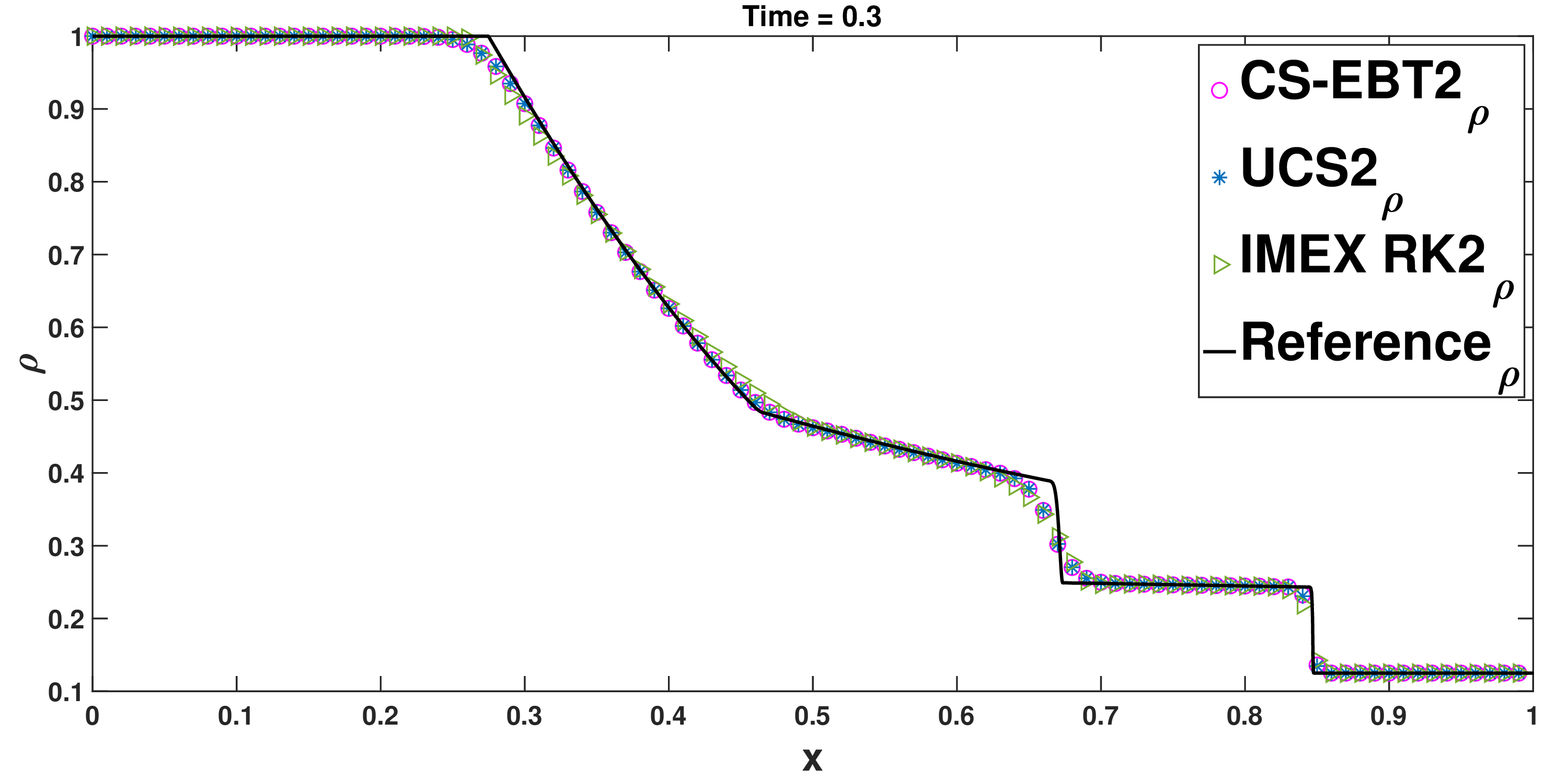}
    \end{minipage}
    \hfill
    \begin{minipage}[b]{0.48\linewidth}
        \includegraphics[width=\linewidth]{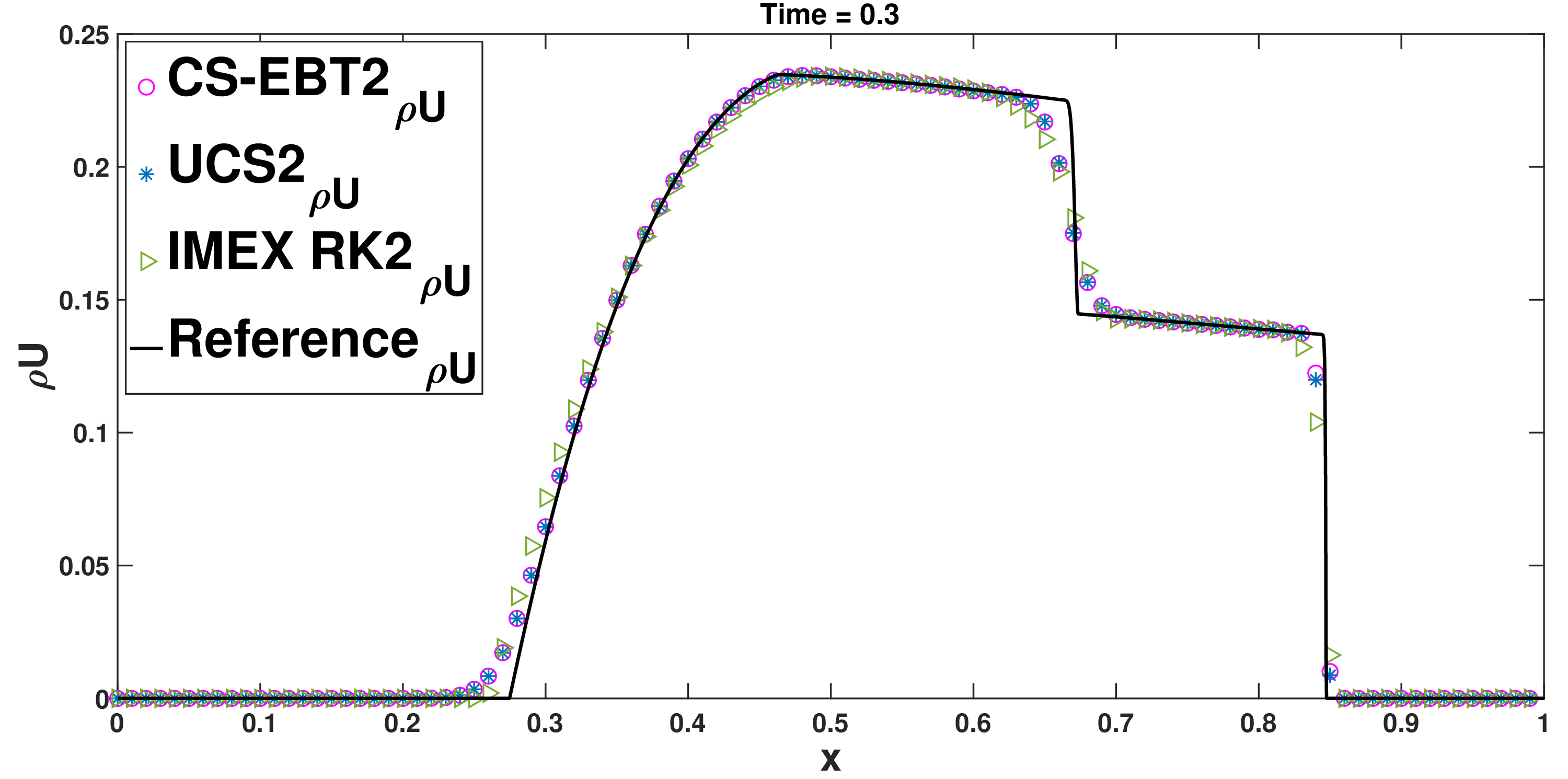}
    \end{minipage}
    \vspace{0.5cm} 
    \centering
    \includegraphics[width=0.5\linewidth]{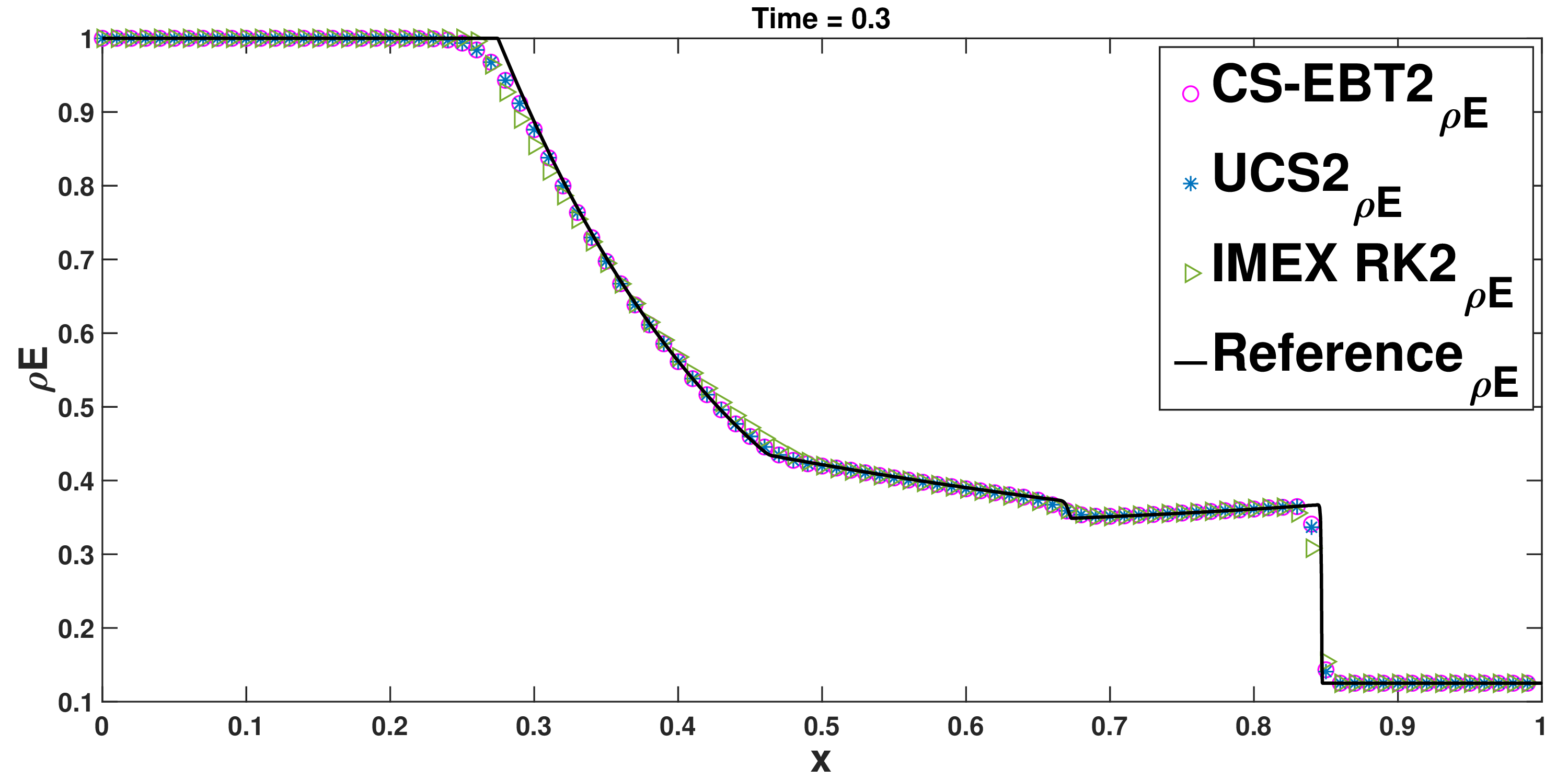}
    \caption{Alternative Euler equations with heat transfer model with non-smooth case: comparison between numerical solution $\rho$(left), $\rho u$(right) and $\rho E$(center) and the reference solution (IMEX RK2) with $\tau = 1$, CFL $0.9$ and $N = 200$.}
    \label{(euler4c)}
\end{figure}
Transmissive boundary conditions are applied at both ends of the computational domain. The governing equations correspond to the alternative Euler system with heat transfer. Numerical solutions are computed using the CS-EBT2, UCS2 and IMEX RK2 schemes on a uniform grid with $N = 200$ points and are compared against a high-resolution reference solution obtained with the IMEX RK2 scheme on a finer mesh with $N = 3200$ points. In Fig.~\ref{(4a)}, the relaxation parameter is set to $\tau = 10^{-8}$, representing a regime close to thermal equilibrium. Figure~\ref{(euler4b)} and \ref{(euler4c)} presents the results for a larger relaxation parameter, $\tau = 0.02$ and $\tau = 1$ respectively, with all other parameters unchanged. Table~\ref{Tabeuler} reports the numerical results for the integral invariant of this model, demonstrating that the CS-EBT2, UCS2 and IMEX RK2 schemes preserve the invariance.
\renewcommand{\arraystretch}{1.2}{
\begin{table}
 \caption{Integral invariance for alternative Euler equations with heat transfer model with initial condition \eqref{discontinuous_euler} for $\rho$ and $\rho u + \tau \rho E$ with $\tau = 10^{-8}$ and CFL $0.9$}
 \centering 
\begin{tabular}{*{7}{c}}
\toprule
\multirow{2}{*}{$N$}
& \multicolumn{2}{c}{CS-EBT2}
& \multicolumn{2}{c}{UCS2}
& \multicolumn{2}{c}{IMEX RK2} \\
\cmidrule(lr){2-3} 
\cmidrule(lr){4-5}
\cmidrule(lr){6-7}
& Error $\rho$ & Error $\rho u + \tau \rho E$ & Error $\rho$ & Error $\rho u + \tau \rho E$ &  Error $\rho$ & Error $\rho u + \tau \rho E$ \\
\midrule

20 & 1.11022302e-16 & 1.05493549e-01 & 1.11022302e-16 & 1.04999968e-01 & 1.20032300e-16 & 1.03064892e-08\\

40 & 1.11022302e-16 & 1.05266833e-01 & 1.11022302e-16 & 1.04999927e-01 & 1.20032300e-16 & 1.03064892e-16\\

80 & 1.11022302e-16 & 1.05138651e-01 & 1.11022302e-16 & 1.0499987e-01 & 1.20032300e-16 & 1.03064892e-16\\

160 & 1.11022302e-16 & 1.05070756e-01 & 1.11022302e-16 & 1.04999777e-01 & 1.20032300e-16 & 1.03064892e-16\\

320 & 1.11022302e-16 & 1.05035791e-01 & 1.11022302e-16 & 1.04999619e-01 & 1.20032300e-16  & 1.03064892e-16\\

640 & 1.11022302e-16 & 1.05017963e-01 & 1.11022302e-16  &1.04999288e-01 & 1.20032300e-16  & 1.03064892e-16\\

1280 & 1.11022302e-16 & 1.05009004e-01 & 1.11022302e-16 & 1.04998607e-01 & 1.20032300e-16 & 1.03064892e-16\\

\bottomrule
\end{tabular}
\label{Tabeuler} 
\end{table}}
\subsubsection{2D Jin-Xin relaxation model}
We consider the alternative relaxation approximation for the 2D Jin--Xin relaxation model~\eqref{modify:Xin_Jin_2D_model}. We first consider a smooth, well-prepared initial condition to assess the accuracy and stability of the proposed numerical scheme.\\
\noindent\textbf{Smooth case (Well-Prepared):}
The smooth initial data is as follows
\begin{equation}
\label{smoothdata:XinJin2d}
\begin{cases}
u(x,y,0) = \sin(2\pi x)\sin(2\pi y), \\[2mm]
v(x,y,0) = a\,u(x,y,0), \\[2mm]
w(x,y,0) = b\,u(x,y,0),
\end{cases}
\end{equation}
where the parameters are fixed as $a = 0.7$ and $b = 0.5$. With this choice, the initial data lies on the equilibrium manifold, leading to a smooth evolution of the solution. This configuration provides an appropriate setting to assess the accuracy and stability of the proposed numerical scheme in the absence of discontinuities. The spatial domain $[0,1]\times[0,1]$ is discretized using a uniform Cartesian mesh with $N_x = N_y = 400$ grid points in the $x$- and $y$-directions. The numerical solution is computed up to the final time $T = 0.35$ with a CFL number of $0.9$. 
\begin{figure}[htbp!]
     \centering
     \begin{subfigure}[b]{0.32\textwidth}
         \centering  
         \includegraphics[scale=0.21, trim=12.5cm 0cm 11.5cm 0cm, clip]{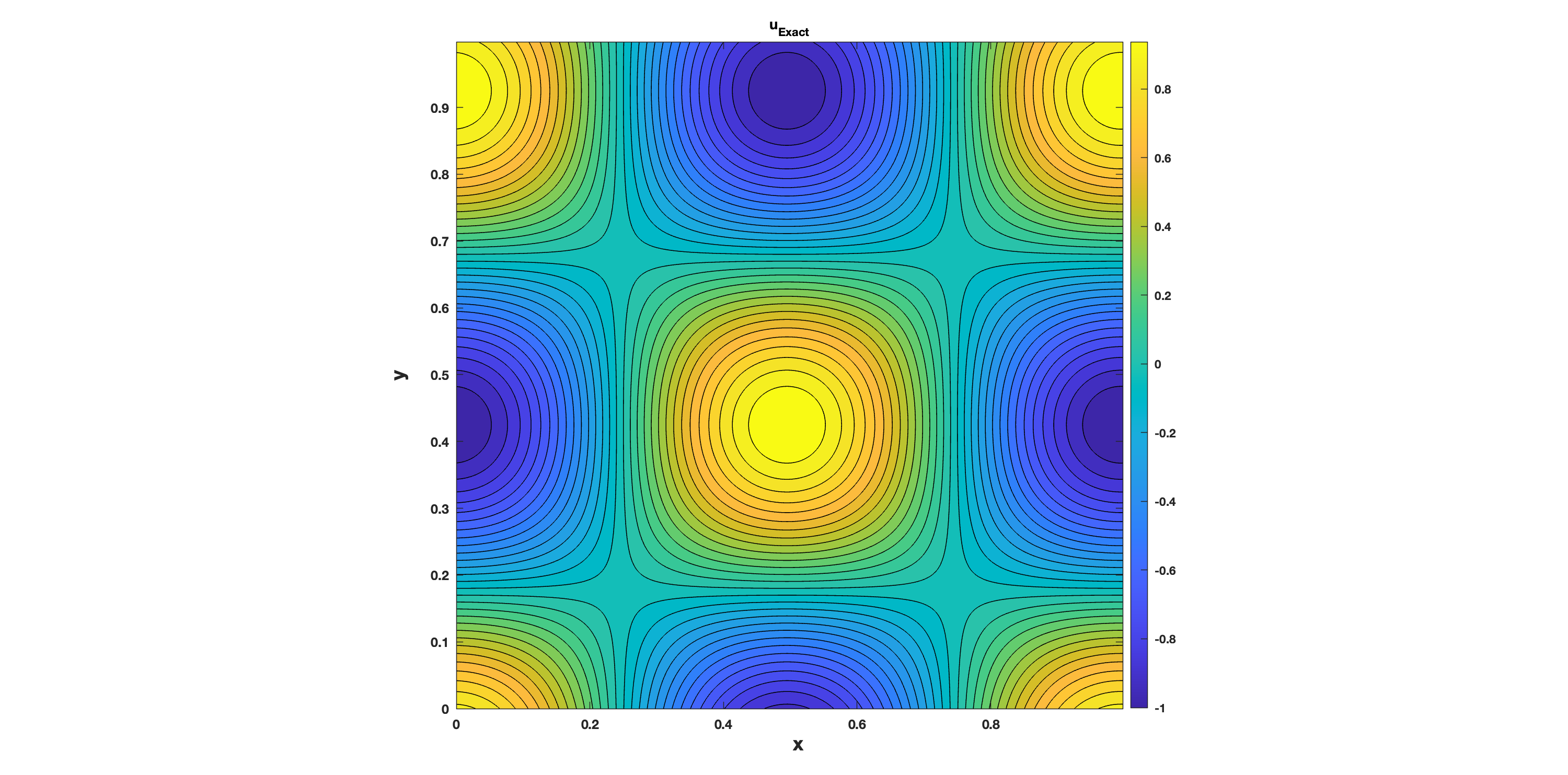}
      \caption{$u$ exact}
         \label{2dxinu:exact}
     \end{subfigure}
     \hfill
    \begin{subfigure}[b]{0.32\textwidth}
     \centering  
      \includegraphics[scale=0.21, trim=12.5cm 0cm 11.5cm 0cm, clip]{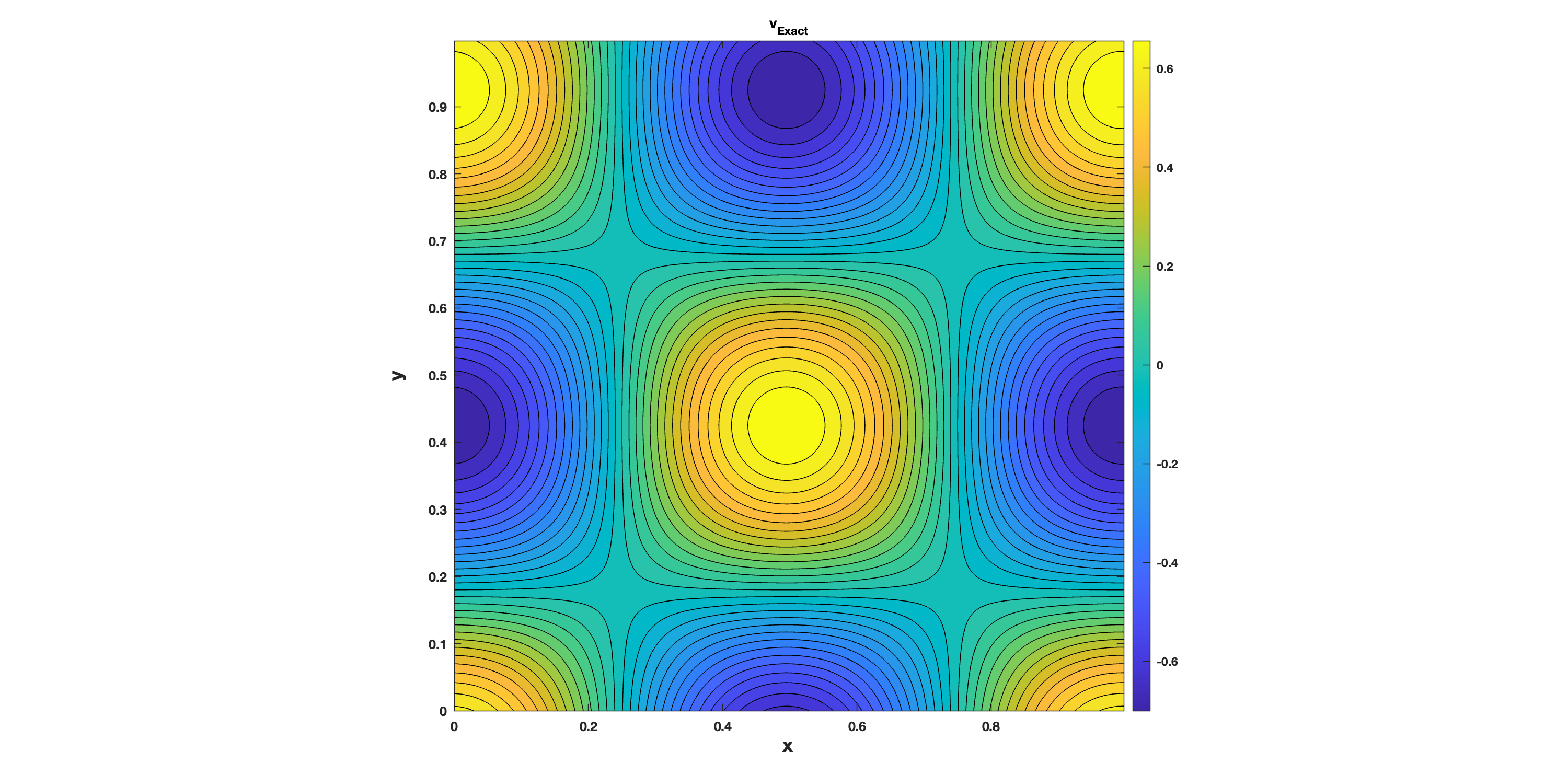}
      \caption{$v$ exact}
    \label{2dxinv:exact}
     \end{subfigure}
     \hfill
     \begin{subfigure}[b]{0.32\textwidth}
         \centering  
         \includegraphics[scale=0.21, trim=12.5cm 0cm 11.5cm 0cm, clip]{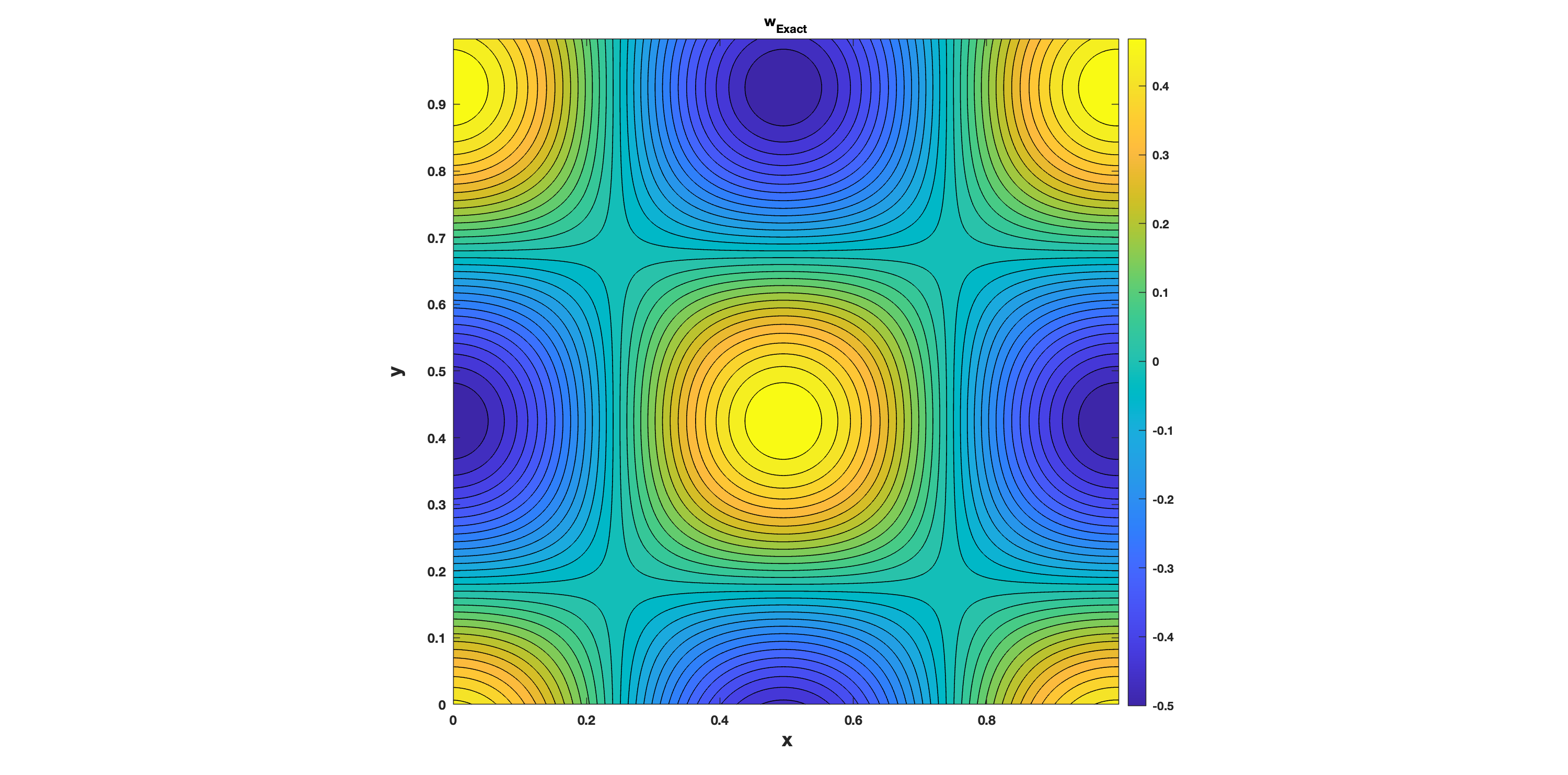}
      \caption{$w$ exact}
         \label{2dxinw:exact}
     \end{subfigure}
     \hfill
      \begin{subfigure}[b]{0.32\textwidth}
         \centering   
         \includegraphics[scale=0.215, trim=12.5cm 0cm 11.5cm 0cm, clip]{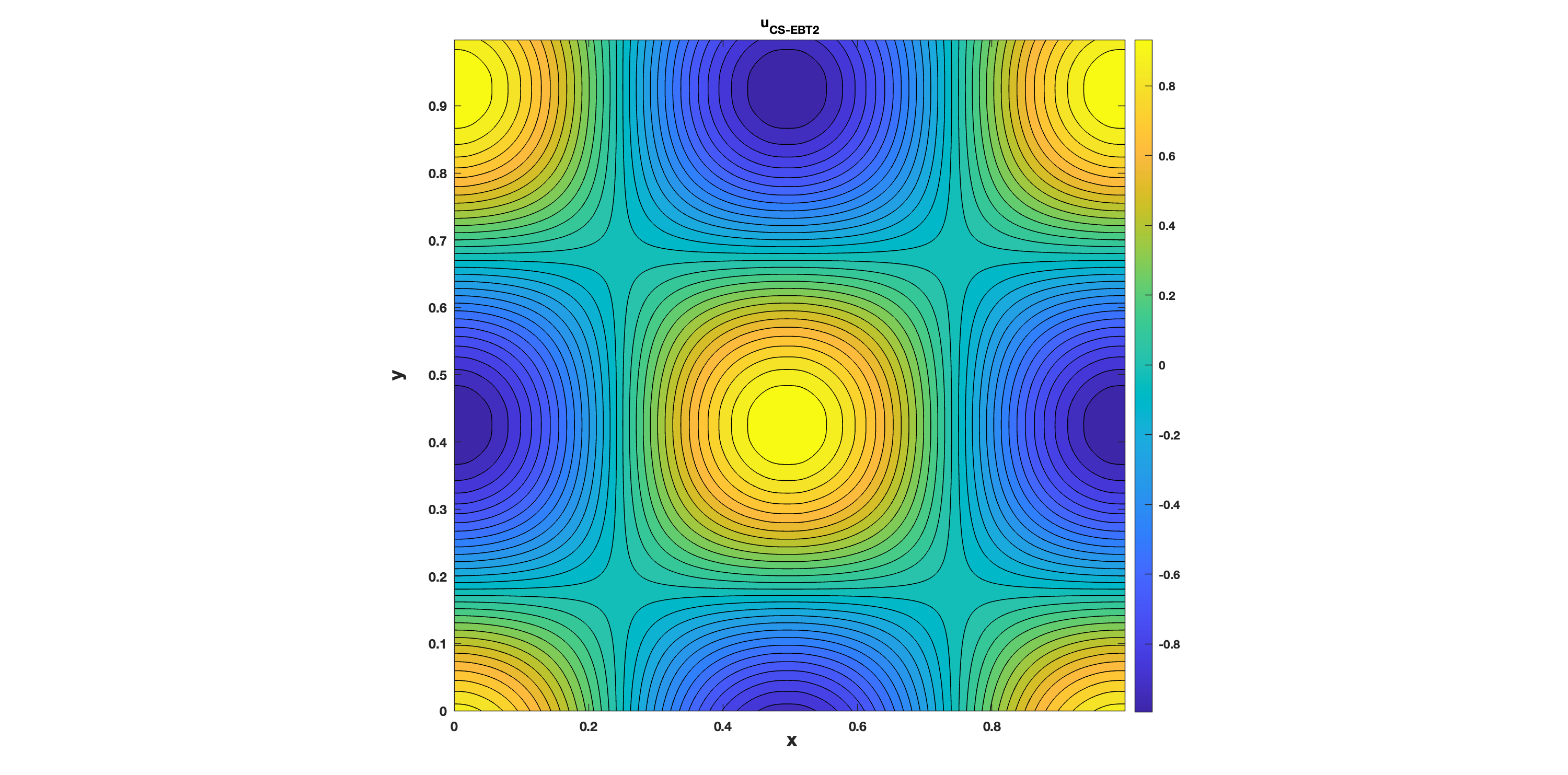}
         \caption{$u$ CS-EBT2}
         \label{2dxinu:csebt}
     \end{subfigure}
     \hfill
     \begin{subfigure}[b]{0.32\textwidth}
         \centering   
         \includegraphics[scale=0.215, trim=12.5cm 0cm 11.5cm 0cm, clip]{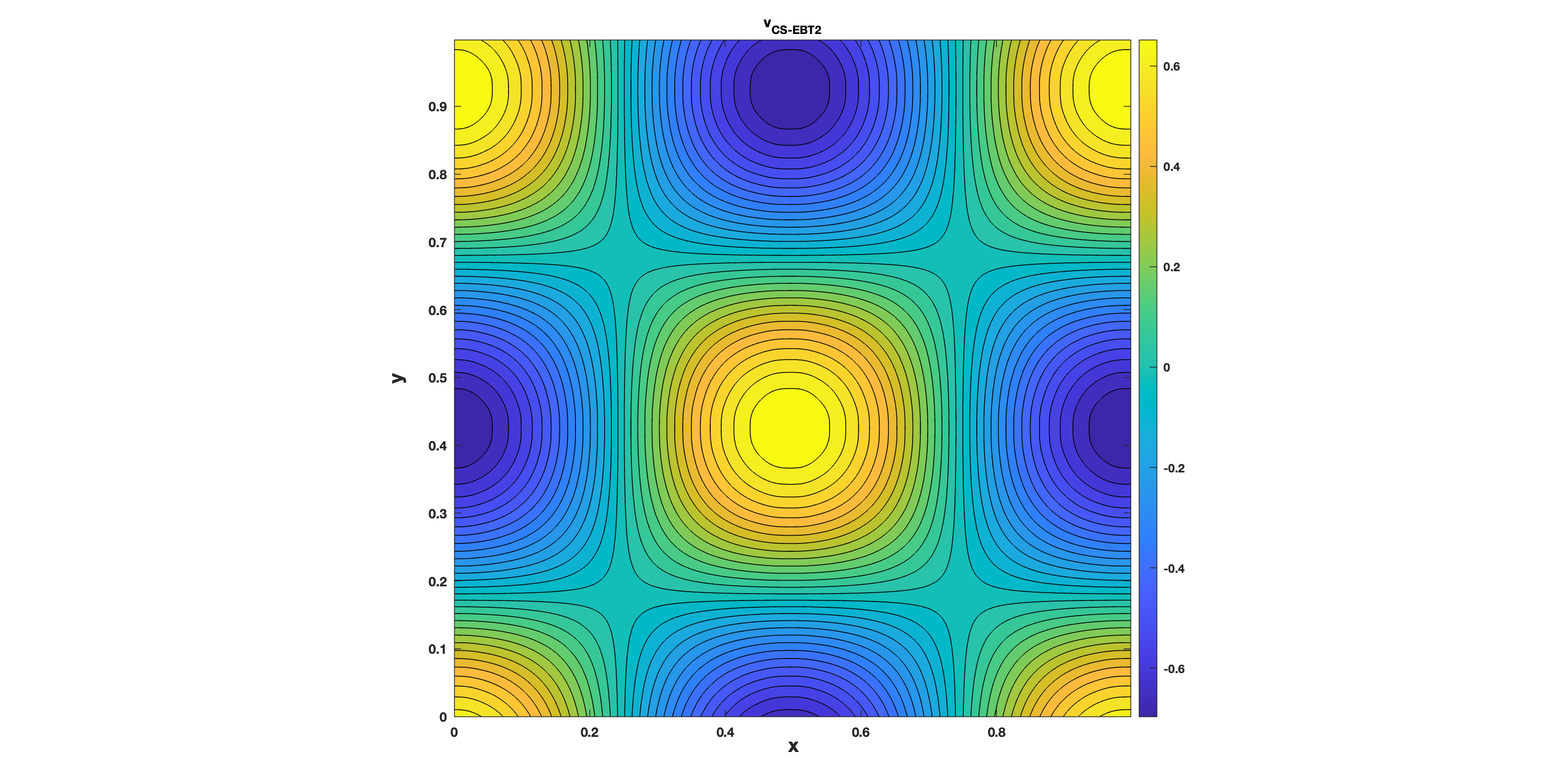}
         \caption{$v$ CS-EBT2}
         \label{2dxinv:csebt}
     \end{subfigure}
     \hfill
     \begin{subfigure}[b]{0.32\textwidth}
         \centering   
         \includegraphics[scale=0.215, trim=12.5cm 0cm 11.5cm 0cm, clip]{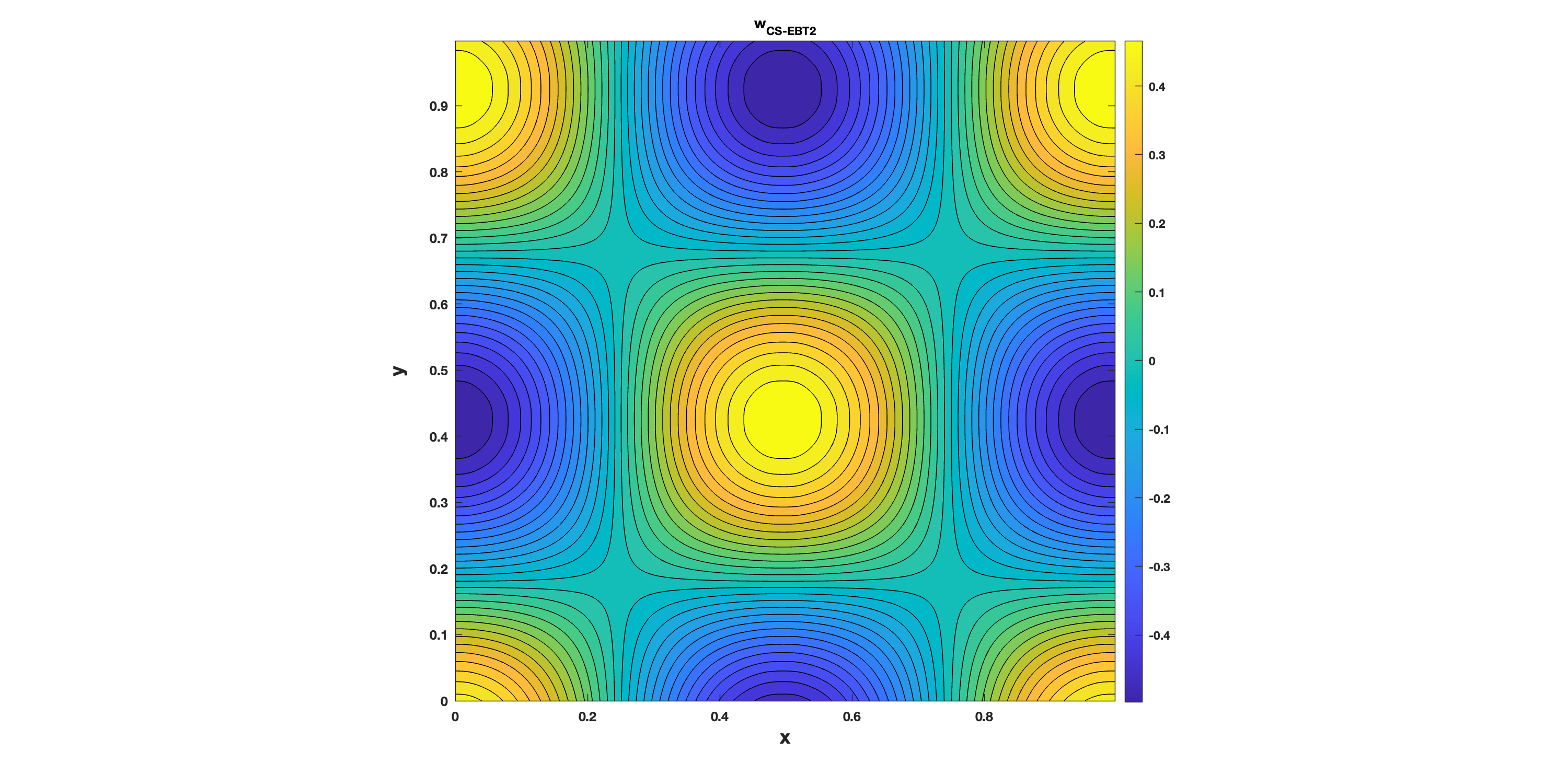}
         \caption{$w$ CS-EBT2}
         \label{2dxinw:csebt}
     \end{subfigure}
 \end{figure}
 \begin{figure} 
    \begin{subfigure}[b]{0.32\textwidth}
         \centering
      \includegraphics[scale=0.215, trim=12.5cm 0cm 11.5cm 0cm, clip]{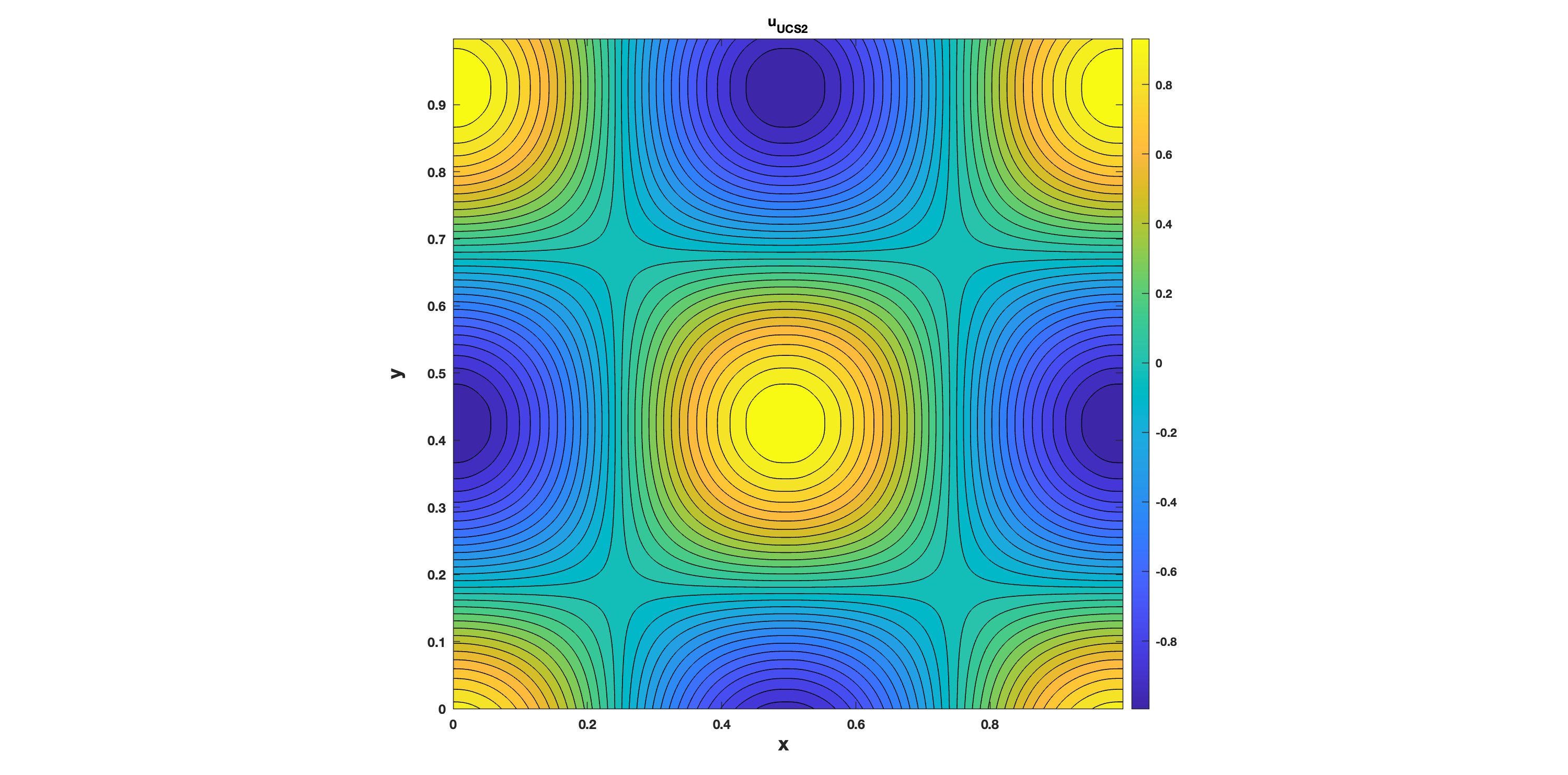}
     \caption{$u$ UCS2}
         \label{2dxinu:ucs}
     \end{subfigure}
     \hfill
    \begin{subfigure}[b]{0.32\textwidth}
         \centering
      \includegraphics[scale=0.215, trim=12.5cm 0cm 11.5cm 0cm, clip]{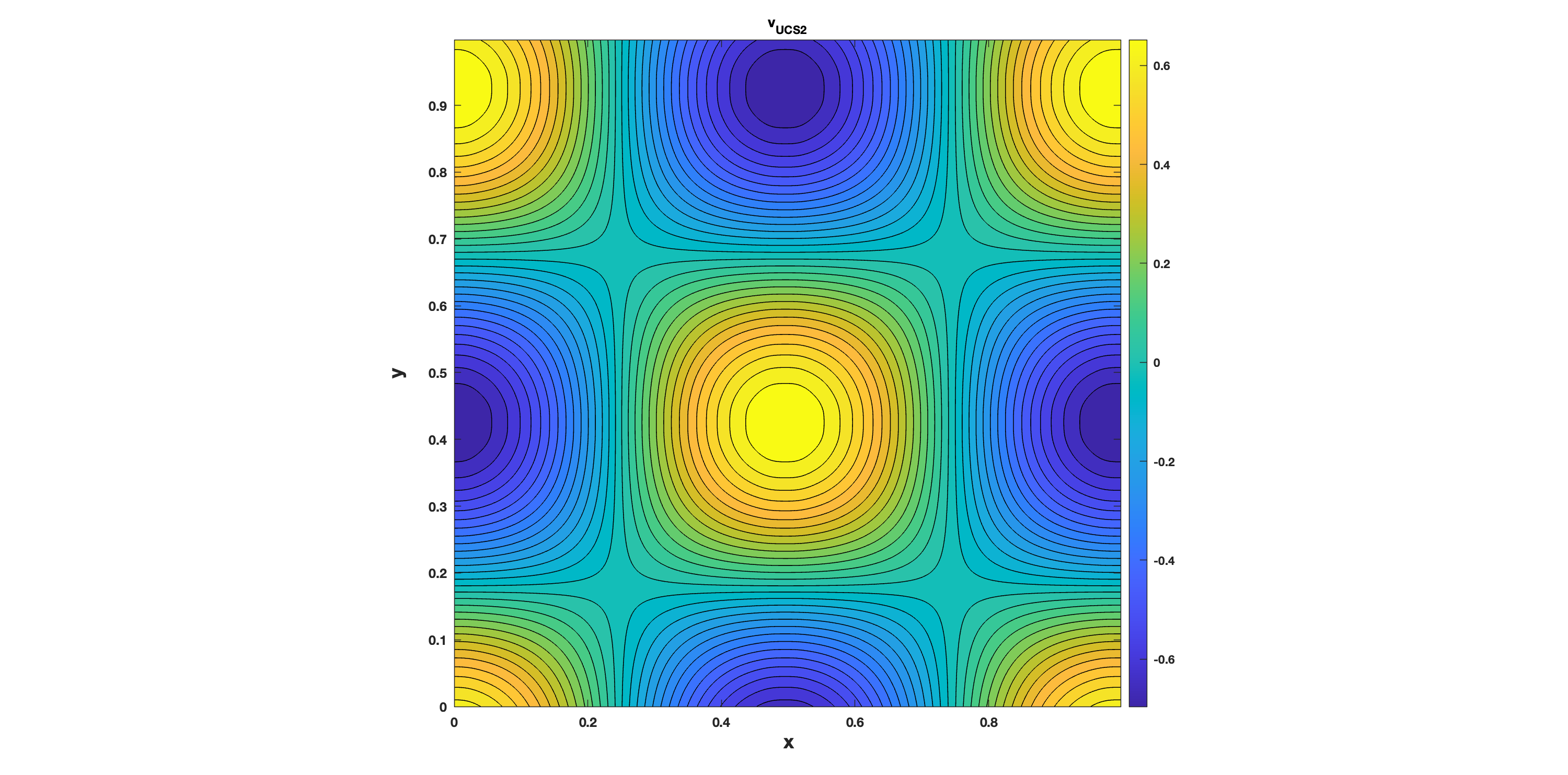}
     \caption{$v$ UCS2}
         \label{2dxinv:ucs}
     \end{subfigure}
     \hfill
    \begin{subfigure}[b]{0.32\textwidth}
         \centering
      \includegraphics[scale=0.215, trim=12.5cm 0cm 11.5cm 0cm, clip]{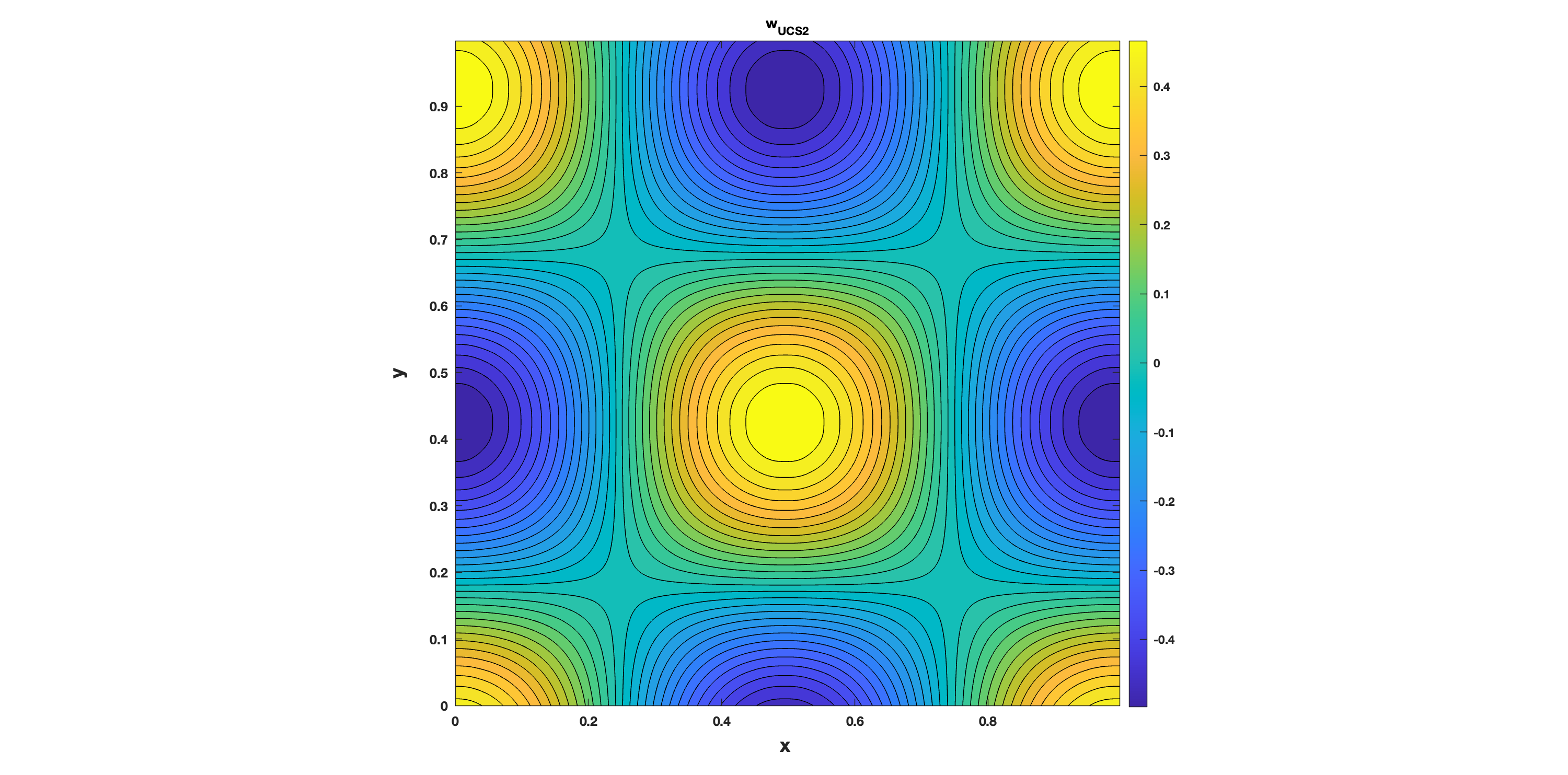}
     \caption{$w$ UCS2}
         \label{2dxinw:ucs}
     \end{subfigure}   
     \begin{subfigure}[b]{0.32\textwidth}
         \centering  
        \includegraphics[scale=0.215, trim=12.5cm 0cm 11.5cm 0cm, clip]{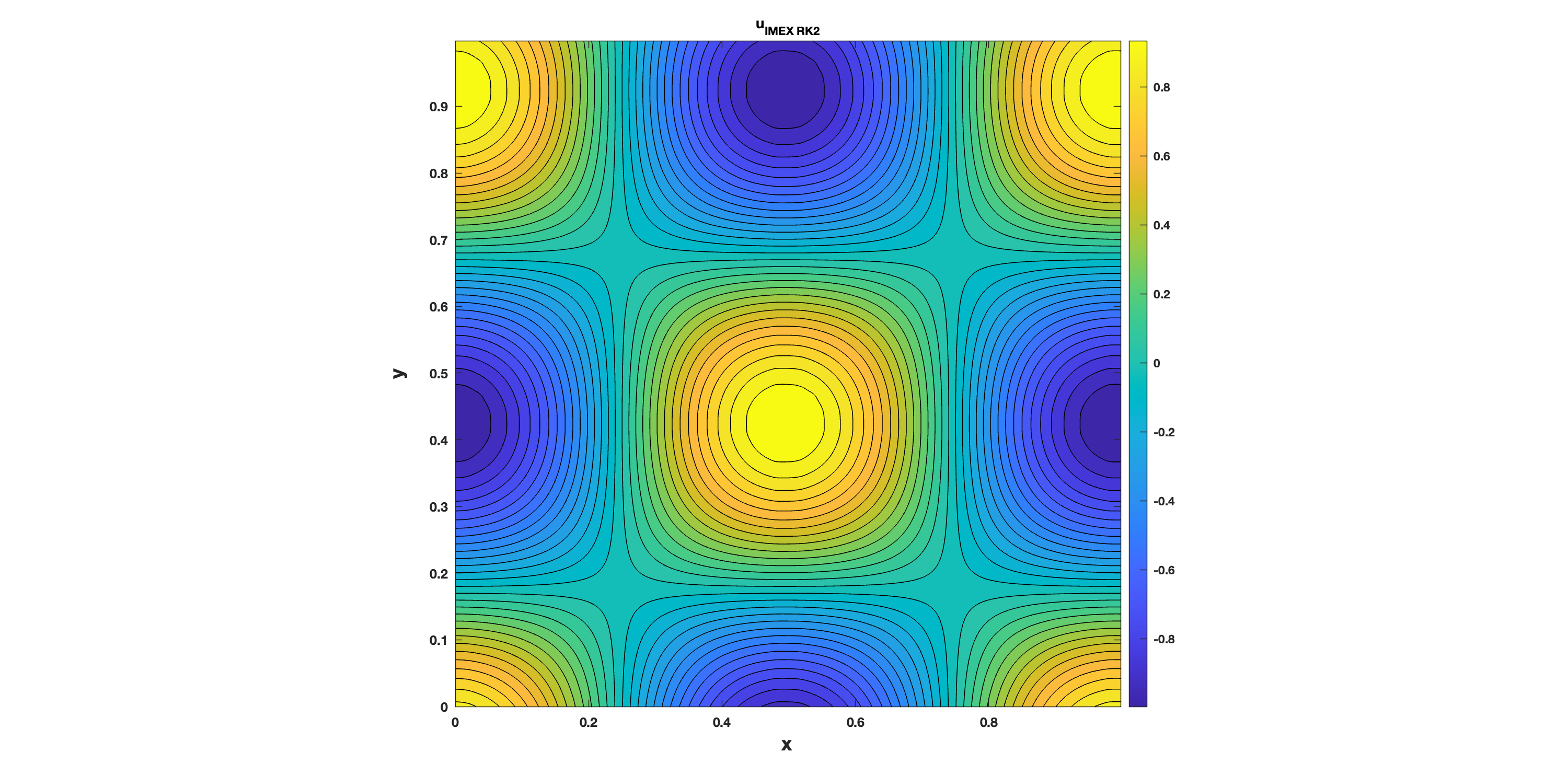}
      \caption{$u$ IMEX RK2}
         \label{2dxinu:rk2}
     \end{subfigure}
     \hfill
    \begin{subfigure}[b]{0.32\textwidth}
         \centering  
         \includegraphics[scale=0.215, trim=12.5cm 0cm 11.5cm 0cm, clip]{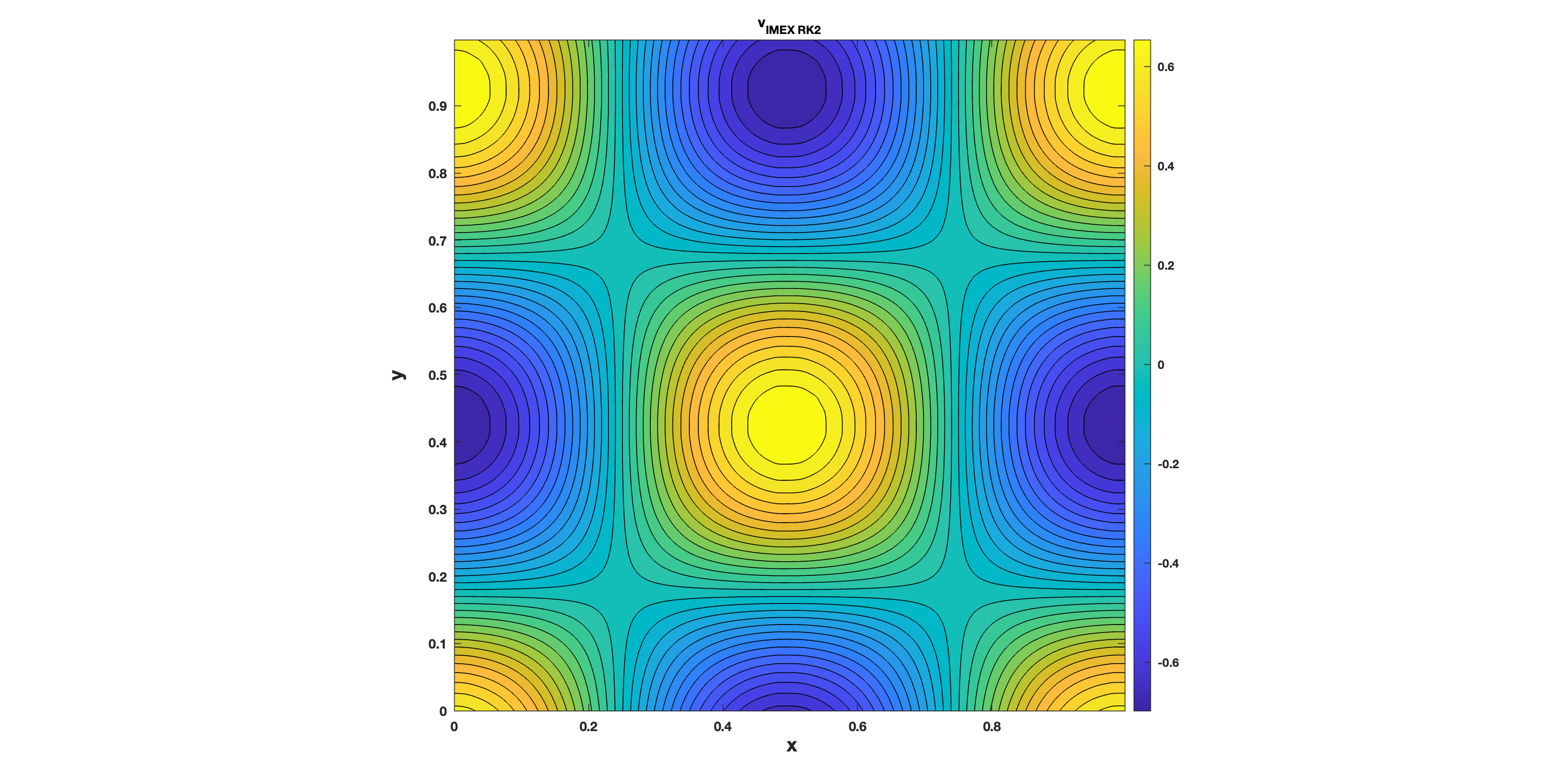}
      \caption{$v$ IMEX RK2}
         \label{2dxinv:rk2}
     \end{subfigure}
     \hfill
     \begin{subfigure}[b]{0.32\textwidth}
         \centering  
         \includegraphics[scale=0.215, trim=12.5cm 0cm 11.5cm 0cm, clip]{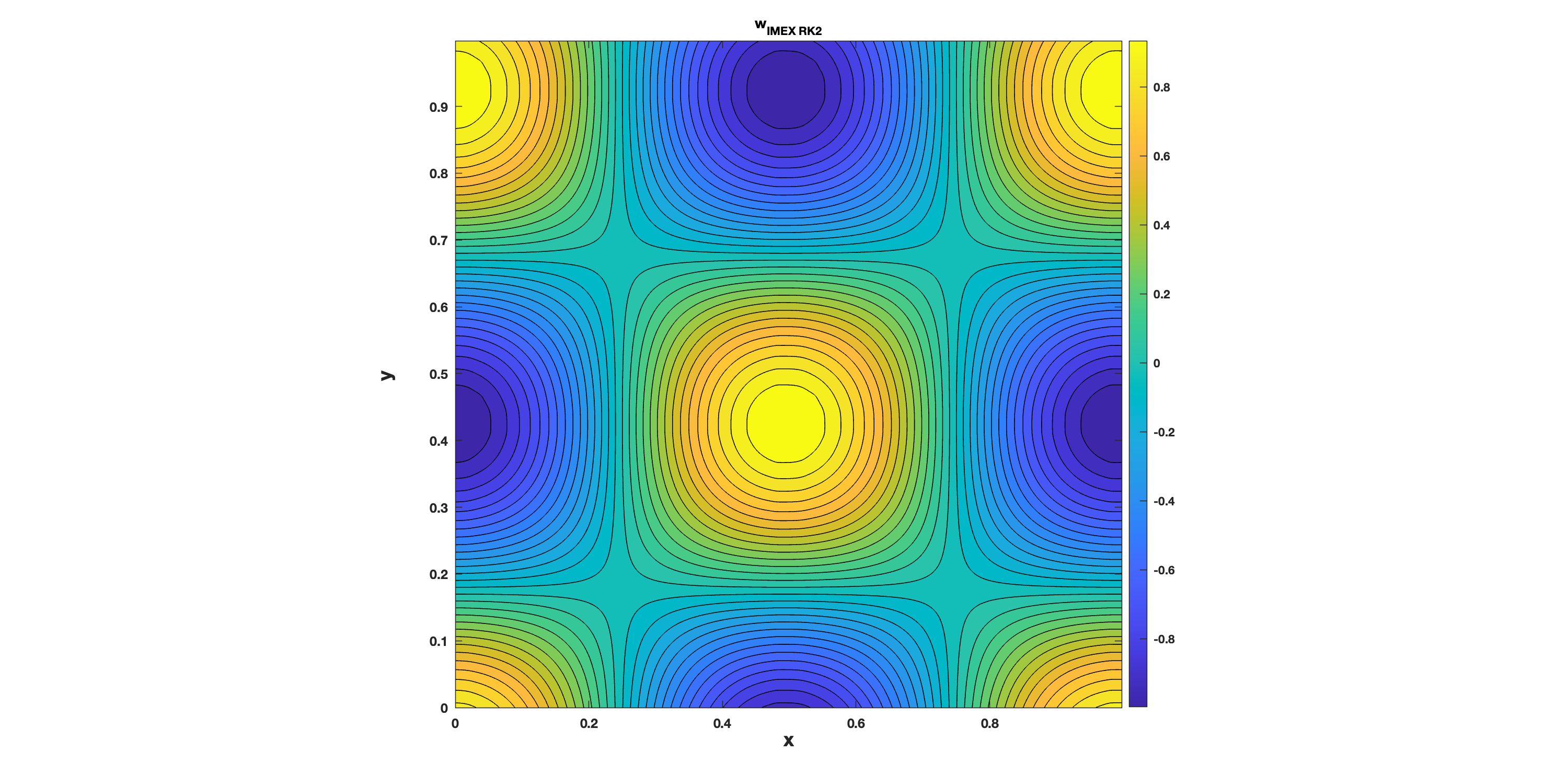}
      \caption{$w$ IMEX RK2}
         \label{2dxinw:rk2}
     \end{subfigure}
         \caption{Exact and numerical contour plots with $30$ levels of the variables $u$, $v$, and $w$ for the two-dimensional smooth test problem computed using the Jin-Xin relaxation system. The computations are carried out on a uniform $400 \times 400$ Cartesian grid with CFL number $0.9$, relaxation parameter $\tau = 10^{-10}$, and final time $T = 0.35$.}
    \label{2dxin:num}
\end{figure}
Figure~\ref{2dxin:num} compares the exact and numerical solutions of the variables $u$, $v$, and $w$ for $\tau = 10^{-10}$. The top row presents the exact solutions, whereas the bottom row shows the numerical results. In each row, the left, center, and right panels correspond to $u$, $v$, and $w$, respectively. All contour plots use 30 levels. The exact solutions are obtained analytically using the Fourier transform method. Tables~\ref{Tabu}, \ref{Tabv}, and \ref{Tabw} present the experimentally observed orders of convergence for the variables $u$, $v$, and $w$, respectively. The results were obtained using a CFL number of $0.9$ and a relaxation parameter of $\tau = 10^{-10}$, with the initial condition given by \eqref{smoothdata:XinJin2d}. Table~\ref{Tabxinjin2d} reports the numerical errors associated with the integral invariant for the two-dimensional linear test problem using the CSEBT2, UCS2, and IMEX RK2 schemes. The integral invariant is assessed by computing the $L^1$-error between the conserved quantity $u+\tau v+\tau w$ and its initial value $u_0+\tau v_0+\tau w_0$.
\renewcommand{\arraystretch}{1.2}{
\begin{table}[ht!]
\caption{$L^{1}$-Error and order of convergence for $u$ (2D Jin-Xin model)}
\centering
\begin{tabular}{*{7}{c}}
\toprule
\multirow{2}{*}{$N \times N$} & 
\multicolumn{2}{c}{CS-EBT2} &
\multicolumn{2}{c}{UCS2} & 
\multicolumn{2}{c}{IMEX RK2} \\
\cmidrule(lr){2-3}
\cmidrule(lr){4-5}
\cmidrule(lr){6-7}
 & $L^{1}$-Error & Order 
 & $L^{1}$-Error & Order 
 & $L^{1}$-Error & Order \\
\midrule
50$\times$ 50 & 5.9475e-04 &-  & 5.8071e-04  & - & 7.4750e-03 & - \\
100 $\times$ 100& 1.4215e-04 & 2.06 & 1.3866e-04  & 2.07 & 2.2717e-03 & 1.72\\
200 $\times$ 200& 3.3952e-05  & 2.07 & 3.2934e-05 & 2.07 &  6.1055e-04 & 1.90 \\
400 $\times$ 400& 8.3170e-06 & 2.03 & 7.9451e-06 & 2.05 & 1.4911e-04 & 2.03 \\
\bottomrule
\end{tabular}
\label{Tabu}
\end{table}
}
 \renewcommand{\arraystretch}{1.2}{
\begin{table}[ht!]
\caption{$L^{1}$-Error and order of convergence for $v$ (2D Jin-Xin model)}
\centering
\begin{tabular}{*{7}{c}}
\toprule
\multirow{2}{*}{$N \times N$} & 
\multicolumn{2}{c}{CS-EBT2} &
\multicolumn{2}{c}{UCS2} & 
\multicolumn{2}{c}{IMEX RK2} \\
\cmidrule(lr){2-3}
\cmidrule(lr){4-5}
\cmidrule(lr){6-7}
 & $L^{1}$-Error & Order 
 & $L^{1}$-Error & Order 
 & $L^{1}$-Error & Order \\
\midrule
50 $\times$ 50 & 4.1633e-04 & - & 4.0650e-04 & - & 5.2325e-03 & - \\
100 $\times$ 100 & 9.9505e-05 & 2.06 & 9.7059e-05 & 2.07 & 1.5902e-03 & 1.72 \\
200 $\times$ 200 & 2.3766e-05 & 2.07 & 2.3054e-05 & 2.07 & 4.2738e-04 & 1.90\\
400 $\times$ 400 & 5.8219e-06 & 2.03 & 5.5615e-06 & 2.05 & 1.0438e-04 & 2.03\\
\bottomrule
\end{tabular}
\label{Tabv}
\end{table}}
\renewcommand{\arraystretch}{1.2}{
\begin{table}[ht!]
\caption{$L^{1}$-Error and order of convergence for $w$ (2D Jin-Xin model)}
\centering
\begin{tabular}{*{7}{c}}
\toprule
\multirow{2}{*}{$N \times N$} & 
\multicolumn{2}{c}{CS-EBT2} &
\multicolumn{2}{c}{UCS2} & 
\multicolumn{2}{c}{IMEX RK2} \\
\cmidrule(lr){2-3}
\cmidrule(lr){4-5}
\cmidrule(lr){6-7}
 & $L^{1}$-Error & Order 
 & $L^{1}$-Error & Order 
 & $L^{1}$-Error & Order \\
\midrule
50 $\times$ 50  &  2.9738e-04 & - & 2.9036e-04 & - & 3.7375e-03 & - \\
100 $\times$ 100 & 7.1075e-05 & 2.06 & 6.9328e-05 & 2.07 & 1.1358e-03 & 1.72 \\
200 $\times$ 200 & 1.6976e-05 & 2.06 & 1.6467e-05 & 2.07 & 3.0527e-04 & 1.89\\
400 $\times$ 400 & 4.1585e-06 & 2.03 & 3.9725e-06 & 2.05 & 7.4557e-05  & 2.03 \\
\bottomrule
\end{tabular}
 \label{Tabw}
\end{table}}
As the mesh is refined, the $L^1$-error remains nearly unchanged and close to machine precision, indicating that the discrete integral invariant is preserved by all three numerical schemes. These numerical results provide strong evidence in support of the theoretical conservation property established for the two-dimensional linear system. The computations were carried out with a CFL number of $0.9$, relaxation parameter $\tau = 10^{-10}$, and parameters $a=0.7$ and $b=0.5$. For all mesh resolutions considered, the CSEBT2, UCS2, and IMEX RK2 schemes preserve the integral invariant with the accuracy.
 \renewcommand{\arraystretch}{1.2}{
\begin{table}[ht!]
\centering
\caption{Integral invariance for alternative Jin-Xin 2D model with initial condition
\eqref{smoothdata:XinJin2d} for $u + \tau v + \tau w$ with $\tau = 10^{-10}$ and CFL $0.9$}
\begin{tabular}{*{5}{c}}
\toprule
\multirow{2}{*}{$N \times N$}
& \multicolumn{1}{c}{CS-EBT2}
& \multicolumn{1}{c}{UCS2}
& \multicolumn{1}{c}{IMEX RK2} \\
\cmidrule(lr){2-2} 
\cmidrule(lr){3-3}
\cmidrule(lr){4-4}
& Error $u + \tau v + \tau w $ & Error $u + \tau v + \tau w $ & Error $u + \tau v + \tau w $ \\
\midrule

300 $\times$ 300& 9.68934763e-06 & 9.68701876e-06  & 2.84217094e-18  \\

400 $\times$ 400& 5.48073811e-06 & 5.47971942e-06  & 1.79890824e-18 \\

500 $\times$ 500& 3.51948454e-06 & 3.51896187e-06  & 2.89901436e-18 \\

600 $\times$ 600& 2.44950790e-06 & 2.44920070e-06 & 5.83422199e-18 \\

700 $\times$ 700& 1.80237124e-06 & 1.80217695e-06  & 3.36420234e-18 \\

800 $\times$ 800& 1.38154450e-06 & 1.38141335e-06 & 1.87428198e-18\\

900 $\times$ 900& 1.09253893e-06 & 1.09244631e-06 &  2.68427256e-18\\
\bottomrule
\end{tabular}
\label{Tabxinjin2d}
\end{table}}
\section{conclusion}\label{sec:conclusion}
In this paper, we have investigated Vasudeva Murthy’s relaxation approach, originally introduced for the Jin--Xin relaxation model, and examined its applicability to several benchmark hyperbolic balance law systems. In the one-dimensional setting, the relaxation framework has been considered for the shallow-water equations, the Broadwell model, and the Euler equations with heat transfer. Although the alternative Jin--Xin relaxation model is available in the literature, corresponding relaxation formulations for these three systems have been constructed and analyzed in the present work. Furthermore, we have extended the relaxation framework to the two-dimensional Jin--Xin model. For all the considered models, the associated integral invariants have been established analytically, including the invariant corresponding to the proposed two-dimensional extension.

To verify the theoretical results computationally, we have implemented second-order numerical schemes, namely the CS-EBT2 scheme, the UCS2 central scheme, and the IMEX-RK2 scheme. The CS-EBT2 and UCS2 schemes are both extensions of the NT central scheme to non-homogeneous balance laws, with different way of treatment of the source terms. In particular, the CS-EBT2 scheme employs an implicit treatment of the source term, while the UCS2 scheme involves predictor-stage evaluations together with the implicit treatment of the stiff contribution. The IMEX RK2 scheme is included as a well-established method for stiff source terms and is also used to generate reference solutions on refined meshes when exact solutions are unavailable.

The numerical experiments provide a comprehensive verification of the analytical results across different relaxation regimes. In particular, the computed solutions confirm the preservation of the derived integral invariants for both the one-dimensional models and the two-dimensional Jin-Xin model. The numerical results also demonstrate the expected second-order accuracy of the implemented schemes. The comparison between CS-EBT2 and UCS2 further indicates that CS-EBT2 provides a larger stability region and allows less restrictive CFL conditions, while requiring fewer predictor stages and lower computational cost. These properties make CS-EBT2 particularly attractive for stiff relaxation problems. Overall, the results provide both theoretical and numerical evidence for the invariant-preserving character of the considered relaxation formulations and demonstrate that the established numerical schemes can be effectively applied to these models, including the proposed two-dimensional extension.


\section*{Data Availability} The datasets generated during and/or analyzed during the current study are available from the corresponding author on reasonable request.
\section*{Declarations}
\textbf{Conflict of interest:} The authors declare that they have no conflict of interest.

\section*{Funding}
The work of both authors is supported by the project NBHM, DAE, India, Ref. No. 02011/46/2021 NBHM(R.P.)/R\&D II/14874.

\bibliographystyle{ieeetr}
\bibliography{biblio}
\end{document}